\documentclass[11pt,a4paper]{article}
\usepackage{amsmath,amssymb,amsthm,mathrsfs}
\usepackage{color}
\usepackage{ascmac}
\usepackage{dsfont}

\usepackage{enumitem}
\usepackage{tikz}
\usepackage{graphicx}
\theoremstyle{definition}
\newtheorem{dfn}{Definition}[section]
\newtheorem{remark}[dfn]{Remark}
\newtheorem{rmk}[dfn]{Remark}

\theoremstyle{plain}
\newtheorem{thm}[dfn]{Theorem}
\newtheorem{prop}[dfn]{Proposition}
\newtheorem{lem}[dfn]{Lemma}
\newtheorem{cor}[dfn]{Corollary}
\newtheorem{ass}[dfn]{Assumption}

\newcommand{\bn}{{\mathbf n}}

\newcommand{\bk}{{\mathbf k}}

\newcommand{\bu}{{\mathbf u}}

\newcommand{\bv}{{\mathbf v}}

\newcommand{\ba}{{\mathbf a}}

\newcommand{\bA}{{\mathbf A}}
\newcommand{\bB}{{\mathbf B}}

\newcommand{\bD}{{\mathbf D}}

\newcommand{\bF}{{\mathbf F}}
\newcommand{\bG}{{\mathbf G}}
\newcommand{\bH}{{\mathbf H}}
\newcommand{\bI}{{\mathbf I}}

\newcommand{\bL}{{\mathbf L}}
\newcommand{\bM}{{\mathbf M}}

\newcommand{\bV}{{\mathbf V}}

\newcommand{\bS}{{\mathbf S}}
\newcommand{\bT}{{\mathbf T}}

\newcommand{\DV}{{\rm Div}\,}
\newcommand{\dv}{{\rm div}\,}

\newcommand{\BR}{{\mathbb R}}
\newcommand{\BC}{{\mathbb C}}

\newcommand{\BN}{{\mathbb N}}

\newcommand{\CA}{{\mathcal A}}
\newcommand{\CB}{{\mathcal B}}

\newcommand{\CD}{{\mathcal D}}

\newcommand{\CF}{{\mathcal F}}

\newcommand{\CK}{{\mathcal K}}
\newcommand{\CL}{{\mathcal L}}

\newcommand{\CN}{{\mathcal N}}

\newcommand{\CR}{{\mathcal R}}
\newcommand{\CS}{{\mathcal S}}
\newcommand{\CT}{{\mathcal T}}

\newcommand{\CP}{{\mathcal P}}

\newcommand{\CW}{{\mathcal W}}
\newcommand{\CX}{{\mathcal X}}
\newcommand{\CY}{{\mathcal Y}}
\newcommand{\CZ}{{\mathcal Z}}

\newcommand{\fh}{{\mathfrak h}}
\newcommand{\fg}{{\mathfrak g}}
\newcommand{\fq}{{\mathfrak q}}
\newcommand{\fp}{{\mathfrak p}}

\newcommand{\fs}{{\mathfrak s}}

\newcommand{\fL}{{\mathfrak L}}

\newcommand{\bg}{{\mathbf g}}
\newcommand{\bh}{{\mathbf h}}
\newcommand{\pd}{\partial}

\usepackage{xcolor}

\newcommand{\ep}{\varepsilon}
\newcommand{\pa}{\partial}
\newcommand{\wt}{\widetilde}
\newcommand{\wh}{\widehat}

\newcommand{\Bf}{{\mathbf f}}
\newcommand{\Bg}{{\mathbf g}}

\newcommand{\Bv}{{\mathbf v}}

\newcommand{\Di}{{\rm Div}\,}
\newcommand{\di}{{\rm div}\,}
\newcommand{\loc}{{\rm loc\,}}
\newcommand{\Hol}{{\rm Hol\,}}

\newcommand{\dOm}{\dot{\Omega}}

\numberwithin{equation}{section} 

\usepackage{scalerel,stackengine}
\stackMath
\newcommand\reallywidehat[1]{%
	\savestack{\tmpbox}{\stretchto{%
			\scaleto{%
				\scalerel*[\widthof{\ensuremath{#1}}]{\kern-.6pt\bigwedge\kern-.6pt}%
				{\rule[-\textheight/2]{1ex}{\textheight}}
			}{\textheight}%
		}{0.5ex}}%
	\stackon[1pt]{#1}{\tmpbox}%
}

\newcommand{\vertiii}[1]{{\left\vert\kern-0.25ex\left\vert\kern-0.25ex\left\vert #1 
		\right\vert\kern-0.25ex\right\vert\kern-0.25ex\right\vert}}
\allowdisplaybreaks[3]

\begin{document}
	
\title{\bf On the solution operators arising from the gas-liquid two-phase problem in unbounded domains with finite depth}
\author{Miao Tu 
\thanks{School of Mathematical Sciences,
Tongji University, 
No.1239, Siping Road, Shanghai (200092), China.
 \endgraf
e-mail address: tm2230925@163.com }
and Xin Zhang
\thanks{School of Mathematical Sciences,
Key Laboratory of Intelligent Computing and Applications (Ministry of Education), 
Tongji University, 
No.1239, Siping Road, Shanghai (200092), China. \endgraf
e-mail address: xinzhang2020@tongji.edu.cn
}}

\date{}
\maketitle

\begin{abstract}
This paper studies some evolution equations arising from the sharp interface problem of the compressible-incompressible Navier-Stokes equations in unbounded domains in $\mathbb{R}^N (N\geq2)$, where the viscous gases initially occupy the upper half space and the viscous liquids below initially lie in the strip-like domain. 
In order to establish the maximal $L_p$-$L_q$ regularity estimates of the evolution problem, 
we construct the $\mathcal{R}$-solver of the resolvent problem associated to the gas-liquid two-phase operator.
The crucial part of our proof lies in the analysis of the \emph{explicit} resolvent operators defined in the unbounded domains with flat boundaries.

\vskip1pc\noindent
subjclass[2020]: Primary: 35Q30; Secondary: 35R35 \vskip0.5pc\noindent
keywords: Maximal $L_p$-$L_q$ regularity;  $\mathcal{R}$-boundedness; Unbounded domains; Free boundary conditions; \\
\phantom{keywords } Compressible-incompressible Navier-Stokes equations.
\end{abstract}

\tableofcontents

\section{Introduction}
\subsection{Model problem}

For any positive constant $b$ and $x'=(x_1,\dots,x_{N-1}) \in \BR^{N-1}$ ($N\geq 2$), let us set
\begin{equation}\label{dfn:Omega}
\begin{aligned}
  	\Omega_+ &=\{x=(x',x_N)\in\BR^N:  x_N>0\}, \\
	\Omega_- &=\{x=(x',x_N)\in \BR^N: -b<x_N<0\},\\
 	\Gamma & =\{x=(x',x_N)\in\BR^N: x_N= 0\},\\
  	S & =\{x=(x',x_N)\in\BR^N: x_N= -b\}.
\end{aligned}
\end{equation}
Let $\BR_+=(0,\infty)$ and $\bn=(0,\cdots,0,1)^\top \in \BR^N.$
This paper concerns the following model problem arising from the study of the sharp interface problem of the gas-liquid two-phase flows in unbounded domains (see  Figure \ref{Fig:Omega2}):
 \begin{equation}\label{eql:1}
	\left\{ \begin{aligned}
		&\pd_t\rho_+ + \gamma_{1+}\, \dv \bv_+=f_+
		&&\quad&\text{in} &\quad \Omega_+ \times \BR_+, \\
		&\gamma_{1+}\pd_t\bv_+ - \DV\bT_+(\bv_+,\gamma_{2+}\rho_+)  = \bg_+
		&&\quad&\text{in}& \quad \Omega_+ \times \BR_+, \\
		&\gamma_{1-}\pd_t\bv_-  - \Di\bT_-(\bv_{-},\fp_-) = \bg_-
		&&\quad&\text{in}& \quad \Omega_- \times \BR_+, \\
		& \di \bv_- = g_d=\di \fg_d
		&&\quad&\text{in}& \quad \Omega_- \times \BR_+, \\
		&\big(\bT_+(\bv_+,\gamma_{2+}\rho_+)- \bT_-(\bv_-,\fp_-)\big)\bn
		= \bh
		&&\quad&\text{on}& \quad \Gamma \times \BR_+,\\
		&\bv_{+}=\bv_{-}
		&&\quad&\text{on}& \quad \Gamma \times \BR_+,  \\
		& \bv_-=0
		&&\quad&\text{on}& \quad   S \times \BR_+, \\
		&(\rho_+, \bv_+)|_{t=0} = (\rho_{+,0}, \bv_{+,0})
		&&\quad&\text{in}& \quad \Omega_+, \\
		&\bv_-\big|_{t=0}=\bv_{-,0}
		&&\quad&\text{in}& \quad \Omega_-.
	\end{aligned}
	\right.
\end{equation}

In \eqref{eql:1}, the constant parameters $\gamma_{1\pm},$ $\gamma_{2+},$ $\mu_{\pm},\nu_+>0,$ and the stress tensors in \eqref{eql:1} are given by
\begin{align*}
\bT_{\pm}(\bv_{\pm},\pi_{\pm}) &=\bS_{\pm}(\bv_{\pm})  -\pi_{\pm}\bI,\\
	\bS_+(\bv_+) &= \mu_+ \bD(\bv_+)  + (\nu_+ - \mu_+)\di\bv_+\bI ,\\
	\bS_-(\bv_-) &= \mu_- \bD(\bv_-),\quad  \bD(\bv_\pm) = \nabla\bv_\pm + (\nabla\bv_\pm)^\top.
\end{align*}
In addition, $\dv \bv = \sum_{j=1}^N \pd_{x_j}v_j$ for $\pd_{x_j}=\pd/\pd x_j,$
the $(i,j)$th entry of the matrix $\nabla \bu$ is $\pd_{x_i} u_j,$ 
$\bI$ is the $N\times N$ identity matrix, 
$\bM^\top$ is the transposed matrix of $\bM=[M_{ij}],$
$\DV \bM$ denotes an $N$-vector of functions 
whose $i$-th component is $\sum_{j=1}^N \pd_{x_j} M_{ij}.$

Furthermore, the system \eqref{eql:1} can be regarded as the linearized model of the free boundary value problem of the compressible-incompressible Navier-Stokes equations in unbounded domains in view of the \emph{Lagrangian} coordinates (see Appendix \ref{app:non}).

\begin{center}
	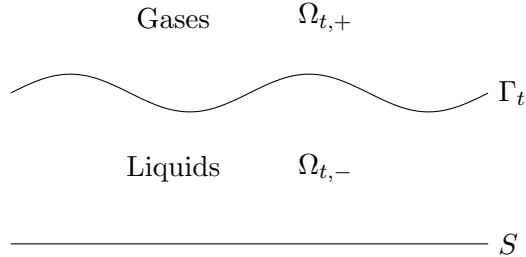
\begin{figure}[h]
		\centering
		\begin{tikzpicture}[>=stealth]
			
			\draw[ domain=-pi:pi, samples=200] plot (\x, {sin(2*\x r)/4}) node[right] {$\Gamma_t$};
			
			\draw (-pi,-2) -- (pi,-2) node[right] {$S$};
			
			\node at (1,1) {$\Omega_{t,+}$};
			\node at (1,-1) {$\Omega_{t,-}$};
			\node at (-1,1) {Gases};
			\node at (-1,-1) {Liquids};
			
		\end{tikzpicture}
		\caption{Domain of the gas-liquid two-phase free flows}
		\label{Fig:Omega2}
	\end{figure}
\end{center}

\subsection{Brief review of the sharp interface problem}
The sharp interface problem is one of fundamental models describing the interfacial phenomena.
For the motion of the viscous fluids, there are a lot of literature on the mathematical theory of one-phase or two-phase flows. To our knowledge, the first breakthrough on the mathematical theory of the sharp interface problem of viscous fluids is due to Solonnikov \cite{Sol1977}
where the author used Lagrangian coordinates and proved the local well-posedness of the one-phase incompressible Navier-Stokes equations in bounded domains within H\"older continuity (see \cite{Sol1987} for the framework of Sobolev spaces).

Another well-known problem on the one-phase incompressible flows is so-called the \emph{viscous surface wave} problem which could trace back to the works \cite{Beale1981, Beale1983} by Beale. The viscous surface wave problem has been revisited by Guo and Tice in a series of works \cite{GT2013b,GT2013a,GT2013c} where the authors developed new high-order energy method (also see the recent results \cite{Gui2021,Wang2020} in this direction).

Based on the $\mathcal{R}$-boundedness theory in Weis' work \cite{Weis2001}, Shibata and Shimizu \cite{shibata2007,ShiShi2008} established a more general approach to handling the solutions of the one-phase incompressible Navier-Stokes in the class of the maximal $L_p$-$L_q$ regularity provided with the indices $p,q\in (1,\infty).$ Very recently, the maximal $L_1$ regularity for the free boundary value problem has been investigated in \cite{DHMT2020,OS2024,SW2024}.

All the aforementioned methods on the incompressible viscous flows can be applied to study the one-phase compressible Navier-Stokes equations with free boundary value conditions. However, due to the more complicated structure of the compressible Navier-Stokes equations, it is not immediate to solve the global well-posedness issue just in view of the results on the incompressible flows. For the interface problem of the compressible Navier-Stokes equations, we refer to the works \cite{DT2023, HL2021,SZ2023,Zaja1993,Zaja1994} and the references therein.
\medskip

For the mathematical study of viscous two-phase problem, the incompressible-incompressible Navier-Stokes problem has attracts significant attention. Assuming the domains are smooth and bounded, the first attempt is due to Tanaka \cite{Tanaka1993} on the global solvability of the model describing the motion of two immiscible inhomogeneous viscous liquids in $L_2$ type Sobolev–Slobodetskii spaces with the capillarity taken into account (also see \cite{Denisova2014}). The extension to the $L_p$ framework, we refer to more recent works on the bounded domains \cite{Pruss2011,PruS2016,SSZ2020}. Furthermore, the progress in the case of unbounded domains is much more involved. In this direction, we refer to the recent works \cite{Saito2024,WTK2014}. Besides, we refer to the works \cite{Deni1997,JTW2016a,KSS2016,Tani1984} for the gas-gas two-phase problem.
\medskip

Compared with the liquid-liquid two-phase problem, there are only a few of works on the gas-liquid two-phase problem. 
Denisova \cite{Denisova14} firstly obtained the local solvability in Sobolev-Slobodetskii space whenever two-phase flows occupy the bounded and exterior domains respectively. Then the more general $L_p$ theory has been developed by Shibata et al. in \cite{Kubo2021, KSS2014,Shibata17} where the key step is to study the resolvent problem of the linearized gas-liquid two-phase model in $\BR^N_+ \cup \BR^N_-$ for $N\geq 2$. 
In this paper, we study the problem \eqref{eql:1} in the unbounded reference domain $\Omega_{+} \cup \Omega_-$ with the fixed bottom $S,$ which reflects the feature of the compressible flow and the viscous surface wave problem.   

\subsection{Main result}

The main goal of this paper is to find the solution $(\rho_{+}, \bv_{\pm},\fp_{-})$ of \eqref{eql:1} in the class of the so-called \emph{maximal $L_p$-$L_q$ regularities} provided with suitable initial data $(\rho_{+,0}, \bv_{\pm,0})$ and source terms $(f_+, \bg_+, \bg_-, g_d,\fg_d,\bh).$ To this end, we introduce the following functional spaces
\footnote{See Subsection \ref{subsec:fs} for the definitions of the standard functional spaces and the operator $\Lambda^{1/2}_{\gamma}$ in this paper.}.
\begin{dfn} \label{def:FD}
Let $1<p,q<\infty$ and $\gamma>0.$
For $\Omega_\pm$ in \eqref{dfn:Omega}, we denote $\dOm=\Omega_+ \cup \Omega_-.$ 
\begin{itemize}

\item We write 
$(f_+, \bg_+, \bg_-, g_d,\fg_d,\bh) \in \CF_{p,q,\gamma}(\dot \Omega)$
if $g_d=\di \fg_d$ and the following regularity assumptions hold 
\begin{gather*}
f_+ \in L_{p,\gamma}\big(\BR; H^1_q(\Omega_{+})\big), \,\,\,
\Bg_\pm \in L_{p,\gamma}\big(\BR; L_q(\Omega_{\pm})^N\big), \\
g_d\in  L_{p,\gamma}\big(\BR; H^1_q(\Omega_-)\big) 
\cap H^{1\slash 2}_{p,\gamma}\big( \BR; L_q(\Omega_-)\big), \,\,\,
\fg_d \in H^1_{p,\gamma}\big(\BR; L_q(\Omega_-)^N\big),\\ 
\bh \in  L_{p,\gamma}\big(\BR; H^1_q(\dOm)^N\big) 
\cap H^{1\slash 2}_{p,\gamma}\big( \BR; L_q(\dOm)^N\big).
\end{gather*}
Moreover, we set 
\begin{equation*}
\begin{aligned}
&\|(f_+, \bg_+, \bg_-, g_d,\fg_d,\bh) \|_{\CF_{p,q,\gamma}(\dot \Omega)}\\
=& \| e^{-\gamma t} f_+\|_{L_p(\BR; H^1_q(\Omega_+))}   
+\| e^{-\gamma t}\bg_+ \|_{L_p(\BR; L_q(\Omega_{+}))}
+\| e^{-\gamma t}\bg_- \|_{L_p(\BR; L_q(\Omega_{-}))}\\
&+\| e^{-\gamma t} g_d\|_{L_p(\BR; H^1_q(\Omega_-))} 
+\| e^{-\gamma t} \Lambda_{\gamma}^{1\slash 2}  g_d\|_{L_p(\BR; L_q(\Omega_-))} 
+ \| e^{-\gamma t} \pd_t \fg_d\|_{L_p(\BR; L_q(\Omega_-))}\\
 &+\| e^{-\gamma t} \bh\|_{L_p(\BR; H^1_q(\dOm))}
+\| e^{-\gamma t} \Lambda_{\gamma}^{1\slash 2} \bh\|_{L_p(\BR; L_q(\dOm))}.
\end{aligned}
\end{equation*}

\item For any given $(f_+, \bg_+, \bg_-, g_d,\fg_d,\bh)$ in $\CF_{p,q,\gamma}(\dot \Omega),$
we say $(\rho_{+,0}, \bv_{+,0},\bv_{-,0}) \in \CD_{p,q}(\dOm)$ if 
\begin{equation*}
    \rho_{+,0} \in H^1_q(\Omega_+), \,\,\,\bv_{+,0} \in B^{2-2/p}_{q,p}(\Omega_+)^N 
    \,\,\,\text{and}\,\,\,\bv_{-,0} \in B^{2-2/p}_{q,p}(\Omega_-)^N
\end{equation*}
satisfying the following compatibility conditions:
\begin{equation}\label{eq:1.4}
\begin{aligned}
& \dv \bv_{-,0}=g_d|_{t=0}\,\,\,\text{in}\,\,\,\Omega_-,
\quad (\bv_{-,0},\nabla \varphi)_{\Omega_-} =(\fg_{d}|_{t=0},\nabla \varphi)_{\Omega_-},
\,\,\,\forall \,\,\varphi \in  H^1_{q',\Gamma}(\Omega_-),\\
&\CT_{\bn} \Big( \big(\mu_+ \bD(\bv_{+,0})-\mu_- \bD(\bv_{-,0}) \big)\bn \Big) \Big|_{\Gamma} =\CT_{\bn} (\bh|_{t=0})|_{\Gamma},\quad
(\bv_{+,0}-\bv_{-,0})\big|_{\Gamma}=0,\quad
\bv_{-,0}|_{S}=0,  
\end{aligned}
\end{equation}
where 
\begin{equation}\label{def:proj}
    \CT_{\bn} \fh = \fh-(\fh \cdot \bn) \bn
\end{equation}
denotes the projection into the tangent space orthogonal to the normal vector $\bn.$ In addition, we also write 
\begin{equation*}
\|(\rho_{+,0}, \bv_{+,0},\bv_{-,0})\|_{\CD_{p,q}(\dOm)} 
= \|\rho_{+,0}\|_{H^1_q(\Omega_+)} 
+\sum_{\fs \in \{+,-\}}\|\bv_{\fs}\|_{B^{2-2/p}_{q,p}(\Omega_\fs)}.
\end{equation*}
\end{itemize}
\end{dfn}

With the notations in Definition \ref{def:FD}, we state the following main result of this paper.
\begin{thm}[Maximal $L_p$-$L_q$ regularity]\label{thm:mr}
Let $N\geq 2$ and $1<p,q<\infty$ with $2/p+N/q\not=1$ or $2.$ Then there exist positive constants $\gamma_0$ and $C$ such that the following assertions hold true.
\begin{itemize}
\item (Existence) 
Let $(f_+, \bg_+, \bg_-, g_d,\fg_d,\bh) \in \CF_{p,q,\gamma_0}(\dot \Omega)$
and $(\rho_{+,0}, \bv_{+,0},\bv_{-,0}) \in \CD_{p,q}(\dOm).$ Then there exists a solution $(\rho_+, \bv_{+}, \bv_{-}, \fp_-)$ of the system \eqref{eql:1} satisfying 
\begin{equation}\label{ss:thm1}
\begin{aligned}
\rho_+ \in H^1_{p,\gamma_0} \big(\BR_+; H^1_q(\Omega_+)\big),\quad
\fp_- \in L_{p,\gamma_0} \big(\BR_+; \wh H^1_q(\Omega_{-})\big),\\ 
\bv_\pm \in L_{p,\gamma_0} \big(\BR_+; H^2_q(\Omega_{\pm})^N\big)
\cap H^1_{p,\gamma_0} \big(\BR_+; L_q(\Omega_\pm)^N\big).
\end{aligned}
\end{equation}
Moreover, there holds
\begin{equation*}
\CN_{p,q,\gamma_0}(\rho_+,\bv_{+}, \bv_{-}, \fp_-) \leq 
C\big(  \|(\rho_{+,0}, \bv_{+,0},\bv_{-,0})\|_{\CD_{p,q}(\dOm)}
+\|(f_+, \bg_+, \bg_-, g_d,\fg_d,\bh) \|_{\CF_{p,q,\gamma_0}(\dot \Omega)} \big)
\end{equation*}
with 
\begin{equation*}
\begin{aligned}
\CN_{p,q,\gamma_0}(\rho_+,\bv_{+}, \bv_{-}, \fp_-)
=& \big\|e^{-\gamma_0 t} (\pd_t \rho_+, \gamma_0 \rho_+)\big\|_{L_p(\BR_+; H^1_q(\Omega_+))} 
+\big\|e^{-\gamma_0 t} \nabla \fp_-\big\|_{L_p(\BR_+; L_q(\Omega_-))}\\
&+\sum_{\fs \in \{+,-\}}
\big\|e^{-\gamma_0 t} (\pd_t \bv_{\fs}, \gamma_0 \bv_{\fs}, \Lambda_{\gamma_0}^{1/2} \nabla \bv_{\fs},\nabla^2 \bv_{\fs})\big\|_{L_p(\BR_+; L_q(\Omega_\fs))}.
\end{aligned}
\end{equation*}

\item (Uniqueness) Suppose that the solution $(\rho_+, \bv_{+}, \bv_{-}, \fp_-)$ of the system \eqref{eql:1} satisfies \eqref{ss:thm1} for  
\begin{equation*}
(\rho_{+,0}, \bv_{+,0},\bv_{-,0})=0
\quad \text{and}\quad
(f_+, \bg_+, \bg_-, g_d,\fg_d,\bh) =0.
\end{equation*}
Then we have $(\rho_+, \bv_{+}, \bv_{-},\fp_-)=0.$
\end{itemize}
\end{thm}

\begin{rmk}
Let us give some comments on Theorem \ref{thm:mr}.
\begin{itemize}
\item In this paper, we mainly focus on the case where the mass density of the gases has the \emph{constant} equilibrium state $\gamma_{1+}$ and the bottom $S$ is \emph{flat}. The advantage of such reduction is to derive the resolvent estimates in very precise manner (see Section \ref{sec:sf}), which sheds light on the challenging global-in-time theory of the gas-liquid two-phase problem. 
On the other hand, Theorem \ref{thm:mr} can help us find the local-in-time solutions of the nonlinear problem \eqref{eq:NS3} with the mass density of gases near $\gamma_{1+}$ by using similar arguments in one-phase problems (see \cite{EvBS2014,Kubo2021,shibata2007} for instance) and thus we omit it here.

\item  Based on the Weis' theory \cite{Weis2001}, the key step to prove Theorem \ref{thm:mr} is to establish the $\CR$-solver of the 
 resolvent problem of \eqref{eql:1}: 
\begin{equation}\label{eq:res_full}
	\left\{ \begin{aligned}
		&\lambda\rho_++\gamma_{1+}\di \bv_+=f_+
		&&\quad&\text{in}& \quad \Omega_+, \\
		&\gamma_{1+}\lambda\bv_+-\Di\bT_+(\bv_+,\gamma_{2+}\rho_+)=\bg_+
		&&\quad&\text{in}& \quad  \Omega_+, \\
		&\gamma_{1-}\lambda\bv_--\Di \bT_-(\bv_-,\fp_-)=\bg_-
		&&\quad&\text{in}& \quad  \Omega_-, \\ 
  		& \di\bv_-=g_d=\di\fg_d
		&&\quad&\text{in}& \quad  \Omega_-, \\ 
		&\bT_+(\bv_+,\gamma_{2+}\rho_+)\bn
		-\bT_-(\bv_-,\fp_-)\bn =\bh
		&&\quad&\text{on}& \quad \Gamma,\\
		& \bv_+-\bv_-=\bk
		&&\quad&\text{on}& \quad \Gamma,	\\
		&\bv_-= 0
		&&\quad&\text{on}& \quad S.   
	\end{aligned}
	\right.
\end{equation}
Here $\lambda$ is some complex number. For the result on the system \eqref{eq:res_full}, we refer to Theorem \ref{thm:main2} in Section \ref{sec:resolvent} for more details.
\end{itemize}
\end{rmk}

To conclude this subsection, we outline the structure of this paper. In Section \ref{sec:sf}, we consider some reduced resolvent problem (see \eqref{eq:2.3}) which eliminates the unknown $\rho_+$ from \eqref{eq:res_full}. By applying the partial Fourier transform with respect to the horizontal variable $x'$ to the reduced resolvent problem, we can explicitly solve the resulting ordinary differential equations of the vertical variable $x_N.$ Moreover, we shall analyze the Lopatinski determinant $\det \bL$ appearing in the solution formula of the reduced resolvent problem in Section \ref{sec:det} which is one of most difficult parts in this paper. In particular, from the derivation of the lower bound of  $\det \bL$ in Subsection \ref{subsec:lemma_3.1}, one may find the effect from both the compressible flows and the viscous surface wave problem. 
Using the result on $\det \bL$ in Section \ref{sec:det}, we construct the $\CR$-solver of the problem \eqref{eq:2.3} in Section \ref{sec:SO} and then return to the problem \eqref{eq:res_full} in Section \ref{sec:resolvent}. 
Let us emphasize that the standard technical results on the $\CR$-boundedness theory in the previous one-phase problems are not sufficient in our analysis of the two-phase problems. 
In this regard, we give some new technical results concerning the multiplier operators in Section \ref{sec:SO} but their proofs are postponed in Appendix \ref{proof}. At last, we apply the results on the resolvent problem \eqref{eq:res_full} to obtain the maximal $L_p$-$L_q$ regularity estimates of the system \eqref{eql:1} in Section \ref{sec:mr}.

\subsection{Notations}
\label{subsec:fs}

Let $\BC$ and $\BN$ denote the set of complex numbers and  natural numbers respectively. For any domain $D\subset\BR^N$, we write $L_q(D)$, $H_q^k(D)$ and $B^s_{q,p}(D)$ for the standard Lebesgue spaces, Sobolev spaces and Besov spaces with $k\in\BN$, $s\in \BR$ and $1\leq p,q\leq \infty$ respectively. In addition, let us set 
\begin{equation*}
   \begin{aligned}
   \wh H^1_q(D)&=\{\theta\in L_{q,\loc}(D):
\nabla\theta\in L_q(D)^N\},  \\
H^1_{q,\Gamma}(\Omega_-)&=\{\theta\in H^1_q(\Omega_-):\,\theta|_\Gamma=0\}
   \end{aligned}
\end{equation*}
with $\Omega_-$ defined in \eqref{dfn:Omega}. 
Moreover, for $\gamma>0,$ $\ell\in\BN$ and some Banach space $X,$ we write 
\begin{align*}
	L_{p,\gamma}(\BR; X)&=\{\theta\in L_{p,\loc}(\BR;X):
    e^{-\gamma t}\theta\in L_p(\BR; X)\},\\
	H^\ell_{p,\gamma}(\BR; X)&=\{\theta\in H^\ell_{p,\loc}(\BR;X):
    e^{-\gamma t}\pd_t^k\theta\in L_p(\BR; X),\,k=0,\dots,\ell\}.
\end{align*}

For $\lambda=\gamma+i\tau\in\BC$, we denote the Laplace transform and its inverse by
\begin{equation}\label{eq:Laplace}
	\begin{aligned}
		\CF_{L} [f] (\lambda)&=\int_{\BR} e^{-\lambda t} f(t) \,dt=\CF[e^{-\gamma t}f](\tau), \\
		\CF_{L}^{-1} [g] (t) &=\frac{1}{2\pi}\int_{\BR}  e^{\lambda t} g(\tau)\,d\tau =e^{\gamma t} \CF^{-1}[g] (t),
	\end{aligned}
\end{equation}
where $\CF$ denote the standard Fourier transform in $\BR.$ 
Then, for $\gamma>0$, let us set
\begin{equation}\label{eq:half}
\Lambda_{\gamma}^{1/2}f(t,x)
=\CF_{L}^{-1}\big[\lambda^{1/2}\CF_{L}[f(t,x)](\tau)\big](t)
=e^{\gamma t} \CF\big[\lambda^{1/2} \CF[e^{-\gamma t}f]\big].    
\end{equation}
The Bessel potential space of order $1/2$ is given by 
\begin{align*}
	H^{1/2}_{p,\gamma}(\BR;L_q(\Omega))&=\big\{f\in L_{p,\gamma}(\BR;L_q(\Omega)):e^{-\gamma t}\Lambda^{1/2}_\gamma f\in L_p(\BR;L_q(\Omega))\big\}
\end{align*}
with its norm 
$$\|f\|_{H^{1/2}_{p,\gamma}(\BR;L_q(\Omega))}=\|e^{-\gamma t}\Lambda_{\gamma}^{1/2}f\|_{L_p(\BR;L_q(\Omega))}.$$

For Banach spaces $X$ and $Y,$ $\CL(X; Y)$ denotes the set of bounded linear transformations from $X$ to $Y.$ 
Besides, $\Hol(\Lambda;X)$ denotes the set of $X$-valued functions defined on domain $\Lambda\subset \BC$.
 
In addition, $A \lesssim B$ stands for $A\leq C B$ whenever $C$ is a harmless constant.

\section{Solution formula for some reduced resolvent problem} \label{sec:sf}
In this section,  we consider the following reduced problem:
	\begin{equation}\label{eq:2.3}
		\left\{ \begin{aligned}
			&\gamma_{1+}\lambda\bv_+-\mu_+\Delta\bv_+-(\nu_++\delta)\nabla\di\bv_+=0
			&&\quad&\text{in}& \quad \Omega_+, \\
			&\gamma_{1-}\lambda\bv_--\mu_-\Delta\bv_-+\nabla \fp_-=0,\quad \di\bv_-=0
			&&\quad&\text{in}& \quad \Omega_-, \\ 
			&\bS_{\delta+}(\bv_+)\bn-\bT_-(\bv_-,\fp_-)\bn=\bh
			&&\quad&\text{on}& \quad \Gamma,\\
			& \bv_+-\bv_-=\bk
			&&\quad&\text{on}& \quad \Gamma ,	\\
			&\bv_-= 0
			&&\quad&\text{on}& \quad S.   
		\end{aligned}
		\right.
	\end{equation}	
where $\Omega_{\pm}, \Gamma, S$ are given in \eqref{dfn:Omega} for some $b>0$ and 
\begin{equation}\label{eq:S_delta}
  \bS_{\delta+}(\bv_+)
	=\mu_+\bD(\bv_+)+(\nu_+-\mu_++\delta)\di\bv_+\bI.  
\end{equation}
In fact, we can obtain the system \eqref{eq:2.3} by the eliminating $\rho_+$ from \eqref{eq:res_full} and neglecting the source terms in the equations of $\bv_{\pm}$.
In addition, let us set 
\begin{equation}\label{eq:KSigma}
\begin{aligned}
	K_{\ep}&=\{\lambda \in \BC: (\Re \lambda +\gamma_{1+}\gamma_{2+} \nu_+^{-1} +\ep)^2 + (\Im \lambda)^2 \geq (\gamma_{1+}\gamma_{2+} \nu_+^{-1} +\ep)^2\}, \\
	\Sigma_{\ep}&= \big\{\lambda \in \BC\backslash \{0\} :
	|\arg \lambda |\leq \pi -\ep \big\},\quad 
	\Sigma_{\ep, \lambda_0} =\{\lambda \in \Sigma_{\ep} : |\lambda| \geq \lambda_0\}
\end{aligned}
\end{equation}
for any $\lambda_0>0$ and $0<\varepsilon<\pi/2$.

\begin{ass}\label{ass:dl}
Let $\delta_0$ be some positive constant.
In \eqref{eq:2.3}, $\lambda,\delta\in \BC$  are some parameters satisfying the following assumption:
\begin{enumerate}[label=$(C\arabic*)$]	
\item $\lambda \in  \Sigma_{\ep, \lambda_0}  \cap K_{\ep},$ 
$\delta=\gamma_{1+}\gamma_{2+}\lambda^{-1};$

\item $\delta \in \Sigma_{\ep}$ 
with $\Re \delta<0$ and $|\delta| \leq \delta_0,$ 
$|\lambda|\geq \lambda_0,$
$\Re \lambda \geq |\Re \delta / \Im \delta| \, |\Im \lambda| ;$

\item $\delta \in \Sigma_{\ep}$ 
with $\Re \delta \geq 0$ and $|\delta| \leq \delta_0,$  
$|\lambda|\geq \lambda_0,$
$\Re \lambda \geq \lambda_0 \, |\Im \lambda|.$
\end{enumerate}

For simplicity, we also write
\begin{equation}\label{def:Gamma}
\Gamma_{\ep,\lambda_0}=
\begin{cases}
\Sigma_{\ep, \lambda_0}  \cap K_{\ep} & 
\text{in case of}\,\,\,\, (C1);\\
\big\{ \lambda\in \BC: |\lambda|\geq \lambda_0,
\Re \lambda \geq |\Re \delta / \Im \delta| \, |\Im \lambda|\big\}
& \text{in case of}\,\,\,\, (C2);\\
\big\{ \lambda\in \BC: |\lambda|\geq \lambda_0,
\Re \lambda \geq \lambda_0 \, |\Im \lambda|\big\}
& \text{in case of}\,\,\,\, (C3).
\end{cases}
\end{equation}
\end{ass}

As the the boundaries of $\Omega_\pm$ are flat, we introduce the (partial) Fourier transform with respect to the tangential variable $x^\prime:$
	$$\wh{\bv}(\xi^\prime,x_N)=\mathcal{F}_{x^\prime}[\bv(\cdot,x_N)](\xi^\prime)=\int_{\mathbb{R}^{N-1}}e^{-ix^\prime\cdot\xi^\prime}\bv(x^\prime,x_N)dx^\prime.$$ 
	Then taking the partial Fourier transform to \eqref{eq:2.3}  yields  
	\begin{equation}\label{eq:2.4}
	\left\{ 
		\begin{aligned}
		&(\gamma_{1+}\lambda+\mu_+|\xi^\prime|^2)\widehat{v}_{+j}
		-\mu_+\pd_N^2\widehat{v}_{+j}-(\nu_++\delta)i\xi_j
		(i\xi^\prime\cdot\widehat{\bv}^\prime_++\pd_N\widehat{v}_{+N})=0 
		& (x_N>0),\\
		&(\gamma_{1+}\lambda+\mu_+|\xi^\prime|^2)\widehat{v}_{+N}
		-\mu_+\pd_N^2\widehat{v}_{+N}-(\nu_++\delta)\pd_N
		(i\xi^\prime\cdot\widehat{\bv}^\prime_++\pd_N\widehat{v}_{+N})=0 
		& (x_N>0),\\	
		&(\gamma_{1-}\lambda+\mu_-|\xi^\prime|^2)\widehat{v}_{-j}
		-\mu_-\pd_N^2\widehat{v}_{-j}+i\xi_j\widehat{\fp}_-=0
		&(-b<x_N<0),\\
		&(\gamma_{1-}\lambda+\mu_-|\xi^\prime|^2)\widehat{v}_{-N}
		-\mu_-\pd_N^2\widehat{v}_{-N}+\pd_N\widehat{\fp}_-=0
		&(-b<x_N<0),\\
		&i\xi^\prime\cdot\widehat{\bv}^\prime_-+\pd_N\widehat{v}_{-N}=0
		&(-b<x_N<0)
	\end{aligned}
		\right.
	\end{equation}
	subject to the conditions
	\begin{equation}\label{eq:2.5}
		\left\{ 
		\begin{aligned}
		&\mu_+(\pd_N\widehat{v}_{+j}+i\xi_j\widehat{v}_{+N})
		-\mu_-(\pd_N\widehat{v}_{-j}+i\xi_j\widehat{v}_{-N})=\widehat{h}_j
		&(x_N=0),\\
		&2\mu_+\pd_N\widehat{v}_{+N}
		+(\nu_+-\mu_++\delta)(i\xi^\prime\cdot\widehat{\bv}_{+}^\prime
		+\pd_N\widehat{v}_{+N})
		-(2\mu_-\pd_N\widehat{v}_{-N}-\widehat{\fp}_-)=\widehat{h}_{N} 
		&(x_N=0),\\
		&\widehat{v}_{+J}-\widehat{v}_{-J}=\widehat{k}_J 
		&(x_N=0),\\
		&\widehat{v}_{-J}=0 &(x_N=-b)
	\end{aligned}
		\right.
	\end{equation}
	for $j=1,\cdots,N-1,$ $J=1,\cdots,N$ and $\pd_N=\pd_{x_N}.$ 
	Moreover, we have also used the equalities
	\begin{equation*}
		i\xi^\prime\cdot\widehat{\bv}_{\pm}^\prime=\sum_{j=1}^{N-1}i\xi_j\widehat{v}_j
	\end{equation*}
	for $\bv_\pm=(\bv_{\pm}',v_{\pm N})^{\top} =(v_{\pm 1},\cdots,v_{\pm (N-1)},v_{\pm N})^{\top}.$
	\medskip
	
Now, let us derive the solution formula from \eqref{eq:2.4} and \eqref{eq:2.5}.
As the discussion of the problem in the whole space in \cite{KSS2014}, we introduce the characteristic roots:
	\begin{equation}\label{eq:2.7}
		\begin{aligned}
			A_+& =\sqrt{\gamma_{1+}(\mu_++\nu_++\delta)^{-1}\lambda+A^2}, & 
			B_{+}& =\sqrt{\gamma_{1+}(\mu_+)^{-1}\lambda+A^2},\\ 
			A_-& =A=|\xi^\prime|, & 
			B_-& =\sqrt{\gamma_{1-}(\mu_-)^{-1}\lambda+A^2}.
		\end{aligned}
	\end{equation}
Suppose that the solution $(\widehat{\bv}_{\pm},\widehat{\fp}_-)$ to the problem \eqref{eq:2.4}-\eqref{eq:2.5} is of the forms 
	\begin{equation}\label{eq:2.8}
		\begin{aligned}
			\widehat{v}_{+J}=&\alpha_{+J}(e^{-B_+x_N}-e^{-A_+x_N})+\beta_{+J}e^{-B_+x_N},\\
			\widehat{v}_{-J}=&\alpha_{-J}^0(e^{B_-x_N}-e^{Ax_N})+\beta_{-J}^0e^{B_-x_N}\\
			&+\alpha_{-J}^{b}(e^{-B_-(x_N+b)}-e^{-A(x_N+b)})+\beta_{-J}^{b}e^{-B_-(x_N+b)},\\
			\widehat{\fp}_-=&\gamma_-^0e^{Ax_N} +\gamma_-^be^{-A(x_N+b)} .
		\end{aligned}
	\end{equation}
Then inserting \eqref{eq:2.8} into \eqref{eq:2.4} and equating the coefficients of the terms $e^{\mp A_{\pm}x_N},$ $e^{\mp B_{\pm}x_N}$, $e^{-A_-(x_N+b)}$ and $e^{-B_-(x_N+b)}$ imply that 
\begin{equation}\label{eq:2.9}
	\left\{ 
\begin{aligned}
		&\mu_+(A^2_+-B^2_+)\alpha_{+j}+(\nu_++\delta)i\xi_j(i\xi^\prime\cdot\alpha^\prime_+-A_+\alpha_{+N})=0,\\
		&\mu_+(A^2_+-B^2_+)\alpha_{+N}-(\nu_++\delta)A_+(i\xi^\prime\cdot\alpha^\prime_+-A_+\alpha_{+N})=0,\\
		&i\xi^\prime\cdot\alpha^\prime_+-B_+\alpha_{+N}
		+i\xi^\prime\cdot\beta^\prime_+-B_+\beta_{+N}=0,\\
				&\mu_-(A^2-B^2_-)\alpha_{-j}^0+i\xi_j\gamma_-^0=0,\quad \mu_-(A^2-B^2_-)\alpha_{-j}^b+i\xi_j\gamma_-^b=0,\\
		&\mu_-(A^2-B^2_-)\alpha_{-N}^0+A\gamma_-^0=0,\quad \mu_-(A^2-B^2_-)\alpha_{-N}^b-A\gamma_-^b=0,\\
			&i\xi^\prime\cdot(\alpha_-^{0})'+B_-\alpha_{-N}^0
	+i\xi^\prime\cdot(\beta_-^{0})'+B_-\beta_{-N}^0=0,\\
	&i\xi^\prime\cdot(\alpha_-^{b})'-B_-\alpha_{-N}^b
	+i\xi^\prime\cdot(\beta_-^{b})'-B_-\beta_{-N}^b=0,
	\\&i\xi^\prime\cdot(\alpha_-^{0})'+A\alpha_{-N}^0=0,\quad i\xi^\prime\cdot(\alpha_-^{b})'-A\alpha_{-N}^b=0.
\end{aligned}
\right.
\end{equation}
From the equalities in \eqref{eq:2.9}, it is not hard to see that
\begin{equation}\label{eq:2.10}
\left\{ 
\begin{aligned}
	i\xi^\prime\cdot\alpha_+^\prime
	&=\frac{A^2}{A_+B_+-A^2}(i\xi^\prime\cdot\beta_+^\prime-B_+\beta_{+N}), \\
	\alpha_{+N} &=\frac{A_+}{A_+B_+-A^2}
	(i\xi^\prime\cdot\beta_+^\prime-B_+\beta_{+N}),\\
i\xi^\prime\cdot(\alpha_-^{0})'
&=\frac{A}{B_--A} 
\big(i\xi^\prime\cdot(\beta_-^{0})'+B_-\beta_{-N}^{0} \big),\\
\alpha_{-N}^0&=-\frac{1}{B_--A} 
\big(i\xi^\prime\cdot(\beta_-^{0})'+B_-\beta_{-N}^{0} \big),\\
i\xi^\prime\cdot(\alpha_-^{b})' &=\frac{A}{B_--A}
\big( i\xi^\prime\cdot(\beta_-^{b})'-B_-\beta_{-N}^{b} \big),\\
\alpha_{-N}^b &=\frac{1}{B_--A}
\big( i\xi^\prime\cdot(\beta_-^{b})'-B_-\beta_{-N}^{b}\big),\\
\gamma_-^0&=-\frac{\mu_-(A+B_-)}{A}
\big( i\xi^\prime\cdot(\beta_-^{0})'+B_-\beta_{-N}^0 \big),\\
\gamma_-^b&=-\frac{\mu_-(A+B_-)}{A}
\big( i\xi^\prime\cdot(\beta_-^{b})'-B_-\beta_{-N}^b \big).
\end{aligned}
\right.
\end{equation}

On the other hand,  we use \eqref{eq:2.8} to get
\begin{equation}\label{2.11-1}
\left\{ 
	\begin{aligned}
		\widehat{v}_{+ J}(0)& =\beta_{+ J},\\
		\pd_N\widehat{v}_{+ J}(0)&=(A_+-B_+)\alpha_{+J}-B_+\beta_{+J},\\
		\widehat{v}_{- J}(0)&=\beta_{-J}^0+(e^{-B_-b}-e^{-Ab})\alpha_{-J}^b
		+e^{-B_-b}\beta_{-J}^b,\\
		\pd_N\widehat{v}_{- J}(0)&=(B_--A)\alpha_{-J}^0+B_-\beta_{-J}^0\\
		&\quad +(-B_-e^{-B_-b}+Ae^{-Ab})\alpha_{-J}^b-B_-e^{-B_-b}\beta_{-J}^b, \\
	\widehat{p}_{-}(0)&=\gamma_-^0+e^{-Ab}\gamma_-^b,\\
	\widehat{v}_{-J}(-b) &=(e^{-B_-b}-e^{-Ab})\alpha_{-J}^0
	+e^{-B_-b}\beta_{-J}^0+\beta_{-J}^{b}.
	\end{aligned}
	\right.
\end{equation}
Substituting \eqref{2.11-1} into \eqref{eq:2.5}, we have
\begin{equation}\label{eq:2.12}
\begin{aligned}
\widehat{h}_j(0)=&\mu_+\big((A_+-B_+)\alpha_{+j}-B_+\beta_{+j}+i\xi_j\beta_{+N}\big)\\
&-\mu_-\big[ (B_--A)\alpha_{-j}^0+B_-\beta_{-j}^0
+(-B_-e^{-B_-b}+Ae^{-Ab})\alpha_{-j}^b
\\& \hspace{1.5cm} -B_-e^{-B_-b}\beta_{-j}^b
+i\xi_j \big( \beta_{-N}^0+(e^{-B_-b}-e^{-Ab})\alpha_{-N}^b
+ e^{-B_-b} \beta_{-N}^b\big) \big] ,\\
\widehat{h}_N(0)=&-2\mu_+\big((B_+-A_+)\alpha_{+N}+B_+\beta_{+N}\big) \\
		&-(\nu_+-\mu_++\delta)\big(-i\xi^\prime\cdot\beta_+^\prime+(B_+-A_+)\alpha_{+N}+B_+\beta_{+N}\big)\\
		&-2\mu_-\big[ (B_--A)\alpha_{-N}^0+B_-\beta_{-N}^0+(-B_-e^{-B_-b}+Ae^{-Ab})\alpha_{-N}^b-B_-e^{-B_-b}\beta_{-N}^b\big]\\
		&+\gamma_-^0+e^{-Ab}\gamma_-^b,\\
\widehat{k}_J(0)=&\beta_{+J}-\beta_{-J}^0
-(e^{-B_-b}-e^{-Ab})\alpha_{-J}^b-e^{-B_-b}\beta_{-J}^b,\\
0=&(e^{-B_-b}-e^{-Ab})\alpha_{-J}^0+e^{-B_-b}\beta_{-J}^0+\beta_{-J}^b.
	\end{aligned}
\end{equation}
\medskip

Let $M_{\pm}(x_N)$ denote the Stokes kernels 
\begin{equation}\label{eq:2.13}
	M_+(y_N)=\frac{e^{-B_+y_N}-e^{-A_+y_N}}{B_+-A_+} \,\,\,(y_N>0),\qquad 
	M_-(z_N)=\frac{e^{B_-z_N}-e^{Az_N}}{B_--A} \,\,\,(z_N<0).
\end{equation}
Using \eqref{eq:2.10} and \eqref{eq:2.12}, we get that
\begin{equation}\label{eq:2.14}
\begin{aligned}
i\xi'\cdot \widehat{h}'(0)=&\mu_+\Big( \frac{(A_+-B_+)A^2}{A_+B_+-A^2}(i\xi'\cdot\beta_+'-B_+\beta_{+N})-B_+i\xi'\cdot\beta_+'-A^2\beta_{+N} \Big)\\
	&-\mu_-\Big\{A\big( i\xi'\cdot(\beta_-^{0})'+B_-\beta_{-N}^0 \big)
	+B_-i\xi' \cdot (\beta_-^{0})' \\
	&\qquad -A\big( B_-M_-(-b)+e^{-Ab} \big) 
	\big( i\xi'\cdot(\beta_-^{b})'-B_-\beta_{-N}^b\big)
	- B_-e^{-B_-b}\,i\xi'\cdot(\beta_-^{b})'\\
	&\qquad -A^2\big[ \beta_{-N}^0+M_-(-b)
	\big( i\xi'\cdot(\beta_-^{b})'-B_-\beta_{-N}^b\big)
	+e^{-B_-b}\beta_{-N}^b \big]\Big\},\\
\widehat{h}_N(0)=&-2\mu_+\Big( \frac{A_+(B_+-A_+)}{A_+B_+-A^2}
(i\xi'\cdot\beta_+'-B_+\beta_{+N})+B_+\beta_{+N} \Big)\\
&-(\nu_+-\mu_++\delta)\Big(-i\xi'\cdot\beta_+'+\frac{A_+(B_+-A_+)}{A_+B_+-A^2}(i\xi'\cdot\beta_+'-B_+\beta_{+N})+B_+\beta_{+N}  \Big)\\
	&-2\mu_-\Big\{-\big( i\xi'\cdot(\beta_-^{0})'+B_-\beta_{-N}^0 \big)+B_-\beta_{-N}^0\\
	&\qquad -\big( B_-M_-(-b)+e^{-Ab} \big) 
	\big( i\xi'\cdot(\beta_-^{b})'-B_-\beta_{-N}^b \big)
	-B_-e^{-B_-b}\beta_{-N}^b \Big\} \\
	&-\frac{\mu_-(A+B_-)}{A}
\big( i\xi^\prime\cdot(\beta_-^{0})'+B_-\beta_{-N}^0 \big)
-\frac{\mu_-(A+B_-)}{A}e^{-Ab}
\big( i\xi^\prime\cdot(\beta_-^{b})'-B_-\beta_{-N}^b \big),\\
i\xi'\cdot\widehat{k}'(0)=&i\xi'\cdot\beta_+'-i\xi'\cdot(\beta_-^{0})'
-AM_-(-b) \big( i\xi'\cdot(\beta_-^{b})'-B_-\beta_{-N}^b \big)
-e^{-B_-b}i\xi'\cdot(\beta_-^{b})',\\
\widehat{k}_N(0)=&\beta_{+N}-\beta_{-N}^0-M_-(-b)
\big( i\xi'\cdot(\beta_-^{b})'-B_-\beta_{-N}^b \big)-e^{-B_-b}\beta_{-N}^b,\\
0=&AM_-(-b)\big(i\xi^\prime\cdot(\beta_-^{0})'+B_-\beta_{-N}^{0} \big)
+e^{-B_-b}i\xi^\prime\cdot(\beta_-^{0})'+
i\xi^\prime\cdot(\beta_-^{b})',\\
	0=&-M_-(-b)\big(i\xi^\prime\cdot(\beta_-^{0})'+B_-\beta_{-N}^{0} \big)
	+e^{-B_-b}\beta_{-N}^0+\beta_{-N}^b.\\
\end{aligned}
\end{equation}
It then follows that
\begin{equation}\label{eq:2.14-2}
\begin{aligned}
	i\xi'\cdot\widehat{h}'&=L_{11}^+ i\xi'\cdot\beta_+'
	+L_{11}^{0-} i\xi'\cdot(\beta_-^{0})' 
	+L_{11}^{b-} i\xi'\cdot(\beta_-^{b})'
	+L_{12}^+\beta_{+N}+L_{12}^{0-}\beta_{-N}^0+L_{12}^{b-}\beta_{-N}^b,\\
	A\widehat{h}_N&=L_{21}^+ i\xi'\cdot\beta_+'
	+L_{21}^{0-} i\xi'\cdot(\beta_-^{0})'
	+L_{21}^{b-} i\xi'\cdot(\beta_-^{b})'
	+L_{22}^+\beta_{+N}+L_{22}^{0-}\beta_{-N}^0+L_{22}^{b-}\beta_{-N}^b
\end{aligned}
\end{equation}
with
\begin{equation}\label{eq:2.15}
	\begin{aligned}
		&L_{11}^+=-\mu_+\frac{A_+(B_+^2-A^2)}{A_+B_+-A^2},
		\qquad L_{11}^{0-}=-\mu_-(A+B_-),\\
		&L_{11}^{b-}=\mu_-\big((A+B_-)e^{-B_-b}+2A^2M_-(-b)\big),\qquad
		L_{12}^+=-\mu_+\frac{A^2(2A_+B_+-A^2-B_+^2)}{A_+B_+-A^2},\\
		& L_{12}^{0-}=-\mu_-(B_--A)A,\qquad 
		 L_{12}^{b-}=-\mu_-\big(2A^2B_-M_-(-b)+(B_--A)Ae^{-B_-b}\big),\\
		&L_{21}^+=-A\Big( 2\mu_+\frac{A_+(B_+-A_+)}{A_+B_+-A^2}-(\nu_+-\mu_++\delta)\frac{A_+^2-A^2}{A_+B_+-A^2}\Big), \\
		& L_{21}^{0-}=-\mu_-(B_--A),\qquad
		L_{21}^{b-}=2\mu_- AB_-M_-(-b)
		-\mu_-(B_--A)e^{-Ab},\\
		&L_{22}^+=-(\mu_++\nu_++\delta)\frac{AB_+(A^2_+-A^2)}{A_+B_+-A^2},
		\qquad L_{22}^{0-}=-\mu_-(A+B_-)B_-,\\
		&L_{22}^{b-}=-2\mu_-A^2B_-M_-(-b)
		+\mu_-(B_-+A)B_-e^{-Ab}.
	\end{aligned}
\end{equation}
Using the equations $\eqref{eq:2.14}_3$ and $\eqref{eq:2.14}_4$, we have
\begin{equation}\label{eq:2.15-2}
\left\{ 
\begin{aligned}
i\xi'\cdot\beta_+'=
&i\xi'\cdot\widehat{k}'(0)-A(A+B_-)M_-^2(-b) 
\big( i\xi'\cdot(\beta_-^{0})'+B_-\beta_{-N}^0\big)\\
&+\big(1-e^{-2B_-b} -2e^{-B_-b}AM_-(-b)\big) i\xi'\cdot(\beta_-^{0})' ,\\
\beta_{+N} =& \widehat{k}_N(0)-(A+B_-)M_-^2(-b)
\big( i\xi'\cdot(\beta_-^{0})'+B_-\beta_{-N}^0\big)\\
&+ \big( 1-e^{-2B_-b}+2e^{-B_-b}B_-M_-(-b) \big)\beta_{-N}^0,\\
i\xi^\prime\cdot(\beta_-^{b})' =& 
-e^{-B_-b} i\xi'\cdot(\beta_-^{0})'
-AM_-(-b)\big( i\xi'\cdot(\beta_-^{0})'+B_-\beta_{-N}^0\big),\\
\beta_{-N}^b =& M_-(-b) \big( i\xi'\cdot(\beta_-^{0})'+B_-\beta_{-N}^0\big)
-e^{-B_-b}\,\beta_{-N}^0.
\end{aligned}\right.
\end{equation}

Now let us set 
\begin{equation*}
\begin{aligned}
\widehat{H}_h=&i\xi'\cdot\widehat{h}'
-L_{11}^+i\xi'\cdot\widehat{k}'-L_{12}^+\widehat{k}_N,\\
\widehat{H}_N=&A\widehat{h}_N-L_{21}^+
i\xi'\cdot\widehat{k}'-L_{22}^+\widehat{k}_N,\\
\fL_1=&-(AL_{11}^++L_{12}^+) (A+B_-)M_-^2(-b)
-L_{11}^{b-}A M_-(-b)+ L_{12}^{b-}M_-(-b),\\
\fL_2=&-(AL_{21}^++L_{22}^+ ) (A+B_-)M_-^2(-b)
-L_{21}^{b-}A M_-(-b)+ L_{22}^{b-}M_-(-b).
\end{aligned}
\end{equation*}
Inserting \eqref{eq:2.15-2} into \eqref{eq:2.14-2} yields  
\begin{equation*}
	\bL\begin{pmatrix}
		i\xi'\cdot(\beta_-^{0})'\\\beta_{-N}^0
	\end{pmatrix}=\begin{pmatrix}
		\widehat{H}_h(0)\\\widehat{H}_N(0)
	\end{pmatrix}
\end{equation*}
where the entries $L_{k\ell}$ of the matrix $\bL$ are given by 
\begin{equation}\label{eq:Lij}
\begin{aligned}
L_{k1}=&L_{k1}^+\big(1-e^{-2B_-b}-2e^{-B_-b}AM_-(-b)\big)
+L_{k1}^{0-}-L_{k1}^{b-}e^{-B_-b}+\fL_k,\\
L_{k2}=&L_{k2}^+\big(1-e^{-2B_-b}+2e^{-B_-b}B_-M_-(-b)\big)
+L_{k2}^{0-}-L_{k2}^{b-}e^{-B_-b}+B_-\fL_k
\end{aligned}
\end{equation}
for $k=1,2.$ 
\medskip

Suppose that
\footnote{See Lemma \ref{lem:3.1} for more details. } 
$\det \bL=L_{11}L_{22}-L_{12}L_{21}\not=0.$ Then the inverse matrix of $\bL$ is written as
\begin{equation*}
	\bL^{-1}=\frac{1}{\det \bL}\begin{pmatrix}
		L_{22}&-L_{12}\\
		-L_{21}&L_{11}
	\end{pmatrix}.
\end{equation*}
Thus we have
\begin{equation}\label{eq:2.16}
\begin{aligned}
i\xi'\cdot(\beta_-^{0})'
=&\frac{1}{\det \bL}\Big(L_{22}i\xi'\cdot\widehat{h}'(0)
		-AL_{12}\widehat{h}_N(0) \\
&\hspace{1.5cm}-(L_{22}L_{11}^+-L_{12}L_{21}^+)i\xi'\cdot\widehat{k}'(0)
		-(L_{22}L_{12}^+-L_{12}L_{22}^+)\widehat{k}_N(0)\Big),\\
\beta_{-N}^{0}
=&\frac{1}{\det \bL}\Big(-L_{21}i\xi'\cdot\widehat{h}'(0)
		+AL_{11}\widehat{h}_N(0)\\
&\hspace{1.5cm}+(L_{21}L_{11}^+-L_{11}L_{21}^+)i\xi'\cdot\widehat{k}'(0)
		+(L_{21}L_{12}^+-L_{11}L_{22}^+)\widehat{k}_N(0)\Big).
	\end{aligned}
\end{equation}
By \eqref{eq:2.15-2} and \eqref{eq:2.16}, we arrive 
\begin{equation}\label{eq:2.16-2}
\begin{aligned}
i\xi'\cdot(\beta_-^{0})'+B_-\beta_{-N}^0
&=\sum_{J=1}^{N}\big(P_{J,0}^{0}\widehat{h}_J(0)
+P_{J,1}^{0}\widehat{k}_{J}(0)\big),\\
i\xi'\cdot(\beta_-^{b})'-B_-\beta_{-N}^b
&=\sum_{J=1}^{N}\big(P_{J,0}^{b}\widehat{h}_J(0)
+P_{J,1}^{b}\widehat{k}_{J}(0)\big),\\
i\xi'\cdot\beta_+^{'}-B_+\beta_{+N}
&=\sum_{J=1}^{N}\big(P_{J,0}^{+}\widehat{h}_J(0)
+P_{J,1}^{+}\widehat{k}_{J}(0)\big).
	\end{aligned}
\end{equation}
Here, we have defined 
\begin{equation*}
\begin{aligned}
P_{j,0}^{0}=&\frac{(L_{22}-B_-L_{21})i\xi_j}{\det\bL},\\ 
P_{N,0}^{0}=&\frac{-AL_{12}+AB_-L_{11}}{\det\bL},\\
P_{j,1}^{0}=&\frac{\big(L_{12}L_{21}^+-L_{22}L_{11}^++B_-(L_{21}L_{11}^+-L_{11}L_{21}^+)\big)i\xi_j}{\det\bL},\\
P_{N,1}^{0}=&\frac{L_{12}L_{22}^+-L_{22}L_{12}^+
+B_-(L_{21}L_{12}^+-L_{11}L_{22}^+)}{\det\bL},\\
P_{j,0}^{b}=&-\frac{\big((A+B_-)M_-(-b)+e^{-B_-b}\big)L_{22}i\xi_j}{\det\bL}\\
&+\frac{\big((A+B_-)M_-(-b)-e^{-B_-b}\big)B_-L_{21}i\xi_j}{\det\bL},\\ 
P_{N,0}^{b}=&\frac{\big((A+B_-)M_-(-b)+e^{-B_-b}\big)AL_{12}}{\det\bL}\\
&-\frac{\big((A+B_-)M_-(-b)-e^{-B_-b}\big)AB_-L_{11}}{\det\bL},\\
P_{j,1}^{b}=&\frac{(L_{22}L_{11}^+-L_{12}L_{21}^+)
\big((A+B_-)M_-(-b)+e^{-B_-b}\big)i\xi_j}{\det\bL}\\
&-\frac{(L_{21}L_{11}^+-L_{11}L_{21}^+)
\big((A+B_-)M_-(-b)-e^{-B_-b}\big)B_-i\xi_j}{\det\bL},\\
P_{N,1}^{b}=&\frac{(L_{22}L_{12}^+-L_{12}L_{22}^+)
\big((A+B_-)M_-(-b)+e^{-B_-b}\big)}{\det\bL}\\
&-\frac{(L_{21}L_{12}^+-L_{11}L_{22}^+)
\big((A+B_-)M_-(-b)-e^{-B_-b}\big)B_-}{\det\bL},\\
	P_{j,0}^+=& N^{+}_{0} P_{j,0}^0+N^+_1 \frac{L_{21}}{\det\bL} i\xi_j,\\
	P_{N,0}^+=&N^{+}_{0} P_{N,0}^0-N^{+}_1\frac{AL_{11}}{\det\bL},  \\
    P_{j,1}^+=& i\xi_j + N^{+}_{0}  P_{j,1}^0
    +N^{+}_1\frac{L_{11}L_{21}^+-L_{21}L_{11}^+}{\det\bL} i\xi_j,\\
	 P_{N,1}^+=&-B_++N^{+}_{0} P_{N,1}^0
		+N^{+}_1 \frac{L_{11}L_{22}^+-L_{21}L_{12}^+}{\det\bL}
	\end{aligned}
\end{equation*}
with
\begin{equation*}
\begin{aligned}
N^{+}_0 &=(B_+-A)(A+B_-)M_-^2(-b)
	+1-e^{-2B_-b}-2e^{-B_-b}AM_-(-b),\\
N^{+}_1 &= (B_++B_-)(1-e^{-2B_-b}) +2B_-(B_+-A)M_-(-b)e^{-B_-b}.
\end{aligned}
\end{equation*}
In particular, the formula of $\widehat{\fp}_-$ reads as
\begin{equation}\label{eq:pressure}
	\widehat{\fp}_-(x_N)=e^{Ax_N}\sum_{J=1}^N
	\big( p_{J,0}^{0}\widehat{h}_J(0)+p_{J,1}^{0}\widehat{k}_J(0)\big)
	+e^{-A(x_N+b)}\sum_{J=1}^N\big(p_{J,0}^{b}\widehat{h}_J(0)
	+p_{J,1}^{b}\widehat{k}_J (0)\big)
\end{equation}
with
\begin{align*}
	p_{J,0}^{0}&=-\frac{\mu_-(A+B_-)}{A}P_{J,0}^{0}, &
	 p_{J,1}^{0}=-\frac{\mu_-(A+B_-)}{A}P_{J,1}^{0},\\
	p_{J,0}^{b}&=-\frac{\mu_-(A+B_-)}{A}P_{J,0}^{b},&
	p_{J,1}^{b}=-\frac{\mu_-(A+B_-)}{A}P_{J,1}^{b}.
\end{align*}
\medskip


Now,  it suffices to derive the formulas of $\wh \bv_{\pm}.$ 
Using \eqref{eq:2.10} and \eqref{eq:2.16-2}, we obtain
\begin{equation}\label{eq:cof_1}
\left\{ 
\begin{aligned}
\alpha_{+N} &=\frac{A_+}{A_+B_+-A^2}
\sum_{J=1}^{N}\big(P_{J,0}^{+}\widehat{h}_J(0)
+P_{J,1}^{+}\widehat{k}_{J}(0)\big),\\
\alpha_{-N}^0&=- \frac{1}{B_--A}\sum_{J=1}^{N}
\big(P_{J,0}^{0}\widehat{h}_J(0)+P_{J,1}^{0}\widehat{k}_{J}(0)\big),\\
\alpha_{-N}^b &=\frac{1}{B_--A}\sum_{J=1}^{N}
\big(P_{J,0}^{b}\widehat{h}_J(0)+P_{J,1}^{b}\widehat{k}_{J}(0)\big),\\
\gamma_-^0&=-\frac{\mu_-(A+B_-)}{A}
\sum_{J=1}^{N}
\big(P_{J,0}^{0}\widehat{h}_J(0)+P_{J,1}^{0}\widehat{k}_{J}(0)\big),\\
\gamma_-^b&=-\frac{\mu_-(A+B_-)}{A}\sum_{J=1}^{N}
\big(P_{J,0}^{b}\widehat{h}_J(0)+P_{J,1}^{b}\widehat{k}_{J}(0)\big).
\end{aligned}
\right.
\end{equation}
On the other hand,  notice from \eqref{eq:2.10} that 
\begin{equation*}
\begin{aligned}
i\xi^\prime\cdot\alpha_+^\prime-A_+\alpha_{+N}
&=-\frac{A_+^2-A^2}{A_+B_+-A^2}
(i\xi^\prime\cdot\beta_+^\prime-B_+\beta_{+N}) \\
&=-\frac{\mu_+ (B_+^2-A_+^2)}{(\nu_++\delta)(A_+B_+-A^2)}
(i\xi^\prime\cdot\beta_+^\prime-B_+\beta_{+N})
\end{aligned}
\end{equation*}
for $(\nu_++\delta)(A_+^2-A^2)=\mu_+ (B_+^2-A_+^2).$
Then \eqref{eq:2.9} and \eqref{eq:2.16-2} yield 
\begin{equation}\label{eq:cof_2}
\left\{
\begin{aligned}
\alpha_{+j}&=-\frac{i\xi_j}{A_+B_+-A^2}
\sum_{J=1}^{N}\big(P_{J,0}^{+}\widehat{h}_J(0)
+P_{J,1}^{+}\widehat{k}_{J}(0)\big),\\
\alpha_{-j}^0&=\frac{i\xi_j \gamma_-^0}{\mu_-(B_-^2-A^2)}
=-\frac{i\xi_j}{A(B_--A)}
\sum_{J=1}^{N}\big(P_{J,0}^{0}\widehat{h}_J(0)
+P_{J,1}^{0}\widehat{k}_{J}(0)\big),\\ 
\alpha_{-j}^b&=\frac{i\xi_j \gamma_-^b}{\mu_-(B_-^2-A^2)}
=-\frac{i\xi_j}{A(B_--A)}
\sum_{J=1}^{N}\big(P_{J,0}^{b}\widehat{h}_J(0)
+P_{J,1}^{b}\widehat{k}_{J}(0)\big)
	\end{aligned} \right.
\end{equation}
for any $j=1,\dots,N-1.$
\medskip

Next, using \eqref{eq:2.16}, we represent $\beta_{-N}^0$ as follows
\begin{equation}\label{eq:cof_3}
\beta_{-N}^0=\sum_{J=1}^N \big( S_{NJ,-1}^-\widehat{h}_J(0)
+S_{NJ,0}^-\widehat{k}_{J}(0) \big),
\end{equation}
where we have set 
\begin{equation*}
\begin{aligned}
S_{Nj,-1}^-&=-\frac{L_{21} i\xi_j}{\det \bL},
&S_{NN,-1}^-&=\frac{AL_{11}}{\det \bL},\\
S_{Nj,0}^-&=\frac{(L_{21}L_{11}^+-L_{11}L_{21}^+)i\xi_j}{\det \bL},
&S_{NN,0}^-&=\frac{L_{21}L_{12}^+-L_{11}L_{22}^+}{\det \bL}
	\end{aligned}
\end{equation*}
for any $j=1,\dots,N-1.$
Then inserting \eqref{eq:cof_3} into \eqref{eq:2.15-2}, we obtain that
\begin{equation}\label{eq:cof_4}
\begin{aligned}
\beta_{+N}=&\sum_{J=1}^N \big( S_{NJ,-1}^+\widehat{h}_J(0)
+S_{NJ,0}^+\widehat{k}_{J}(0) \big),\\
\beta_{-N}^b=&\sum_{J=1}^N \big( S_{NJ,-1}^b\widehat{h}_J(0)
+S_{NJ,0}^b\widehat{k}_{J}(0) \big)
	\end{aligned}
\end{equation}
with
\begin{equation*}
\begin{aligned}
S_{NJ,-1}^+=&-(A+B_-)M_-^2(-b)P^{0}_{J,0}\\
&+\big( 1-e^{-2B_-b}+2e^{-B_-b}B_-M_-(-b) \big)S_{NJ,-1}^-,\\
S_{NJ,0}^+=&\delta_{JN}-(A+B_-)M_-^2(-b)P^{0}_{J,1}\\
&+\big( 1-e^{-2B_-b}+2e^{-B_-b}B_-M_-(-b) \big)S_{NJ,0}^-,\\ 
S_{NJ,-1}^b=& M_-(-b)P^{0}_{J,0}-e^{-B_-b} S_{NJ,-1}^-, \\
S_{NJ,0}^b =& M_-(-b)P^{0}_{J,1}-e^{-B_-b} S_{NJ,0}^-
\end{aligned}
\end{equation*}
for $J=1,\dots,N.$ 
\medskip

From the last two equations of \eqref{eq:2.12},  we have
\begin{equation}\label{eq:2.17-1}
\begin{aligned}
\beta_{-j}^b=&-e^{-B_-b}\beta_{-j}^0-(B_--A)M_-(-b)\alpha_{-j}^0,\\
\beta_{+j}=&(1-e^{-2B_-b})\beta_{-j}^0+\widehat{k}_j(0)
-e^{-B_-b}(B_--A)M_-(-b)\alpha_{-j}^0\\
&+(B_--A)M_-(-b)\alpha_{-j}^b
	\end{aligned}
\end{equation}
for $j=1,\dots,N-1.$  Then inserting \eqref{eq:2.17-1} into $\eqref{eq:2.12}_1$ implies that 
\begin{equation*}
\begin{aligned}
\beta_{-j}^0=\frac{1}{D(B_{\pm})}\Big(&
-\widehat{h}_j(0)-\mu_+B_+\widehat{k}_j(0)
-\mu_+(B_+-A_+)\alpha_{+j}+\mu_+i\xi_j\beta_{+N}\\
&+\big( (\mu_+B_+-\mu_-B_-)e^{-B_-b}M_-(-b)-\mu_-\big)(B_--A) \alpha_{-j}^0
-\mu_-i\xi_j  \beta_{-N}^0\\
&  +((\mu_{-}B_--\mu_+B_+)M_-(-b)+\mu_{-}e^{-Ab})(B_--A) \alpha_{-j}^b \\
&-\mu_- i\xi_j M_-(-b) (B_--A) \alpha_{-N}^b 
-\mu_-i\xi_j e^{-B_-b} \beta_{-N}^b \Big),
	\end{aligned}
\end{equation*}
where we denote for simplicity 
$$D(B_{\pm})=\mu_+B_+(1-e^{-2B_-b})+\mu_-B_-(1+e^{-2B_-b}).$$
Then \eqref{eq:cof_1}, \eqref{eq:cof_2}, \eqref{eq:cof_3} and \eqref{eq:cof_4} furnish that 
\begin{equation}\label{eq:cof_5}
\beta_{-j}^0=\sum_{\ell=1}^N \big( S_{j\ell,-1}^-\widehat{h}_\ell(0)
+S_{j\ell,0}^-\widehat{k}_{\ell}(0) \big)
\end{equation}
for $j=1,\dots,N-1,$ provided with
\begin{equation*}
\begin{aligned}
S_{j\ell,-1}^- =\frac{1}{D(B_{\pm})}\Big(&
-\delta_{j\ell}
+\mu_+\frac{i\xi_j(B_+-A_+)}{A_+B_+-A^2}P_{\ell,0}^{+}
+\mu_+i\xi_j S_{N\ell,-1}^+\\
&-\big( (\mu_+B_+-\mu_-B_-)e^{-B_-b}M_-(-b)-\mu_-\big)\frac{i\xi_j}{A}
P_{\ell,0}^{0}-\mu_-i\xi_j   S_{N\ell,-1}^-\\
&
 -\big(  (\mu_{-}B_--\mu_+B_+)M_-(-b)+\mu_{-}e^{-Ab}\big) \frac{i\xi_j}{A}P_{\ell,0}^{b}\\
&-\mu_- i\xi_j M_-(-b)P_{\ell,0}^{b} 
-\mu_-i\xi_j e^{-B_-b} S_{N\ell,-1}^b\Big),\\
S_{j\ell,0}^-=\frac{1}{D(B_{\pm})}\Big(&
-\mu_+B_+\delta_{j\ell}
+\mu_+\frac{i\xi_j(B_+-A_+)}{A_+B_+-A^2}P_{\ell,1}^{+}
+\mu_+i\xi_j S_{N\ell,0}^+\\
&-\big( (\mu_+B_+-\mu_-B_-)e^{-B_-b}M_-(-b)-\mu_-\big)\frac{i\xi_j}{A}
P_{\ell,1}^{0}-\mu_-i\xi_j  S_{N\ell,0}^- \\
&
 -\big( (\mu_{-}B_--\mu_+B_+)M_-(-b)+\mu_{-}e^{-Ab}\big) \frac{i\xi_j}{A}
P_{\ell,1}^{b}\\
&-\mu_- i\xi_j M_-(-b)P_{\ell,1}^{b}
-\mu_-i\xi_j e^{-B_-b} S_{N\ell,0}^b\Big).
\end{aligned}
\end{equation*}
Then combining \eqref{eq:2.17-1} and \eqref{eq:cof_5} implies that
\begin{equation}\label{eq:cof_6}
\begin{aligned}
\beta_{-j}^b&=\sum_{\ell=1}^N \big( S_{j\ell,-1}^b\widehat{h}_\ell(0)
+S_{j\ell,0}^b\widehat{k}_{\ell}(0) \big),\\
\beta_{+j}&=\sum_{\ell=1}^N \big( S_{j\ell,-1}^+\widehat{h}_\ell(0)
+S_{j\ell,0}^+\widehat{k}_{\ell}(0) \big)
\end{aligned}
\end{equation}
with 
\begin{equation*}
\begin{aligned}
S_{j\ell,-1}^b=&-e^{-B_-b} S_{j\ell,-1}^-
+M_-(-b)\frac{i\xi_j}{A}
P_{\ell,0}^{0},\\
S_{j\ell,0}^b=&-e^{-B_-b}
S_{j\ell,0}^-
+M_-(-b)\frac{i\xi_j}{A}P_{\ell,1}^{0},\\
S_{j\ell,-1}^+=&(1-e^{-2B_-b}) S_{j\ell,-1}^-
+e^{-B_-b}M_-(-b)\frac{i\xi_j}{A} P_{\ell,0}^{0}
-M_-(-b)\frac{i\xi_j}{A}P_{\ell,0}^{b},\\
S_{j\ell,0}^+=&\delta_{j\ell}
+(1-e^{-2B_-b}) S_{j\ell,0}^-
+e^{-B_-b}M_-(-b)\frac{i\xi_j}{A} P_{\ell,1}^{0}
-M_-(-b)\frac{i\xi_j}{A} P_{\ell,1}^{b}.
\end{aligned}
\end{equation*}
\medskip

At last, let us introduce the symbols
\begin{gather*}
	R_{j\ell,0}^+=\frac{(B_+-A_+)i\xi_j}{A_+B_+-A^2}P_{\ell,0}^+,
	\quad R_{N\ell,0}^+=\frac{(B_+-A_+)A_+}{A_+B_+-A^2}P_{\ell,0}^+,\\
	R_{j\ell,1}^+=\frac{(B_+-A_+)i\xi_j}{A_+B_+-A^2}P_{\ell,1}^+,
	\quad R_{N\ell,1}^+=\frac{(B_+-A_+)A_+}{A_+B_+-A^2}P_{\ell,1}^+,\\
  R_{j\ell,0}^-=-\frac{i\xi_j}{A}P_{\ell,0}^{0},\quad R_{N\ell,0}^-=-P_{\ell,0}^{0},\quad R_{j\ell,1}^-=-\frac{i\xi_j}{A}P_{\ell,1}^{0},\quad R_{N\ell,1}^-=-P_{\ell,1}^{0},\\
T_{j\ell,0}^-=-\frac{i\xi_j}{A}P_{\ell,0}^{b},\quad T_{N\ell,0}^-=P_{\ell,0}^{b},\quad T_{j\ell,1}^-=-\frac{i\xi_j}{A}P_{\ell,1}^{b},\quad T_{N\ell,1}^-=P_{\ell,1}^{b}
\end{gather*}
for $j=1,\dots,N-1$ and $\ell=1,\dots,N.$
Then substituting \eqref{eq:cof_1},  \eqref{eq:cof_2}, \eqref{eq:cof_3},
\eqref{eq:cof_4}, \eqref{eq:cof_5} and \eqref{eq:cof_6} into \eqref{eq:2.8} yields 
\begin{equation}\label{eq:2.27}
\widehat{v}_{+J}=\sum_{k=1}^{2}\widehat{v}_{Jk}^+,\quad \widehat{v}_{-J}=\sum_{k=1}^{4}\widehat{v}_{Jk}^-
\end{equation}
with
\begin{equation}\label{eq:all_symbol}
\begin{aligned}
&\widehat{v}_{J1}^{\pm}=	M_{\pm}(x_N)\sum_{\ell=1}^N\big[R_{J\ell,0}^{\pm}\widehat{h}_\ell(0)+R_{J\ell,1}^{\pm}\widehat{k}_\ell(0)\big],\\
		&\widehat{v}_{J2}^{\pm}=e^{\mp B_{\pm}x_N}\sum_{\ell=1}^N\big[S_{J\ell,-1}^{\pm}\widehat{h}_\ell(0)+S_{J\ell,0}^{\pm}\widehat{k}_\ell(0)\big],\\
		&\widehat{v}_{J3}^-=M_-(-x_N-b)\sum_{\ell=1}^N\big[T_{J\ell,0}^-\widehat{h}_\ell(0)+T_{J\ell,1}^-\widehat{k}_\ell(0)\big],\\
		&\widehat{v}_{J4}^-=e^{-B_-(x_N+b)}\sum_{\ell=1}^N\big[S_{J\ell,-1}^b\widehat{h}_\ell(0)+S_{J\ell,0}^b\widehat{k}_\ell(0)\big].
	\end{aligned}
\end{equation}
Thus the derivation of the solution formula of \eqref{eq:2.4} is complete up to taking the inverse Fourier transformation  with respect to $\xi^\prime$.

\section{Analysis of the Lopatinski determinant} 
\label{sec:det}
Recall the definition of the matrix 
\begin{equation*}
\bL= \begin{pmatrix}
L_{11} & L_{12} \\
L_{21} & L_{22} 
\end{pmatrix}
\end{equation*}
whose the entries $L_{jk}$ ($j,k=1,2$) are given by \eqref{eq:Lij}.
For simplicity, let us write 
\begin{equation}\label{def:LN}
	L_{jk}=\wt L_{jk} + N_{jk} 
\end{equation}
with
\begin{equation}\label{def:N_jk}
	\begin{aligned}
	\wt L_{jk}=&L_{jk}^++L_{jk}^{0-},\\
		N_{11}=&L_{11}^+\big(-e^{-2B_-b}-2e^{-B_-b}AM_-(-b)\big)
		-L_{11}^{b-}e^{-B_-b}+\fL_1,\\
		N_{12}=&L_{12}^+\big(-e^{-2B_-b}+2e^{-B_-b}B_-M_-(-b)\big)
		-L_{12}^{b-}e^{-B_-b}+B_-\fL_1,\\
		N_{21}=&L_{21}^+\big(-e^{-2B_-b}-2e^{-B_-b}AM_-(-b)\big)
		-L_{21}^{b-}e^{-B_-b}+\fL_2,\\
		N_{22}=&L_{22}^+\big(-e^{-2B_-b}+2e^{-B_-b}B_-M_-(-b)\big)
		-L_{22}^{b-}e^{-B_-b}+B_-\fL_2.
	\end{aligned}
\end{equation}
Then we set $\det \wt \bL=\wt L_{11}\wt L_{22}-\wt L_{12}\wt L_{21}$
for $\wt\bL=[\wt L_{jk}].$
 Moreover, it is not hard to see that
\begin{equation}\label{eq:det_whole}
	\det \bL = \wt L_{11}\wt L_{22}-\wt L_{12}\wt L_{21}+h(A_{\pm},B_{\pm})
\end{equation}
with 
\begin{equation*}
	h(A_{\pm},B_{\pm})= \wt L_{11} N_{22}+ \wt L_{22}N_{11}
	+N_{11}N_{22}-\wt L_{12}N_{21}-\wt L_{21}N_{12}-N_{12}N_{21}.
\end{equation*}
\smallbreak

The main goal of this section is to derive the lower bound of $\det \bL,$
which plays a fundamental role in our analysis. Indeed, we shall prove that
\begin{lem}\label{lem:3.1}
Let $\det \bL$ be given in \eqref{eq:det_whole}. 
Then there exists a positive $C$ depending on 
$\mu_{\pm},$ $\nu_+,$ $\varepsilon,$ 
$\gamma_{1\pm},$ $\gamma_{2+}$,  $\lambda_0$ and $b$ such that
	\begin{equation}\label{eq:3.3}
		|\det\bL|\geq C\big(|\lambda|^{1/2}+A\big)^3
	\end{equation}
	for any $(\lambda,\xi')\in\widetilde{\Gamma}_{\varepsilon,\lambda_0}
    =\Gamma_{\varepsilon,\lambda_0}\times \big(\BR^{N-1}\backslash\{0\}\big)$.
\end{lem}

\subsection{Some results on the multipliers}
In this subsection, we begin with the definition of the classes of symbols (see \cite[p.592]{ShiShi2012} for instance).
\begin{dfn}
	Let $\Lambda\subset \mathbb{C}, s\in \BR $ and $\kappa^\prime\in \BN_0^{N-1}.$
	Let $m(\lambda,\xi^\prime)$ be a smooth function defined on the domain $\widetilde{\Lambda}=\Lambda\times (\BR^{N-1}\backslash\{0\})$.
	\begin{enumerate}[label=(\roman*)]
		\item The symbol $m(\lambda,\xi^\prime)$ is called a multiplier of order $s$ with type $1$ on $\widetilde{\Lambda}$, denoted by $m\in\bM_{s,1}(\widetilde{\Lambda})$, if there exists a constant $C_{\kappa^\prime}$ such that
		$$\big|\pd_{\xi^\prime}^{\kappa^\prime}(\tau\pd_{\tau})^{\ell}
		m(\lambda,\xi^\prime)\big|\leq C_{\kappa^\prime}
		(|\lambda|^{1/2}+|\xi^\prime|)^{s-|\kappa^\prime|}$$
		for any $(\lambda,\xi^\prime)\in\widetilde{\Lambda}$ and $\ell=0,1.$
		
		\item The symbol $m(\lambda,\xi^\prime)$ is called a multiplier of order $s$ with type $2$ on $\widetilde{\Lambda}$, denoted by $m\in\bM_{s,2}(\widetilde{\Lambda})$, if there exists a constant $C_{\kappa^\prime}$ such that
		$$\big|\pd_{\xi^\prime}^{\kappa^\prime}(\tau\pd_{\tau})^{\ell}
		m(\lambda,\xi^\prime)\big|\leq C_{\kappa^\prime}
		(|\lambda|^{1/2}+|\xi^\prime|)^{s}|\xi^\prime|^{-|\kappa^\prime|}$$
		for any $(\lambda,\xi^\prime)\in\widetilde{\Lambda}$ and $\ell=0,1.$	
	\end{enumerate}
For any $m\in\bM_{s,i}(\widetilde{\Lambda})$, $i=1,2$, we denote $M(m,\Lambda) =\max_{|\kappa'|\leq N}C_{\kappa'}$.
\end{dfn}

\begin{remark}\label{rmk:4.2}
Let $s_1,s_2 \in \BR$ and $i=1,2$. It is easy to see that 
\begin{enumerate}[label=(\roman*)]
\item $m_1m_2\in\bM_{s_1+s_2,1}(\widetilde{\Lambda})$ for $m_i\in\bM_{s_i,1}(\widetilde{\Lambda})$;
\item $m_1m_2\in\bM_{s_1+s_2,2}(\widetilde{\Lambda})$ for $m_i\in\bM_{s_i,2}(\widetilde{\Lambda})$;
\item  $m_1m_2\in\bM_{s_1+s_2,2}$ for $m_1\in\bM_{s_1,1}(\widetilde{\Lambda})$ and $m_2\in\bM_{s_2,2}(\widetilde{\Lambda}).$
\end{enumerate}
\end{remark}

To study the characteristic roots in \eqref{eq:2.7}, we recall the following lemma.
\begin{lem}[{\cite[Lemma 3.1]{GS2014}}]
\label{lem:B_comp}
Let $0<\ep<\pi/2,$ $\lambda_0,$ $\mu_{\pm},$ $\nu_+>0$, $s\geq1$. Then the following assertions hold true.
\begin{enumerate}[label=(\arabic*)]	
\item For any $\lambda\in\Sigma_{\ep},$ $\alpha,\beta>0$, 
we have 
$$|\alpha\lambda+\beta|\geq \sin(\ep/2)(\alpha|\lambda|+\beta).$$

\item There exists a constant $0<\ep^\prime<\pi/2$ such that 
	$$(s\mu_++\nu_++\delta)^{-1}\lambda\in\Sigma_{\ep^\prime},\quad \forall\,\lambda\in\Gamma_{\ep,\lambda_0},$$
where the choice of $\ep^\prime$ depends on 
$s,$ $\mu_+,$ $\nu_+,$ $\gamma_{1+},$ $\gamma_{2+},$ $\lambda_0,$ $\delta_0$ and $\ep.$

\item For any $\lambda\in\Gamma_{\ep,\lambda_0}$, there exist constants $c,C>0$ depending on 
$s,$ $\mu_+,$ $\nu_+,$ $\gamma_{1+},$ $\gamma_{2+},$ $\lambda_0,$ $\delta_0$ and $\ep$ such that
	$$c(|\lambda|+|\xi^\prime|^2)\leq \big|(s\mu_++\nu_++\delta)^{-1}\lambda+|\xi^\prime|^2\big|\leq C(|\lambda|+|\xi^\prime|^2).$$
\end{enumerate}
\end{lem}

\begin{lem}\label{lem:4.3}
Let $0<\varepsilon<\pi/2$ and $\lambda_0,$ $\mu_{\pm},\nu_+>0.$ Then the following assertions hold true.	
\begin{enumerate}[label=(\arabic*)]	
\item Let the letter $M$ belong to $\{A_+,B_{\pm}, A_++B_+, \mu_+B_++\mu_-B_-\}.$
There exist positive constants $c$ and $C$ depending on $\mu_\pm,\nu_+,\gamma_{1\pm}, \gamma_{2+} ,\lambda_0$ and $\varepsilon$ such that
\begin{equation*}
c(|\lambda|^{1/2}+|\xi^\prime|)\leq |M|\leq C (|\lambda|^{1/2}+|\xi^\prime|)
\end{equation*}
for any $(\lambda,\xi^\prime)\in\widetilde{\Gamma}_{\varepsilon,\lambda_0}.$ 
Moreover, $M^s\in \bM_{s,1}(\widetilde{\Gamma}_{\varepsilon,\lambda_0})$ for any $s\in\BR.$
		
\item Let the multi-index $\alpha'\in\BN_0^{N-1}, \ell=0,1,$ and $a>0.$ 
   \begin{enumerate}
     \item There exist positive constants 
      $C_{\alpha',-}=C_{\alpha',\varepsilon,\mu_{-},\gamma_{1-}}$ and
       $c_{-}=c_{\varepsilon,\mu_{-},\gamma_{1-}}$ such that
      \begin{equation*}
          \begin{aligned}
         \big|\pd_{\xi^\prime}^{\alpha^\prime}(\tau\pd_{\tau})^{\ell}e^{-B_-a}\big|
      &\leq C_{\alpha',-}\big(|\lambda|^{1/2}+A\big)^{-|\alpha'|}e^{-c_-(|\lambda|^{1/2}+A)a}, \\
    \big|\pd_{\xi^\prime}^{\alpha^\prime}e^{-Aa}\big|
    &\leq C_{\alpha'}A^{-|\alpha'|}e^{-(1/2)Aa},\\
      \big|\pd_{\xi^\prime}^{\alpha^\prime}(\tau\pd_{\tau})^{\ell}M_-(-a)\big|
      &\leq C_{\alpha',-} \big(a \,\,\text{or}\,\,|\lambda|^{-1/2}\big) A^{-|\alpha'|}e^{-c_- Aa}
      \end{aligned}
 \end{equation*}
for any $(\lambda,\xi')\in\wt \Sigma_{\ep,\lambda_0}.$ 
 \item There exist positive constants 
 $C_{\alpha',+}=C_{\alpha',\varepsilon,\mu_{+},\nu_+,\gamma_{1+},\gamma_{2+}}$ and $c_+=c_{\varepsilon,\mu_{+},\nu_+,\gamma_{1+},\gamma_{2+}}$ such that 
 \begin{equation*}
            \begin{aligned}
     \big|\pd_{\xi^\prime}^{\alpha^\prime}(\tau\pd_{\tau})^{\ell}(e^{-B_+a})\big|
     &\leq C_{\alpha',+}(|\lambda|^{1/2}+A)^{-|\alpha'|}e^{-c_+ (|\lambda|^{1/2}+A)a},\\
   \big|\pd_{\xi^\prime}^{\alpha^\prime}(\tau\pd_{\tau})^{\ell}M_+(a)\big|
   &\leq C_{\alpha',+} a A^{-|\alpha'|}e^{- c_+(|\lambda|^{1/2}+A)a}
            \end{aligned}
        \end{equation*}
for any $(\lambda,\xi')\in\wt \Gamma_{\ep,\lambda_0}.$ 
    \end{enumerate}

\item Let $i,j=1,2$ and let the letter $k$ belong to $\{+,0-,b-\}.$ Recalling the notations in \eqref{eq:2.15} and \eqref{eq:det_whole}, we have 
\begin{equation*}
   L_{i1}^k\in\bM_{1,2}(\widetilde{\Gamma}_{\varepsilon,\lambda_0}),\quad 
  L_{i2}^k\in\bM_{2,2}(\widetilde{\Gamma}_{\varepsilon,\lambda_0})
  \quad \text{and}\quad 
  \det \bL\in\bM_{3,2}(\widetilde{\Gamma}_{\varepsilon,\lambda_0}).
\end{equation*}
		 Moreover, $M(L_{ij}^k,\Gamma_{\varepsilon,\lambda_0})$ and $M(\det\bL,\Gamma_{\varepsilon,\lambda_0})$ depend solely on $\varepsilon,\mu_{\pm},\nu_+,\gamma_{1\pm},\gamma_{2+},\lambda_0$.
	\end{enumerate}	
\end{lem}
\begin{proof}
The assertion (1) is standard in view of Lemma \ref{lem:B_comp} above, while the bounds in assertion (2) have been proved by \cite[Lemma 5.3]{ShiShi2012} and \cite[p.468]{ES2013} respectively. Finally, using Remark \ref{rmk:4.2} and the assertions (1)-(2), it is not hard to prove the last assertion and we omit the details.
\end{proof}

\subsection{Proof of Lemma \ref{lem:3.1}}
\label{subsec:lemma_3.1}
The proof of Lemma \ref{lem:3.1} will be divided into four steps. 
\medskip

{\bf Step 1.} We shall verify the bound \eqref{eq:3.3} when $(\lambda,\xi')$ locates in the area
$$\CD_1 = \{(\lambda,\xi')\in
\widetilde{\Gamma}_{\varepsilon,\lambda_0}:
|\xi'|\leq r_{1}\}$$
for some sufficiently small positive constant $r_{1}<1.$ 
\medskip

Firstly, we claim that there exists some (large enough) positive constant
$d_1$ such that 
\begin{equation}\label{eq:c11}
	|\det\bL| \gtrsim (|\lambda|^{1/2}+A)^3  \quad 
	\text{for}\,\,\, 0<|\xi'|\leq 1
	\,\,\,\text{and}\,\,\, |\lambda|\geq d_1.
\end{equation}
To prove \eqref{eq:c11},  we rewrite $\det \bL$ as follows:
\begin{equation}\label{eq:rdet}
	\det\bL =H_1(A_{\pm},B_{\pm})+h_1(A_{\pm},B_{\pm})
\end{equation}
with
\begin{equation}\label{eq:h(A,B)}
	\begin{aligned}
		H_1(A_{\pm},B_{\pm})=&(L_{11}^++L_{11}^{0-})
		\big( L_{22}^{0-}+B_-M_-(-b)L_{22}^{b-} \big),\\
		h_1(A_{\pm},B_{\pm})=
		&L_{11}\Big( L_{22}^+(1-e^{-2B_-b}+2e^{-B_-b}B_-M_-(-b))
		-L_{22}^{b-}e^{-B_-b}\\
		&\quad\quad-(AL_{21}^++L_{22}^+)(A+B_-)B_-M_-^2(-b)
		-L_{21}^{b-}AB_-M_-(-b) \Big)\\
		&+\big(L_{22}^{0-}+B_-M_-(-b)L_{22}^{b-}\big)
		N_{11}-L_{12}L_{21}.
	\end{aligned}
\end{equation}
Above, the formula of $N_{11}$ is given by \eqref{def:N_jk}.
\smallbreak

Notice the fact that
\begin{equation*}
	B_-M_-(-b)=e^{-B_-b}+AM_-(-b)-e^{-Ab}. 
\end{equation*}
Then we use \eqref{eq:2.15} to obtain that
\begin{equation*}
	\begin{aligned}
			L_{11}^++L_{11}^{0-}
		=&-(\mu_{+}B_++\mu_{-}B_-)
		-\mu_{+}\frac{(B_+-A_+)A^2}{A_+B_+-A^2}-\mu_{-}A,\\
		L_{22}^{0-}+B_-M_-(-b)L_{22}^{b-}
		=&\,\, L_{22}^{0-}-e^{-Ab} L_{22}^{b-}
		+ \big( e^{-B_-b}+AM_-(-b) \big)  L_{22}^{b-}\\
		=&-\mu_-(1+e^{-2Ab})B_-^2+R_1(A_\pm,B_\pm),
	\end{aligned}
\end{equation*}
where $R_1(A_\pm,B_\pm)$ is given by
\begin{equation*}
	\begin{aligned}
		R_1(A_\pm,B_\pm)=&-\mu_{-}AB_-+2\mu_-A^2B_-M_-(-b)e^{-Ab}\\
		&-\mu_-AB_-e^{-2Ab}
		+\big( e^{-B_-b}+AM_-(-b) \big)  L_{22}^{b-}.
	\end{aligned}
\end{equation*}
Thus we see that 
\begin{equation}\label{eq:rwdet}
	\begin{aligned}
		H_1(A_{\pm},B_{\pm})
		=&\mu_{-}\big(1+e^{-2Ab}\big)B_-^2
		\big(\mu_{+}B_++\mu_{-}B_-\big)+R_2(A_\pm,B_\pm)
	\end{aligned}
\end{equation}
with
\begin{equation*}
	\begin{aligned}
		R_2(A_\pm,B_\pm)
		=&\Big( \mu_{+}\frac{(B_+-A_+)A^2}{A_+B_+-A^2}+\mu_{-}A \Big)
		\mu_-(1+e^{-2Ab})B_-^2\\
		&+(L_{11}^++L_{11}^{0-}) R_1(A_\pm,B_\pm).
	\end{aligned}
\end{equation*}
Moreover, Lemma \ref{lem:4.3} yields
\begin{equation*}
	\big| \mu_{-}\big(1+e^{-2Ab}\big)B_-^2
	\big(\mu_{+}B_++\mu_{-}B_-\big) \big| \gtrsim (|\lambda|^{1/2}+A)^3.
\end{equation*}
\smallbreak 

To find the upper bound of the remainder term $R_2(A_\pm,B_\pm),$ Lemma \ref{lem:4.3} yields 
\begin{equation}\label{eq:est}
	\begin{aligned}
		&|e^{-aB_-b}| \lesssim  e^{-ca(|\lambda|^{1/2}+A)b}
		\lesssim  e^{-ca|\lambda|^{1/2}b}   \quad (a>0),\\ 
		&|M_-(-b)| \lesssim \frac{1}{|\lambda|^{1/2}}e^{-cAb},\\
		&|A_+M_-(-b)|+|B_\pm M_-(-b)|\lesssim \frac{|\lambda|^{1/2}+A}{|\lambda|^{1/2}}e^{-cAb}.
	\end{aligned}
\end{equation}
for some positive constant $c=c_{\varepsilon,\mu_-,\gamma_{1-}}.$ 
As in \cite[Section 4]{KSS2014}, we rewrite $L^+_{jk}$ in \eqref{eq:2.15} by
\begin{equation}\label{eq:Lkks}
	\begin{aligned}
		&L_{11}^+=-\frac{\mu_{+}(\mu_{+}+\nu_++\delta)}{2\mu_{+}+\nu_++\delta}A_+P,\quad L_{12}^+=-\mu_{+}A^2(2-\frac{\mu_{+}+\nu_++\delta}{2\mu_{+}+\nu_++\delta}P),\\
		&L_{21}^+=-\Big(\frac{2\mu_{+}(\nu_++\delta)}{2\mu_{+}+\nu_++\delta}\frac{A_+}{B_++A_+}-\frac{\mu_{+}(\nu_+-\mu_{+}+\delta)}{2\mu_{+}+\nu_++\delta}\Big)AP,\\
		&L_{22}^+=-\frac{\mu_{+}(\mu_{+}+\nu_++\delta)}{2\mu_{+}+\nu_++\delta}AB_+P
	\end{aligned}
\end{equation}
with $$P(\lambda,\xi')=\frac{A_+B_++A^2}{\gamma_{1+}(2\mu_{+}+\nu_++\delta)\lambda+A^2}\in\bM_{0,1}.$$ 

For $s\in \{\nu_+,\mu_++\nu_+,\nu_+-\mu_{+}\}$, we have
\begin{equation*}
	\Big|\frac{s+\delta}{2\mu_++\nu_++\delta} \Big|\lesssim 1.
\end{equation*}
Then we see from \eqref{eq:2.15}, \eqref{eq:est} and \eqref{eq:Lkks} that
\begin{equation}\label{est:Ljk_1}
	\begin{aligned}
		& |L_{22}^{0-}| + |L_{22}^{b-}|\lesssim  |\lambda|,\\
		& |L_{11}^+|+|L_{22}^+| +|L_{11}^{0-}|  +|L_{12}^{0-}|
		+|L_{21}^{0-}|+ |L_{21}^{b-}|+|\fL_2| \lesssim |\lambda|^{1/2},\\
		& |L_{12}^+| +|L_{21}^+|+|L_{12}^{b-}|+|\fL_1| \lesssim 1,\\
		& |L_{11}^{b-}|   \lesssim |\lambda|^{-1/2}
	\end{aligned}
\end{equation}
for any $0<|\xi'|\leq 1 \leq  |\lambda|,$
which yield
\begin{equation*}
\begin{aligned}
& |L_{11}|+|L_{12}|+|L_{21}|+|R_1(A_\pm,B_\pm)|\lesssim |\lambda|^{1/2}, \\
&   |L_{22}| + |R_2(A_\pm,B_\pm)| \lesssim |\lambda|.
\end{aligned}
\end{equation*}
Furthermore, recall the definition of $h_1(A_{\pm},B_{\pm})$ in \eqref{eq:h(A,B)}, and we use \eqref{eq:est} and \eqref{est:Ljk_1} to obtain that
\begin{equation*}
	|h_1(A_{\pm},B_{\pm})| \lesssim |\lambda|+ |\lambda| |N_{11}| \lesssim  |\lambda|.
\end{equation*}
Thus summing up the estimates above, we can conclude that
\begin{equation}\label{eq:c11-1}
|\det\bL |\geq |H_1(A_{\pm},B_{\pm})|-C_1 |\lambda|\geq C_2^{-1}(|\lambda|^{1/2}+A)^3-C_2 |\lambda|
\end{equation}
for any $0<|\xi'|\leq 1 \leq  |\lambda|.$
Then we easily see \eqref{eq:c11} follows from 
\eqref{eq:c11-1} by choosing the constant $d_1$ large enough. 
\medskip

Now, we fix the constant $d_1$ chosen in \eqref{eq:c11}.
Let us introduce the function 
$$F(A,\lambda)=\det\bL$$ 
for any $\lambda\in \CK=\{z\in\Gamma_{\ep,\lambda_0}:\,|z|\leq d_1\}.$
It is easy to see that 
\begin{gather*}
	\fL_1|_{A=0}=L_{12}|_{A=0}=0, \quad  
	 \fL_2|_{A=0}=(\mu_{-}\gamma_{1-})^{1/2}\lambda^{1/2}(e^{-b\sqrt{\gamma_{1-}(\mu_{-})^{-1}\lambda}}-1),\\
	L_{11}|_{A=0}=-(\mu_+\gamma_{1+})^{1/2} \lambda^{1/2}
		\big( 1-e^{-2b\sqrt{ \gamma_{1-}(\mu_-)^{-1}\lambda}}  \big)
		-(\mu_-  \gamma_{1-})^{1/2} \lambda^{1/2}
		\big(1+e^{-2b\sqrt{ \gamma_{1-}(\mu_-)^{-1}\lambda}}\big),\\
	L_{22}|_{A=0} = -2\gamma_{1-}\lambda.
\end{gather*}
Thus we have 
\begin{equation*}
	\begin{aligned}
		F(0,\lambda)& = (L_{11} L_{22})|_{A=0} \\
		&=2  \gamma_{1-} \lambda^{3/2} \Big( (\mu_+\gamma_{1+})^{1/2} 
		\big( 1-e^{-2b\sqrt{ \gamma_{1-}(\mu_-)^{-1}\lambda}}  \big)
		+(\mu_-  \gamma_{1-})^{1/2} 
		\big(1+e^{-2b\sqrt{ \gamma_{1-}(\mu_-)^{-1}\lambda}} \big) \Big).
	\end{aligned}
\end{equation*}
 As $|\sqrt{\lambda}|\geq \lambda_0^{1/2},$ we have 
\begin{equation*}
\Big| e^{-2b\sqrt{ \gamma_{1-}(\mu_-)^{-1}\lambda}} \Big| 
\leq c' \quad \text{and}\quad
|F(0,\lambda)| \geq 2c''
\end{equation*}
with 
\begin{equation*}
\begin{aligned}
c' &=\exp \Big(  -2b
	\big(\gamma_{1-}(\mu_-)^{-1}\lambda_0\big)^{1/2}\Big),
	\\
c''&=\gamma_{1-} \lambda_0^{3/2}   (1-c') 
	\big( (\mu_+\gamma_{1+})^{1/2}
	+(\mu_-  \gamma_{1-})^{1/2} \big) .	
\end{aligned}
\end{equation*} 
By the continuity of $F(A,\lambda)$ and the compactness of the set $\CK,$
there exists a small positive constant $r_1$ such that 

\begin{equation}\label{es:C1-2}
	|\det \bL| =|F(A,\lambda)|\geq  c''
	\geq c''\big( d_{1}^{1/2}+r_1 \big)^{-3}\big(|\lambda|^{1/2}+A\big)^3
\end{equation}
for any $(\lambda,A) \in [0, r_1]\times \CK.$ 
\smallbreak

Therefore, the bound \eqref{eq:3.3} holds for any $(\lambda,\xi')\in \CD_1$ thanks to \eqref{eq:c11} and \eqref{es:C1-2}.
\medskip

{\bf Case 2.}
For the parameter $r_1$ given in Case 1, we introduce the set 
$$\CD_2 = \{(\lambda,\xi')\in
\widetilde{\Gamma}_{\varepsilon,\lambda_0}:
|\xi'|\geq  r_{1},\,\,\, R_2|\xi'|\leq |\lambda|^{1/2}\}.$$
Here, $R_2$ is a fixed very large positive number. 
In this case,  we know from \cite[Section 5]{KSS2014} that
\begin{equation}\label{es:C2-1}
	\begin{aligned}
		&A_+=(\mu_++\nu_++\delta)^{-1/2}(\gamma_{1+}\lambda)^{1/2}\big(1+O(R_2^{-1})\big),\\ 
		&B_{\pm}=(\mu_{\pm})^{-1/2}(\gamma_{1\pm}\lambda)^{1/2}\big(1+O(R_2^{-1})\big),\\
		&L_{11}^+=-(\mu_+\gamma_{1+})^{1/2}\lambda^{1/2}\big(1+O(R_2^{-1})\big),
		\quad     L_{12}^+=-A^2 O(1),\\
		&L_{21}^+=-A\,O(1),\quad 
		L_{22}^+=-\lambda O(R_2^{-1}),\\
		&L_{11}^{0-}=L_{21}^{0-}=-(\mu_{-}\gamma_{1-})^{1/2}\lambda^{1/2}
		\big(1+O(R_2^{-1})\big),\quad   
		L_{12}^{0-}=-\lambda O(R_2^{-1}),\\
		& L_{22}^{0-}=-\gamma_{1-}\lambda\big(1+O(R_2^{-1})\big).
	\end{aligned}
\end{equation}
On the other hand, according to \eqref{eq:est} and \eqref{es:C2-1}, we have
\begin{equation}\label{es:C2-2}
\begin{aligned}
 L_{11}^{b-}&=(\mu_{-}\gamma_{1-}\lambda)^{1/2}e^{-B_-b}\big(1+O(R_2^{-1})\big)
+\lambda^{1/2}e^{-cAb}O(R_2^{-2}), \\
L_{12}^{b-}&=\lambda  O(R_2^{-1}),\\
L_{21}^{b-}&=-(\mu_{-}\gamma_{1-}\lambda)^{1/2}e^{-Ab}
+\lambda^{1/2}e^{-cAb} O(R_2^{-1}),\\
L_{22}^{b-}&=\gamma_{1-}\lambda e^{-Ab}\big(1+O(R_2^{-1})\big)+\lambda e^{-cAb} O(R_2^{-2}).
	\end{aligned}
\end{equation}

Next, we use \eqref{eq:rdet}, \eqref{eq:rwdet} and \eqref{es:C2-1} to obtain 
\begin{equation*}
	\begin{aligned}
		|\lambda|^{-3/2} |\det\bL| & \geq |\lambda|^{-3/2} |H_1(A_{\pm},B_{\pm})|
		- |\lambda|^{-3/2}|h_1(A_{\pm},B_{\pm})| \\
		&\geq 2c-|\lambda|^{-3/2} \big( |h_1(A_{\pm},B_{\pm})| +|R_2(A_{\pm},B_{\pm})|  \big)
	\end{aligned}
\end{equation*}
for some constant $c>0.$
In fact, we can prove that there exists a constant $C_2$ such that
\begin{equation}\label{clm:c2}
	|\lambda|^{-3/2} \big( |h_1(A_{\pm},B_{\pm})| +|R_2(A_{\pm},B_{\pm})|  \big)
	\leq C_2 \big(R_2^{-1}+e^{-C_2^{-1} R_2}\big).
\end{equation}
Then choosing $R_2$ large enough, we see that 
\begin{equation*}
	|\det\bL|  \geq c  |\lambda|^{3/2} 
	\geq c ( 1+R_2^{-1} )^{-3}\big(|\lambda|^{1/2}+A\big)^3 
\end{equation*}
for any $(\lambda,\xi')\in \CD_2.$
\medskip

To prove \eqref{clm:c2},  we first note from 
\eqref{eq:est}, \eqref{es:C2-1} and \eqref{es:C2-2} that
\begin{equation}\label{es:C2-3}
	\begin{aligned}
		&|\lambda|^{-1/2} |\fL_1|+|\lambda|^{-1}|L_{12}|=O(R_2^{-1}), \\
		&|\lambda|^{-1/2} (|\fL_2|+|L_{11}|+|L_{21}|)+|\lambda|^{-1}|L_{22}| =O(1).
	\end{aligned}
\end{equation}
Since $|\lambda|^{1/2} \geq r_1 R_2$ for any $(\lambda,\xi')\in\CD_2,$  we infer from  
\eqref{eq:est}, \eqref{es:C2-1}, \eqref{es:C2-2} and \eqref{es:C2-3} that
\begin{equation*}
	\begin{aligned}
		& |\lambda|^{-1} |R_1(A_{\pm},B_{\pm})|
		+|\lambda|^{-3/2} |R_2(A_{\pm},B_{\pm})|
		\lesssim  R_2^{-1}+e^{-c_2r_1 R_2},\\
				&	 |\lambda|^{-3/2} |h_1(A_{\pm},B_{\pm})|
	\lesssim { R_2^{-1}+e^{-c|\lambda|^{1/2}b} 
		+|\lambda|^{-1/2}|N_{11}| }
	\lesssim R_2^{-1}+e^{-c_2 r_1 R_2}
	\end{aligned}
\end{equation*}
for some constant $c_2>0.$ Then the proof of \eqref{clm:c2} is complete.
\medskip

{\bf Case 3.} For the fixed number $R_2$ given in {\bf Case 2}, let us consider the domain
\begin{equation*}
	\CD_3 = \{(\lambda,\xi')\in\widetilde{\Gamma}_{\varepsilon,\lambda_0}: |\xi'|\geq R_3,\,\,R_2|\xi'|\geq |\lambda|^{1/2} \}
\end{equation*}
where $R_3$ is a sufficiently large positive number.
Moreover, by Lemma \ref{lem:4.3}, there exists a positive constant $c=c_{\varepsilon,\gamma_{1-},\mu_-}$ such that
\begin{equation}\label{es:exp_c3}
	\begin{aligned}
		|A_+|+ |B_\pm|&\lesssim |\lambda|^{1/2} + A \lesssim R_2A,\\
		|e^{-aB_-b}|& \lesssim  e^{-ca(|\lambda|^{1/2}+A)b}
		\lesssim  e^{-caA b},\\
		|M_-(-b)|& \lesssim b e^{-cb A},\\ 
		|AM_-(-b)|&\lesssim  bAe^{-cAb} \lesssim e^{-cb A/2},\\
		|B_-M_-(-b)|&\lesssim  R_2bA e^{-cb A}
		\lesssim  R_2 e^{-cb A/2}
	\end{aligned}
\end{equation}
for any $(\lambda,\xi')\in \CD_3,$ $a>0$ and $R_2 \geq 1.$
\smallbreak 

To derive the lower bound of $|\det \bL|,$ we first notice from \eqref{eq:2.15}, \eqref{def:N_jk} and \eqref{es:exp_c3} that
\begin{equation}\label{es:L_3}
	\begin{aligned}
		& |L_{11}^+| + |L_{11}^{0-}| +|L_{21}^{0-}|
		+ { |\wt L_{11}|}+ |\wt L_{21}|\lesssim R_2 A, \quad
		|L_{12}^+| \lesssim A^2,\quad
		|L_{21}^+| \lesssim A, \\
		& |L_{22}^+| + |L_{12}^{0-}|+|\wt L_{12}|  \lesssim R_2 A^2, \qquad
		|L_{22}^{0-}|+|\wt L_{22}| \lesssim  R_2^2 A^2,\\ 
	 &	|L_{11}^{b-}|+|L_{12}^{b-}| + |L_{21}^{b-}| \lesssim R_2 e^{-d_3 A}, \qquad
		 |L_{22}^{b-}|\lesssim R_2^2 e^{-d_3 A} 
	\end{aligned}
\end{equation}
for some constant $d_3$ satisfying $0<2 d_3<b \min\{1,c\}.$
Then \eqref{def:N_jk}, \eqref{es:exp_c3} and \eqref{es:L_3} yield
\begin{equation}\label{es:N_3}
	\begin{aligned}
		|\fL_1|+|N_{11}|+|N_{12}| \lesssim &  
		\big( |L_{11}^+|+|L_{12}^+|
		+|L_{11}^{b-}|+| L_{12}^{b-} |\big) { R_2^2} e^{-cbA/2} 
		\lesssim  {  R_2^3 }  e^{-cbA/4},\\
		|\fL_2|+|N_{21}|+|N_{22}| \lesssim & 
		\big(|L_{21}^+|+|L_{22}^+ | 
		+|L_{21}^{b-}|+ |L_{22}^{b-}|\big) { R_2^2 } e^{-cbA/2}
		\lesssim  {  R_2^4 }  e^{-cbA/4}.
	\end{aligned}
\end{equation}
Combining the estimates \eqref{es:L_3} and \eqref{es:N_3},  we see that
\begin{equation}\label{es:C3-1}
	|h(A_{\pm},B_{\pm})|\lesssim  { R_2^7 } A^2  e^{-cbA/4}.
\end{equation}
By \cite[Lemma 5.1]{KSS2014}, there exists a constant $c_3>0$ such that
\begin{equation*}
|\wt L_{11}\wt L_{22}-\wt L_{12}\wt L_{21}|
	\geq  2c_3 (|\lambda|^{1/2}+A)^3,
\end{equation*}
which together with \eqref{eq:det_whole} and \eqref{es:C3-1} yields
\begin{equation*}
	|\det \bL|\geq 2c_3 (|\lambda|^{1/2}+A)^3-C R_2^7 A^2  e^{-cbA/4} 
	\geq c_3 (|\lambda|^{1/2}+A)^3
\end{equation*}
for large enough $R_{3}$.
\medskip

{\bf Case 4.}
In this step, we shall prove $\det \bL\neq 0$ by the contradiction argument provided with $(\lambda,A)\in \CD_4$ for
\begin{equation*}
	\begin{aligned}
		\CD_4 &= \{(\lambda,\xi')\in\widetilde{\Gamma}_{\varepsilon,\lambda_0}: r_{1}\leq |\xi'|\leq R_3,
		\,\,R_2|\xi'|\geq  |\lambda|^{1/2} \}\\
		&  \,\, \subset \{(\lambda,\xi')\in\widetilde{\Gamma}_{\varepsilon,\lambda_0}:
		r_{1}\leq |\xi'|\leq R_3,
		\,\lambda_0^{1/2}\leq|\lambda|^{1/2}\leq R_2R_3\}.  
	\end{aligned}
\end{equation*}
Here, the constants $R_2$ and $R_3$ are fixed by {\bf Case 2} and {\bf Case 3} respectively.
\smallbreak 

Suppose that $\det\bL=0$ for some $(\lambda_{\star},\xi'_{\star})\in \CD_4$. Thus the problem \eqref{eq:2.4}-\eqref{eq:2.5} admits some non-trivial solution 
$$(\bu_{\pm},\fq)(x_N)=(\widehat{\bv}_{\pm},\widehat{\fp}_-)(\xi',x_N)$$
for $\widehat{h}_J(0)=\widehat{k}_J(0)=0$ ($J=1,\dots,N$). 
For simplicity, we still write $(\lambda,\xi')=(\lambda_{\star},\xi'_{\star})$ in what follows.
Therefore, $(\bu_{\pm},\fq)$ satisfies the equations
\begin{equation}\label{eq:ho}
	\begin{aligned}
		&\gamma_{1+}\lambda u_{+j}-\sum_{\ell=1}^{N-1}\mu_+i\xi_\ell(i\xi_ju_{+\ell}+i\xi_\ell u_{+j})-\mu_+\pd_N(i\xi_ju_{+N}+\pd_Nu_{+j}) \\
		&\quad-(\nu_+-\mu_++\delta)i\xi_j(i\xi'\cdot \bu_+'+\pd_Nu_{+N})=0
		& \text{for}\,\,x_N>0,\\
		&\gamma_{1+}\lambda u_{+N}-\sum_{\ell=1}^{N-1}\mu_+i\xi_\ell(\pd_Nu_{+\ell}+i\xi_\ell u_{+N})-2\mu_+\pd_N^2u_{+N}\\
		&\quad-(\nu_+-\mu_++\delta)\pd_N(i\xi'\cdot \bu_+'+\pd_Nu_{+N})=0
		& \text{for}\,\,x_N>0,\\
		&\gamma_{1-}\lambda u_{-j}-\sum_{\ell=1}^{N-1}\mu_-i\xi_\ell(i\xi_ju_{-\ell}+i\xi_\ell u_{-j})\\
		&\quad-\mu_-\pd_N(i\xi_ju_{-N}+\pd_Nu_{-j})+i\xi_j\fq=0
		&  \text{for}\,\,-b<x_N<0,\\
		&\gamma_{1-}\lambda u_{-N}-\sum_{\ell=1}^{N-1}\mu_{-}i\xi_\ell(\pd_Nu_{-\ell}+i\xi_\ell u_{-N})\\
		&\quad-2\mu_{-}\pd_N^2u_{-N}+\pd_N\fq=0
		&  \text{for}\,\, -b<x_N<0,\\
		&i\xi'\cdot \bu_-'+\pd_Nu_{-N}=0
		&  \text{for}\,\,-b<x_N<0
	\end{aligned}
\end{equation}
with the boundary conditions
\begin{equation}\label{bc:unique}
	\begin{aligned}
		&\mu_+(\pd_Nu_{+j}+i\xi_ju_{+N})\big|_{x_N=0}
		=\mu_{-}(\pd_Nu_{-j}+i\xi_ju_{-N})\big|_{x_N=0},\\
		&2\mu_+\pd_Nu_{+N}+(\nu_+-\mu_++\delta)
		(i\xi'\cdot \bu_+'+\pd_Nu_{+N})\big|_{x_N=0}
		=(2\mu_{-}\pd_Nu_{-N}-\fq)\big|_{x_N=0},\\
		& u_{+J}\big|_{x_N=0}=u_{-J}\big|_{x_N=0},\quad 
		u_{-J}\big|_{x_N=-b}=0.
	\end{aligned}
\end{equation}

Let us set  $\|f_{\pm}\|_{\pm}=(f,f)_\pm^{1/2}$ with
$$(f_{+},g_{+})_{+}=\int_{0}^{\infty} f_+(t)\overline{g_+(t)}\,dt
\quad \text{and}\quad
(f_{-},g_{-})_{-}=\int^{0}_{-b} f_-(t)\overline{g_-(t)}\,dt.$$
Then multiplying the equations in \eqref{eq:ho} by $\overline{u_{\pm J}}$ and using integration by parts and interface conditions in \eqref{bc:unique},  we obtain that

\begin{equation}\label{eq:l2}\allowdisplaybreaks 
	\begin{aligned}
		0=&\lambda\big(\gamma_{1+}\sum_{J=1}^N\|u_{+J}\|_+^2
		+\gamma_{1-}\sum_{J=1}^N\|u_{-J}\|_-^2\big)\\
		&+\mu_+\Big( \sum_{j,\ell=1}^{N-1}\|i\xi_\ell u_{+j}\|_+^2
		+\|i\xi'\cdot \bu_+'\|_+^2+\sum_{j=1}^{N-1}\|\pd_Nu_{+j}\|_+^2
		+\sum_{\ell=1}^{N-1}\|i\xi_{\ell}u_{+N}\|_+^2
		+2\|\pd_Nu_{+N}\|_+^2\\
		&\qquad +\sum_{j=1}^{N-1}(i\xi_ju_{+N},\pd_Nu_{+j})_+		
		+\sum_{\ell=1}^{N-1}(\pd_Nu_{+\ell},i\xi_{\ell}u_{+N})_+\Big)\\
		&+(\nu_+-\mu_++\delta)\Big( \|i\xi'\cdot \bu_+'\|_+^2+\|\pd_Nu_{+N}\|_+^2
		+(\pd_Nu_{+N},i\xi'\cdot \bu_+')_+
		+(i\xi'\cdot \bu_+',\pd_Nu_{+N})_+\Big)\\
		&+\mu_-\Big( \sum_{j,\ell=1}^{N-1}\|i\xi_\ell u_{-j}\|_-^2+\|i\xi'\cdot \bu_-'\|_-^2+\sum_{j=1}^{N-1}\|\pd_Nu_{-j}\|_-^2  
		+\sum_{\ell=1}^{N-1}\|i\xi_{\ell}u_{-N}\|_-^2
		+2\|\pd_N u_{-N}\|_-^2\\
		&\qquad \qquad+\sum_{\ell=1}^{N-1}(\pd_Nu_{-\ell},i\xi_{\ell}u_{-N})_-
		+\sum_{j=1}^{N-1}(i\xi_ju_{-N},\pd_Nu_{-j})_-	\Big).
	\end{aligned}
\end{equation}
Thus it is not hard to find from \eqref{eq:l2} that
\begin{equation}\label{eq:l22}
	\begin{aligned}
		0=&\lambda(\gamma_{1+}\|\bu_+\|_+^2+\gamma_{1-}\|\bu_-\|_-^2)\\
		&+\mu_+\Big(\sum_{j,\ell=1}^{N-1}\|i\xi_\ell u_{+j}\|_+^2+\|i\xi'\cdot \bu_+'\|_+^2
		+\sum_{j=1}^{N-1}\|\pd_Nu_{+j}+i\xi_ju_{+N}\|_+^2+2\|\pd_Nu_{+N}\|_+^2\Big)\\
		&+(\nu_+-\mu_++\delta)\|i\xi'\cdot \bu_+'+\pd_Nu_{+N}\|_+^2\\
		&+\mu_-\Big(\sum_{j,\ell=1}^{N-1}\|i\xi_\ell u_{-j}\|_-^2
		+\|i\xi'\cdot \bu_-'\|_-^2
		+\sum_{j=1}^{N-1}\|\pd_Nu_{-j}+i\xi_ju_{-N}\|_-^2
		+2\|\pd_Nu_{-N}\|_-^2\Big).
	\end{aligned}
\end{equation}

Next, taking the imaginary part of \eqref{eq:l22}, we see that 
\begin{equation}\label{eq:IM}
	\begin{aligned}
	 0=& \Im\lambda(\gamma_{1+}\|\bu_+\|_+^2+\gamma_{1-}\|\bu_-\|_-^2)
		 +\Im \delta 
		 \|i\xi'\cdot \bu_+'+\pd_Nu_{+N}\|_+^2.
	\end{aligned}
\end{equation}
On the other hand, taking the real part of \eqref{eq:l22} and using the inequality
\begin{equation*}
	\begin{aligned}
		\sum_{j,\ell=1}^{N-1}\|i\xi_ju_{+\ell}\|_+^2
		+\|i\xi'\cdot \bu_+'\|_+^2+2\|\pd_Nu_{+N}\|_+^2 
		&\geq 2(\|i\xi'\cdot \bu_+'\|_+^2+\|\pd_Nu_{+N}\|_+^2)\\
		&\geq \|i\xi'\cdot \bu_+'+\pd_Nu_{+N}\|_+^2,
	\end{aligned}
\end{equation*}
we get that
\begin{equation*}
	\begin{aligned}
		0\geq &\Re\lambda(\gamma_{1+}\|\bu_+\|_+^2
		+\gamma_{1-}\|\bu_-\|_-^2)
		+\mu_+\sum_{j=1}^{N-1}\|\pd_Nu_{+j}+i\xi_ju_{+N}\|_+^2\\
		&+  ( \nu_++\Re \delta ) 
		\|i\xi'\cdot \bu_+'+\pd_Nu_{+N}\|_+^2\\
		&+\mu_-\Big(\sum_{j,\ell=1}^{N-1}\|i\xi_\ell u_{-j}\|_-^2+\|i\xi'\cdot \bu_-'\|_-^2+\sum_{j=1}^{N-1}\|\pd_Nu_{-j}+i\xi_ju_{-N}\|_-^2+2\|\pd_Nu_{-N}\|_-^2\Big).
	\end{aligned}
\end{equation*}
In particular, we have from the inequality above that
\begin{equation}\label{eq:con}
	\begin{aligned}
		0\geq &\Re\lambda(\gamma_{1+}\|\bu_+\|_+^2
		+\gamma_{1-}\|\bu_-\|_-^2)+ ( \nu_++\Re \delta ) 
		\|i\xi'\cdot \bu_+'+\pd_Nu_{+N}\|_+^2.
	\end{aligned}
\end{equation}

In fact, according to \eqref{eq:IM} and \eqref{eq:con}, we can prove that
\begin{equation}\label{cl:u=0}
\gamma_{1+}\|\bu_+\|_+^2
		+\gamma_{1-}\|\bu_-\|_-^2=0
\end{equation}
which immediately gives $\bu_{\pm}=0.$ Then we have $\fq=0$ from \eqref{eq:ho}.
This is a contradiction to  our choice of $(\bu_{\pm},\fq)$. 
\medskip

The proof of \eqref{cl:u=0} will be divided into three cases thanks to Assumption \ref{ass:dl}.
\begin{enumerate}
\item[(C1).] 
In this situation,  we have $\delta=\gamma_{1+}\gamma_{2+}\lambda^{-1}.$
\begin{itemize}

\item If $\Im \lambda\not=0,$ then by \eqref{eq:IM} we have
\begin{equation}\label{c4:1}
	\gamma_{1+}\|\bu_+\|_+^2+\gamma_{1-}\|\bu_-\|_-^2
	=\gamma_{1+}\gamma_{2+}|\lambda|^{-2}\|i\xi'\cdot \bu_+'+\pd_Nu_{+N}\|_+^2,
\end{equation}
which combined with \eqref{eq:con} furnishes that
\begin{equation*}
	\begin{aligned}
		0\geq& \nu_+|\lambda|^{-2}(2\gamma_{1+}\gamma_{2+}\nu_+^{-1}\Re\lambda+|\lambda|^2)\|i\xi'\cdot \bu_+'+\pd_Nu_{+N}\|_+^2\\
		=&\nu_+|\lambda|^{-2}\big((\Re\lambda+\gamma_{1+}\gamma_{2+}\nu_+^{-1})^2+(\Im\lambda)^2-(\gamma_{1+}\gamma_{2+}\nu_+^{-1})^2\big)
		\|i\xi'\cdot \bu_+'+\pd_Nu_{+N}\|_+^2,
	\end{aligned}
\end{equation*}
as $(\Re\lambda+\gamma_{1+}\gamma_{2+}\nu_+^{-1})^2+(\Im\lambda)^2-(\gamma_{1+}\gamma_{2+}\nu_+^{-1})^2>0$ for $\lambda \in K_\varepsilon.$ 
Then  we have 
\begin{equation*}
0=\|i\xi'\cdot \bu_+'+\pd_Nu_{+N}\|_+^2.
\end{equation*}
Thus \eqref{cl:u=0} holds in view of \eqref{c4:1}.
\item For the case where $\Im \lambda=0$ and $\lambda\in \Sigma_{\ep, \lambda_0}$ , we have $\lambda>0$ and $\delta=\gamma_{1+}\gamma_{2+}\lambda^{-1}>0.$ 
Thus we can conclude \eqref{cl:u=0} from \eqref{eq:con}.  

\end{itemize}

\item[(C2).] Suppose that $\delta\in\Sigma_{\varepsilon}$ with $\Re\delta<0$, $\lambda\in\mathbb{C}$ with $|\lambda|\geq \lambda_0$ and $\Re\lambda\geq|\Re\delta/\Im\delta||\Im\lambda|$. 

\begin{itemize}
\item If $\Im\lambda=0$,  we have $\|i\xi'\cdot \bu_+'+\pd_Nu_{+N}\|_+^2=0$ from \eqref{eq:IM}.
Then \eqref{eq:con} yields
\begin{equation*}
	0\geq \Re\lambda(\gamma_{1+}\|\bu_+\|_+^2+\gamma_{1-}\|\bu_-\|_-^2).
\end{equation*}
As $|\lambda|\geq\lambda_0$, we have $\Re\lambda>0,$ and then \eqref{cl:u=0} holds.

\item If $\Im\lambda\neq 0$ and $\Im\lambda\cdot \Im\delta>0$,
then it is easy to see \eqref{cl:u=0} from \eqref{eq:IM}. 

\item Suppose that $\Im\lambda\neq 0$ and $\Im\lambda\cdot \Im\delta<0.$ 
Inserting the following equality 
\begin{equation*}
	\gamma_{1+}\|\bu_+\|_+^2+\gamma_{1-}\|\bu_-\|_-^2=\Big|\frac{\Im\delta}{\Im\lambda}\Big|\|i\xi'\cdot \bu_+'+\pd_Nu_{+N}\|_+^2
\end{equation*}
into \eqref{eq:con} furnishes that
\begin{equation}\label{c4:2}
	0\geq \Big( \nu_++\Re\delta+\Re\lambda\Big|\frac{\Im\delta}{\Im\lambda}\Big| \Big)
	\|i\xi'\cdot \bu_+'+\pd_Nu_{+N}\|_+^2.
\end{equation}
By the fact that $$\Re\delta+\Re\lambda\Big|\frac{\Im\delta}{\Im\lambda}\Big|\geq 0,$$
 we obtain $\|i\xi'\cdot \bu_+'+\pd_Nu_{+N}\|_+^2=0$ and then the equality \eqref{cl:u=0}.
\end{itemize}

\item[(C3).]
Assume that $\delta\in\Sigma_{\varepsilon}$ with $\Re\delta\geq0$, $\lambda\in\mathbb{C}$ with $|\lambda|\geq \lambda_0$ and $\Re\lambda\geq\lambda_0|\Im\lambda|$. 
\begin{itemize}
\item If $\Im\lambda=0$ and $|\lambda|\geq\lambda_0,$  then $\Re\lambda>0.$
Thus \eqref{cl:u=0} holds true according to \eqref{eq:con} and $\Re\delta\geq0.$

\item For the case where $\Im\lambda\neq 0$ and $\Im\lambda\cdot \Im\delta>0$, 
we can easily obtain \eqref{cl:u=0} from \eqref{eq:IM}.

\item Suppose that $\Im\lambda\neq 0$ and $\Im\lambda\cdot \Im\delta\leq 0.$ According to \eqref{eq:IM} and \eqref{eq:con}, we also obtain the inequality \eqref{c4:2} as before. 
Since 
$$\Re\delta+\Re\lambda\Big|\frac{\Im\delta}{\Im\lambda}\Big|
\geq \lambda_0|\Im\delta|+\Re\delta\geq 0,$$
we have $\|i\xi'\cdot \bu_+'+\pd_Nu_{+N}\|_+^2=0$ and then \eqref{cl:u=0}.

\end{itemize}

\end{enumerate}

Summing up, we have $\det\bL\neq0$ for any $(\lambda,A)\in \CD_4.$
By the compactness of $\CD_4$, there exists a positive constant $C_5$ such that
$$\inf_{(\lambda,A)\in \CD_4}|\det\bL|\geq C_5 
\geq { \frac{C_5}{R_{3}^3 (1+R_2)^3} 
(|\lambda|^{1/2}+A)^3. }$$
This completes the proof of Lemma \ref{lem:3.1}.
\medskip


Let us end this section with some corollary from Lemma \ref{lem:3.1}, which concerns the orders of the multipliers appearing in \eqref{eq:pressure} and \eqref{eq:all_symbol}.
\begin{cor}\label{cor1}
Under the same assumptions in Lemma \ref{lem:4.3},  the following assertions hold true.

\begin{enumerate}[label=(\arabic*)]	
\item For the matrix $\bL$ defined by \eqref{def:LN},
$(\det\bL)^{-1}\in\bM_{- 3,2}(\widetilde{\Gamma}_{\varepsilon,\lambda_0})$.

\item Let $i=0,1$, $j=-1,0$ and $J,\ell=1,\dots,N$. 
For the symbols appearing in \eqref{eq:pressure} and \eqref{eq:all_symbol}, we have 
\begin{equation*}
\begin{aligned}
 & p_{\ell,i}^0\in\bM_{i,2}(\widetilde{\Gamma}_{\varepsilon,\lambda_0}), \quad 
p_{\ell,i}^b\in\bM_{i,2}(\widetilde{\Gamma}_{\varepsilon,\lambda_0}),\quad 
 R_{J\ell,i}^\pm\in\bM_{i,2}(\widetilde{\Gamma}_{\varepsilon,\lambda_0}), \quad
 R_{J\ell,i}^b\in\bM_{i,2}(\widetilde{\Gamma}_{\varepsilon,\lambda_0}),\\
&  S_{J\ell,j}^\pm\in\bM_{j,2}(\widetilde{\Gamma}_{\varepsilon,\lambda_0}),  \quad
T_{J\ell,i}^-\in\bM_{i,2}(\widetilde{\Gamma}_{\varepsilon,\lambda_0}),  \quad 
S_{J\ell,j}^b\in\bM_{j,2}(\widetilde{\Gamma}_{\varepsilon,\lambda_0}).
\end{aligned}
\end{equation*}
\end{enumerate}
\end{cor}
\begin{proof}
As the proof of Corollary \ref{cor1} is quite elementary, we omit the details.
\end{proof}

\section{Analysis of the solution operators of \eqref{eq:2.3}}
\label{sec:SO}
In this section, we shall establish the property of the $\CR$-boundedness of the solution operators of \eqref{eq:2.3}. 
To this end, let us recall the definition of $\CR$-boundedness of operator families.
\begin{dfn}\label{def:R}
	Let $X,Y$ be Banach spaces. A family of operators $\tau\subset\CL(X;Y)$ is $\CR$-bounded if for any $N\in\BN, T_j\in\tau, x_j\in X$ and the \emph{Rademacher} functions $r_j(u)={\rm sign }(\sin \,2^j\pi u)$ defined for $u\in [0,1]$, there holds
	$$\Big\|\sum_{j=1}^{N}r_jT_jx_j\Big\|_{L_p([0,1];Y)}\leq C_p\Big\|\sum_{j=1}^{N}r_jx_j\Big\|_{L_p([0,1];X)} \quad \text{for some}\,\,p\in [0,\infty).$$
	The smallest constant $C_p$ above is called the $\CR$-bound of $\tau$, denoted by $\CR_{\CL(X;Y)}(\tau)$.
\end{dfn}
\begin{rmk}\label{rmk:R-bdd}
	Let $X,Y,Z$ be Banach spaces.
	\begin{enumerate}[label=(\roman*)]
		\item If $\CT$ and $\CS$ are two $\CR$-bounded families in $\CL(X;Y)$, then the set
		$$\CT+\CS=\{T+S:T\in\CT,S\in\CS\}$$
		is also $\CR$-bounded and 
		$\CR_{\CL(X;Y)}(\CT+\CS)\leq \CR_{\CL(X;Y)}(\CT)+\CR_{\CL(X;Y)}(\CS).$
		\item If the families $\CT\subset\CL(X;Y)$ and $\CS\subset\CL(X;Z)$ are $\CR$-bounded, then the set $$\CS\CT=\{S\circ T:T\in\CT,S\in\CS\}\subset\CL(X;Z)$$ 
        is also $\CR$-bounded and $\CR_{\CL(X;Z)}(\CS\CT)\leq \CR_{\CL(X;Y)}(\CT)\CR_{\CL(Y;Z)}(\CS).$
	\end{enumerate}
\end{rmk}

For \eqref{eq:2.3}, we can prove the following result. 
\begin{thm}\label{thm:5.1}
	Assume that $0<\varepsilon<\pi/2$, $1<q<\infty$ and the constants $\mu_{\pm},$ $\nu_+,$ $\gamma_{1\pm},$ $\gamma_{2+}>0.$ Set that
	\begin{equation*}
		Z_q(\dot\Omega)=H^1_q(\dot\Omega)^N \times H^2_q(\dot\Omega)^N,\quad  
		\CZ_q(\dot\Omega)=L_q(\dot\Omega)^N\times H^1_q(\dot\Omega)^{N}
		\times L_q(\dot\Omega)^N \times H^1_q(\dot\Omega)^{N}\times H^2_q(\dot\Omega)^{N}.
	\end{equation*}
	For any $(\bh,\bk)\in Z_q(\dot\Omega),$ there exist constants $\lambda_0, r_b$ and the operator families
	\begin{equation*}
		\begin{aligned}
			\CS_{\pm}(\lambda)&\in \Hol\Big(\Gamma_{\varepsilon,\lambda_0};
			\CL\big( \CZ_q(\dot\Omega);\,H^2_q(\Omega_{\pm})^N\big)\Big),\\
			\CP_-(\lambda)&\in\Hol \Big(\Gamma_{\varepsilon,\lambda_0};
			\CL\big( \CZ_q(\dot\Omega);\,\wh H^1_q(\Omega_{-})\big)\Big)
		\end{aligned}	
	\end{equation*}
	such that 
	$$(\bv_{\pm},\fp_-)=\big(\CS_{\pm}(\lambda),\CP_-(\lambda)\big)
	(\lambda^{1/2}\bh,\bh,\lambda\bk,\lambda^{1/2}\bk,\bk)$$
	solves \eqref{eq:2.3}.  Here the set $\Gamma_{\ep,\lambda_0}$ is defined by \eqref{def:Gamma} for some $\delta_0>0.$
	Moreover, we have
 \begin{equation}\label{res:key_1}
     \begin{aligned}
		\CR_{\CL\big(\CZ_q (\dot\Omega)\,;\,H^{2-j}_q(\Omega_\pm)^N \big)}
		\Big(\Big\{ (\tau\pd_{\tau})^\ell \big( \lambda^{j\slash 2} \CS_{\pm}(\lambda) \big):
		\lambda\in\Gamma_{\varepsilon,\lambda_0} \Big\}\Big) & \leq r_b,\\
		\CR_{\CL\big( \CZ_q (\Omega_-)\,;\, L_q(\Omega_{-})^{N} \big)}
		\Big( \Big\{ (\tau\pd_{\tau})^\ell \nabla\CP_-(\lambda) :
		\lambda\in\Gamma_{\varepsilon,\lambda_0} \Big\}\Big) & \leq r_b
     \end{aligned}
 \end{equation}
	for $\ell =0,1,$ $j=0,1,2,$ and $\tau = \Im \lambda.$ 
	Above the choice of  $r_b$ depends solely on 
	$\varepsilon,$ $q,$ $N,$ $\gamma_{1\pm},$ $\gamma_{2+},$ $\mu_{\pm},$ 
	$\nu_+,$  $\lambda_0$, $\delta_0$ and $b.$
\end{thm}

In what follows, we first address some technical results on the boundedness of Fourier multipliers in Subsection \ref{subsec:TL}. 
Then we construct the operator families $\CS_{\pm}(\lambda)$ and $\CP_-(\lambda)$ and derive the bounds in  \eqref{res:key_1} in Subsections \ref{subsec:S+} and \ref{subsec:S-P} respectively.

\subsection{Technical lemmas}
\label{subsec:TL}
To prove Theorem \ref{thm:5.1},  we first recall some technical results used in the one-phase problems \cite{KSS2014, Saito2015}.
\begin{lem}[{\cite[Lemma 3.1]{KSS2014}}]
\label{lem:4.4}
	Assume that $0<\varepsilon<\pi/2$, $1<q<\infty$ and the constants $\gamma_{1+},$ $\gamma_{2+},$ $\mu_{+},$ 
	$\nu_+,$  $\lambda_0>0.$ Let the set $\Gamma_{\ep,\lambda_0}$ be defined by \eqref{def:Gamma} for some $\delta_0>0,$ and let the multiplier $m(\lambda,\xi')$ belong to $\bM_{-2,2}(\widetilde{\Gamma}_{\varepsilon,\lambda_0}).$ For $0<x_N<+\infty$, set that
	\begin{equation*}
		\begin{aligned}
			\big[K^1_{+,+}(\lambda)f\big](x',x_N)=&\int_{0}^\infty\mathcal{F}_{\xi'}^{-1}\Big[
			m(\lambda,\xi')AA_+M_+(x_N+y_N)\mathcal{F}_{y'}[f](\xi',y_N)\Big](x')dy_N,\\
			\big[K^2_{+,+}(\lambda)f\big](x',x_N)=&\int_{0}^\infty\mathcal{F}_{\xi'}^{-1}\Big[
			m(\lambda,\xi')Ae^{-B_+(x_N+y_N)}\mathcal{F}_{y'}[f](\xi',y_N)\Big](x')dy_N.
		\end{aligned}
	\end{equation*}
Then there exists a constant $r_b$ depending solely on 
$\varepsilon,$ $q,$ $N,$ $\gamma_{1+},$ $\gamma_{2+},$ $\mu_{+},$ 
	$\nu_+,$  $\lambda_0$, $\delta_0$ and $b$ such that
	$$\mathcal{R}_{\mathcal{L}\big(L_q(\Omega_+),H^{2-j}_q(\Omega_+)\big)}
	\Big(\Big\{(\tau\pd_{\tau})^\ell \lambda^{j/2} K^k_{+,+}(\lambda):
	\lambda\in \Gamma_{\varepsilon,\lambda_0} \Big\}\Big)\leq r_b$$
	for $j=0,1,2,$ $k=1,2,$ and $\ell=0,1.$ 
\end{lem}

For the results in the layer type domains,
let us introduce the cut-off functions $\varphi_{\pm}(x_N)\in C^\infty(\BR)$ satisfying 
\begin{equation*}
	0\leq \varphi_\pm(x_N)\leq 1, \quad 
	\varphi_+(x_N)=
	\begin{cases}
		1 & \text{for}\,\,\,x_N\leq 1/3,\\
		0 & \text{for}\,\,\, x_N\geq 2/3
	\end{cases}	
    \quad\text{and}\quad 
    \varphi_-(x_N)=
	\begin{cases}
		1 & \text{for}\,\,\,x_N\leq -2/3,\\
		0 & \text{for}\,\,\, x_N\geq -1/3.
	\end{cases}
\end{equation*}
Then we denote
\begin{equation}\label{def:cut-off-1}
\begin{aligned}
    \varphi_{-b}(x_N)&=\varphi_-(x_N/b), &  \varphi_{-0}(x_N)&=1-\varphi_-(x_N/b),\\
    \varphi_{+0}(x_N)&=\varphi_+(x_N/b), & \varphi_{b}(x_N)&=1-\varphi_+(x_N/b),
\end{aligned}
\end{equation}
and 
\begin{equation*}
 \begin{aligned}
 E_k^- (y_N)
=\begin{cases}
		\varphi_{-0}(y_N) & \text{for}\,\,\, k=1,\\
		\varphi'_{-0}(y_N) & \text{for}\,\,\, k=2,	
	\end{cases} \qquad 
    E_k^+ (y_N) =\begin{cases}
		\varphi_{+0}(y_N) & \text{for}\,\,\, k=1,\\
		\varphi'_{+0}(y_N) & \text{for}\,\,\, k=2.
	\end{cases}       
    \end{aligned}
\end{equation*}

\begin{lem}[{\cite[Lemmas 5.4 and 5.5]{Saito2015}}]\label{lem:4.5}
Assume that $0<\varepsilon<\pi/2$, $1<q<\infty$ and the constants $\lambda_0$, $\mu_{-},$ $\gamma_{1-}>0.$ Let the set $\Gamma_{\ep,\lambda_0}$ be defined by \eqref{def:Gamma} for some $\delta_0>0,$ and let the multiplier $m_1(\lambda,\xi')\in \bM_{-2,2}(\widetilde{\Gamma}_{\varepsilon,\lambda_0})$
    and $m_2(\lambda,\xi') \in \bM_{0,2}(\widetilde{\Gamma}_{\varepsilon,\lambda_0}).$
	For $X,Y\in \{A,B_-\}$, $d_1(x_N)=-x_N, d_2(x_N)=x_N+b$ with $-b<x_N<0$ and $k,m=1,2$, we define $K_{-,-}^{i,n}\,(i=1,2,3, n=1,2)$ by
	\begin{equation*}
		\begin{aligned}
			\big[K_{-,-}^{1,n}(\lambda)f\big](x',x_N)=&\int_{-b}^0\CF_{\xi'}^{-1}\Big[
			E_n^-(y_N)m_1(\lambda,\xi')Ae^{-Xd_k(x_N)}e^{-Yd_m(y_N)}\widehat{f}(\xi',y_N)\Big](x')dy_N,\\
			\big[K_{-,-}^{2,n}(\lambda)f\big](x',x_N)=&\int_{-b}^0\CF_{\xi'}^{-1}\Big[
			E_n^-(y_N)m_1(\lambda,\xi')AB_-M_-(-d_k(x_N))e^{-Xd_m(y_N)}\widehat{f}(\xi',y_N)\Big](x')dy_N,\\
            \big[K_{-,-}^{3,n}(\lambda)f\big](x',x_N)=&\int_{-b}^0\CF_{\xi'}^{-1}\Big[E_n^-(y_N)m_2(\lambda,\xi')e^{-A(d_k(x_N)+d_m(y_N))}\widehat{f}(\xi',y_N)\Big](x')dy_N.
		\end{aligned}
	\end{equation*}
Then there exists a constant $r_b$ depending solely on 
$N,$ $q,$ $\varepsilon,$ $\gamma_{1-},$ $\mu_{-},$ 
$\lambda_0$  and $b$ such that
\begin{align*}
\mathcal{R}_{\mathcal{L}\big(L_q(\Omega_{-}),H^{2-j}_q(\Omega_{-})\big)}
\Big(\Big\{(\tau\pd_{\tau})^\ell \lambda^{j/2}K_{-,-}^{1,n}(\lambda):
\lambda\in {\Gamma_{\varepsilon,\lambda_0} } \Big\}\Big)&\leq r_b,\\
\mathcal{R}_{\mathcal{L}\big(L_q(\Omega_{-}),H^{2-j}_q(\Omega_{-})\big)}
\Big(\Big\{(\tau\pd_{\tau})^\ell \lambda^{j/2}K_{-,-}^{2,n}(\lambda):
\lambda\in \Gamma_{\varepsilon,\lambda_0} \Big\}\Big)&\leq r_b,\\
\mathcal{R}_{\mathcal{L}\big(L_q(\Omega_{-}),L_q(\Omega_{-})\big)}
\Big(\Big\{(\tau\pd_{\tau})^\ell \nabla K_{-,-}^{3,n}(\lambda):
\lambda\in \Gamma_{\varepsilon,\lambda_0} \Big\}\Big)
&\leq r_b
\end{align*}
for $\ell=0,1$, $j=0,1,2$.
\end{lem}

Now, we state some technical results which arise from the gas-liquid two-phase problem in this paper.
For the proofs of Lemmas \ref{lem:-+}, \ref{lem:+-} and \ref{lem:+-2}, we refer to Appendix \ref{proof}.
\begin{lem}
\label{lem:-+}
	Assume that $0<\varepsilon<\pi/2$, $1<q<\infty$ the constants $\gamma_{1+},$ $\gamma_{2+},$ $\mu_{+},$ 
	$\nu_+,$  $\lambda_0>0.$  Let the set $\Gamma_{\ep,\lambda_0}$ be defined by \eqref{def:Gamma} for some $\delta_0>0,$ and let the multiplier $m(\lambda,\xi')$ belong to $\bM_{-2,2}(\widetilde{\Gamma}_{\varepsilon,\lambda_0}).$ For $0<x_N<+\infty$, we define $K_{-,+}^{i,n}\,(i,n=1,2)$ by
	\begin{equation*}
		\begin{aligned}
			\big[K_{-,+}^{1,n}(\lambda)f\big](x',x_N)=&\int_{-b}^0\CF_{\xi'}^{-1}\Big[
			E_n^-(y_N)m(\lambda,\xi')Ae^{-B_+x_N}e^{Ay_N}\widehat{f}(\xi',y_N)\Big](x')dy_N,\\
			\big[K_{-,+}^{2,n}(\lambda)f\big](x',x_N)=&\int_{-b}^0\CF_{\xi'}^{-1}\Big[
			E_n^-(y_N)m(\lambda,\xi')A^2M_+(x_N)e^{Ay_N}\widehat{f}(\xi',y_N)\Big](x')dy_N.
		\end{aligned}
	\end{equation*}
	Then there exists a constant $r_b$ depending solely on 
	$\varepsilon,$ $q,$ $N,$ $\gamma_{1+},$ $\gamma_{2+},$ $\mu_{+},$ 
	$\nu_+,$  $\lambda_0$, $\delta_0$ and $b$ such that
	$$\mathcal{R}_{\mathcal{L}\big(L_q(\Omega_{-}),H^{2-j}_q(\Omega_{+})\big)}
	\Big(\Big\{(\tau\pd_{\tau})^\ell \lambda^{j/2}K_{-,+}^{i,n}(\lambda):
	\lambda\in \Gamma_{\varepsilon,\lambda_0} \Big\}\Big)\leq r_b$$
	for $\ell=0,1$, $j=0,1,2.$ 
\end{lem}

Let us introduce the strip region for any $b>0:$
\begin{equation*}
    \Omega_{+b}=\{x=(x',x_N)\in \BR^N: 0<x_N<b\}.
\end{equation*}
\begin{lem}\label{lem:+-}
	Assume that $0<\varepsilon<\pi/2$, $1<q<\infty$ and the constants $\mu_{-},$ $\gamma_{1-}>0.$ Let the set $\Gamma_{\ep,\lambda_0}$ be defined by \eqref{def:Gamma} for some $\delta_0>0,$ and let the multiplier $m(\lambda,\xi')$ belong to $\bM_{-2,2}(\widetilde{\Gamma}_{\varepsilon,\lambda_0}).$ For $-b<x_N<0$, we define $K_{+,-}^{i,n}\,(i=1,2,3,4,\,n=1,2)$ by
	\begin{equation*}
		\begin{aligned}
			\big[K^{1,n}_{+,-}(\lambda)f\big](x',x_N)=&\int_{0}^b\mathcal{F}_{\xi'}^{-1}\Big[
			E_n^+(y_N)m(\lambda,\xi')A^2e^{-Ay_N}M_-(x_N)\mathcal{F}_{y'}[f](\xi',y_N)\Big](x')dy_N,\\
			\big[K^{2,n}_{+,-}(\lambda)f\big](x',x_N)=&\int_{0}^b\mathcal{F}_{\xi'}^{-1}\Big[
			E_n^+(y_N)m(\lambda,\xi')Ae^{-Ay_N}e^{B_-x_N}\mathcal{F}_{y'}[f](\xi',y_N)\Big](x')dy_N,\\
			\big[K^{3,n}_{+,-}(\lambda)f\big](x',x_N)=&\int_{0}^b\mathcal{F}_{\xi'}^{-1}\Big[
			E_n^+(y_N)m(\lambda,\xi')A^2e^{-Ay_N}M_-(-x_N-b)\mathcal{F}_{y'}[f](\xi',y_N)\Big](x')dy_N,\\
			\big[K^{4,n}_{+,-}(\lambda)f\big](x',x_N)=&\int_{0}^b\mathcal{F}_{\xi'}^{-1}\Big[
			E_n^+(y_N)m(\lambda,\xi')Ae^{-Ay_N}e^{B_-(-x_N-b)}\mathcal{F}_{y'}[f](\xi',y_N)\Big](x')dy_N.
		\end{aligned}
	\end{equation*}
	Then there exists a constant $r_b$ depending solely on 
	$\varepsilon,$ $q,$ $N,$ $\gamma_{1-},$  $\mu_{-},$ 
	  $\lambda_0$, $\delta_0$ and $b$ such that
	$$\mathcal{R}_{\mathcal{L}\big(L_q(\Omega_{+b}),H^{2-j}_q(\Omega_-)\big)}
	\Big(\Big\{(\tau\pd_{\tau})^\ell \lambda^{j/2} K^{i,n}_{+,-}(\lambda):
	\lambda\in \Gamma_{\varepsilon,\lambda_0} \Big\}\Big)\leq r_b$$
	for $\ell=0,1$, $j=0,1,2$.
\end{lem}

\begin{lem}\label{lem:+-2}
	Assume that $0<\varepsilon<\pi/2$, $1<q<\infty$ and the constants $\mu_{-},$ $\gamma_{1-}>0.$ Let the set $\Gamma_{\ep,\lambda_0}$ be defined by \eqref{def:Gamma} for some $\delta_0>0,$ and let the multiplier $m(\lambda,\xi')$ belong to $\bM_{0,2}(\widetilde{\Gamma}_{\varepsilon,\lambda_0}).$ For $-b<x_N<0$, we define $K_{+,-}^{i,n}\,(i=5,6,\,n=1,2)$ by
	\begin{equation*}
		\begin{aligned}
			\big[K^{5,n}_{+,-}(\lambda)f\big](x',x_N)=&\int_{0}^b\mathcal{F}_{\xi'}^{-1}\Big[
			E_n^+(y_N)m(\lambda,\xi')e^{-Ay_N}e^{Ax_N}\mathcal{F}_{y'}[f](\xi',y_N)\Big](x')dy_N,\\
			\big[K^{6,n}_{+,-}(\lambda)f\big](x',x_N)=&\int_{0}^b\mathcal{F}_{\xi'}^{-1}\Big[
			E_n^+(y_N)m(\lambda,\xi')e^{-Ay_N}e^{A(-x_N-b)}\mathcal{F}_{y'}[f](\xi',y_N)\Big](x')dy_N.
		\end{aligned}
	\end{equation*}
	Then there exists a constant $r_b$ depending solely on 
	$\varepsilon,$ $q,$ $N,$ $\gamma_{1-},$ $\mu_{-},$ 
  $\lambda_0$, $\delta_0$ and $b$ such that
	$$\mathcal{R}_{\mathcal{L}\big(L_q(\Omega_{+b}),L_q(\Omega_-)\big)}
	\Big(\Big\{(\tau\pd_{\tau})^\ell  \nabla K^{i,n}_{+,-}(\lambda):
	\lambda\in \Gamma_{\varepsilon,\lambda_0} \Big\}\Big)\leq r_b$$
	for $\ell=0,1$.
\end{lem}

\subsection{Bound of the operator family $\CS_{+}(\lambda)$}
\label{subsec:S+}
In this and next subsection, we shall construct the solution operators in Theorem \ref{thm:5.1} from \eqref{eq:2.27}. Let us take $(\bh,\bk)\in Z_q(\dot\Omega).$ Then there exist $\bh_\pm\in H_q^1(\Omega_{\pm})^N$ and $\bk_\pm\in H_q^2(\Omega_{\pm})^N$ such that
\begin{equation*}
	\bh=\bh_+\mathds{1}_{\Omega_{+}} + \bh_-\mathds{1}_{\Omega_{-}} 
	\quad \text{and}\quad 
	\bk=\bk_+\mathds{1}_{\Omega_{+}} + \bk_-\mathds{1}_{\Omega_{-}},  
\end{equation*}
with $\mathds{1}_{\Omega_{\pm}}$ standing for the indicator functions of $\Omega_{\pm}.$ In particular, taking traces on the interface $\Gamma,$ we have 
\begin{equation*}
	\bh(x',0)=\bh_+(x',0)+\bh_-(x',0)  
	\quad \text{and}\quad 
	\bk(x',0)=\bk_+(x',0)+\bk_-(x',0).
\end{equation*}

\begin{lem}[Volevich's trick]
\label{lem:Vol-1}
For any function $a_+(\xi',x_N)$ vanishing as $x_N\to \infty,$ we have 
\begin{equation}\label{Vol:+1}
\begin{aligned}
a_+(\xi',x_N)\widehat{\bh}_+(\xi',0)
           =  &-\int_{0}^{\infty}(\pd_N a_+)(\xi',x_N+y_N)
    \frac{\gamma_{1+}\lambda^{1/2}}{\mu_{+}B_{+}^2}
    \widehat{\lambda^{1/2}\bh_+}(\xi',y_N) dy_N\\
           &+\sum_{\ell=1}^{N-1}\int_{0}^{\infty}(\pd_N a_+)(\xi',x_N+y_N)
\frac{i\xi_\ell}{B_{+}^2} \widehat{\pd_\ell \bh_+}(\xi',y_N) dy_N\\
    &-\int_{0}^{\infty} a_+(\xi',x_N+y_N)\widehat{\pd_N\bh_+}(\xi',y_N)dy_N,
	\end{aligned}
\end{equation}
and 
\begin{equation}\label{Vol:+2}
	\begin{aligned}
		a_+(\xi',x_N)\widehat{\bk}_+(\xi',0)=
  		&-\int_{0}^{\infty}(\pd_N a_+)(\xi',x_N+y_N)
  \frac{\gamma_{1+}}{\mu_{+}B_{+}^2} \widehat{\lambda \bk_+}(\xi',y_N) dy_N\\
  &-\int_{0}^{\infty}a_+(\xi',x_N+y_N)
  \frac{\gamma_{1+}\lambda^{1/2}}{\mu_{+}B_{+}^2}
  \widehat{\lambda^{1/2}\pd_N\bk_+}(\xi',y_N) dy_N\\
    		&+\int_{0}^{\infty}(\pd_N a_+ )(\xi',x_N+y_N)
    \frac{1}{B_{+}^2}\widehat{\Delta'\bk_+}(\xi',y_N)dy_N\\
    &+\sum_{\ell=1}^{N-1}\int_{0}^{\infty}a_+(\xi',x_N+y_N)
    \frac{i\xi_{\ell}}{B_{+}^2} \widehat{\pd_{\ell}\pd_N\bk_+}(\xi',y_N)dy_N
	\end{aligned}
\end{equation}
for $\Delta'=\sum_{\ell=1}^{N-1} \pd_{\ell}^2.$
\smallbreak 

In addition, let the cut-off function $\varphi_{-0}$ be defined in \eqref{def:cut-off-1}. Then we have 
\begin{equation}\label{Vol:-1}
	\begin{aligned}
		\widehat{\bh}_-(\xi',0)
  =&\int_{-b}^0 \big( \varphi'_{-0}(y_N) +\varphi_{-0}(y_N)A\big)
  e^{Ay_N} \frac{\gamma_{1+}\lambda^{1/2}}{\mu_+ B_+^2} \widehat{\lambda^{1/2}\bh_-}(\xi',y_N)dy_N\\
    &  - \sum_{\ell=1}^{N-1} \int_{-b}^0 \big( \varphi'_{-0}(y_N) +\varphi_{-0}(y_N)A\big)
  e^{Ay_N} \frac{i\xi_{\ell}}{B_+^2} \widehat{\pd_{\ell}\bh_-}(\xi',y_N)dy_N\\
		&+\int_{-b}^0\varphi_{-0}(y_N)e^{Ay_N}\widehat{\pd_N \bh_-}(\xi',y_N)dy_N.
	\end{aligned}
\end{equation}
\begin{equation}\label{Vol:-2}
	\begin{aligned}
		\widehat{\bk}_-(\xi',0)
  =&\int_{-b}^0 \big( \varphi'_{-0}(y_N) +\varphi_{-0}(y_N)A\big)
  e^{Ay_N} \frac{\gamma_{1+}}{\mu_+ B_+^2} \widehat{\lambda\bk_-}(\xi',y_N)dy_N\\
 &+\int_{-b}^0\varphi_{-0}(y_N)e^{Ay_N}\frac{\gamma_{1+}\lambda^{1/2}}{\mu_+ B_+^2}
 \widehat{\lambda^{1/2}\pd_N \bk_-}(\xi',y_N)dy_N\\
  &  - \int_{-b}^0 \big( \varphi'_{-0}(y_N) +\varphi_{-0}(y_N)A\big) e^{Ay_N} 
    \frac{1}{B_+^2} \widehat{\Delta'\bk_-}(\xi',y_N)dy_N\\
		&-\sum_{\ell=1}^{N-1}\int_{-b}^0\varphi_{-0}(y_N)e^{Ay_N}
  \frac{i\xi_{\ell}}{B_+^2} \widehat{\pd_{\ell}\pd_N \bk_-}(\xi',y_N)dy_N.
	\end{aligned}
\end{equation}
\end{lem}

\begin{proof}
For any function $f_+$ defined on $\Omega_+,$ it is not hard to see that 
\begin{equation}\label{eq:vh}
\begin{aligned}
a_+(\xi',x_N)\widehat{f_+}(\xi',0) = & -\int_{0}^{\infty}
	(\pd_N a_+)(\xi^\prime,x_N+y_N)\widehat{f_+}(\xi',y_N) d y_N\\
	& -\int_{0}^{\infty} a_+(\xi^\prime,x_N+y_N)\widehat{\pd_N f_+}(\xi',y_N) d y_N.
\end{aligned}
\end{equation}
Then taking $f_+$ in \eqref{eq:vh} as the components of $(\bh_+,\bk_+)$ and inserting the identity 
\begin{equation}\label{eq:1_+}
    1=\frac{\gamma_{1+}\lambda}{\mu_{+}B_{+}^2}
-\sum_{\ell=1}^{N-1}\frac{(i\xi_{\ell})\cdot (i\xi_{\ell})}{B_{+}^2}
\end{equation}
into \eqref{eq:vh} yield \eqref{Vol:+1} and \eqref{Vol:+2}. 

On the other hand, for any function $f_-$ defined on $\Omega_-,$ we have 
\begin{equation}\label{eq:-}
	\begin{aligned}
		\widehat{f_-}(\xi',0)=&\int_{-b}^{0}\frac{d}{dy_N}\big[\varphi_{-0}(y_N)e^{Ay_N}
  \widehat{f_-}(\xi',y_N)\big]dy_N\\
		=&\int_{-b}^0 \big( \varphi'_{-0}(y_N) +\varphi_{-0}(y_N)A\big)
  e^{Ay_N}\widehat{f_-}(\xi',y_N)dy_N\\
		&+\int_{-b}^0\varphi_{-0}(y_N)e^{Ay_N}\widehat{\pd_N f_-}(\xi',y_N)dy_N.
	\end{aligned}
\end{equation}
Then using \eqref{eq:1_+} and \eqref{eq:-} implies \eqref{Vol:-1} and \eqref{Vol:-2}.
\end{proof}


In what follows, for $(\bh,\bk)\in Z_q(\dot\Omega),$ we construct the solution operator of $\bv_+,$ namely,
\begin{equation*}
		\bv_+=\CS_+(\lambda)(\lambda^{1/2}\bh,\bh,\lambda\bk,\lambda^{1/2}\bk,\bk).
\end{equation*}
with $\CS_+(\lambda)=\big(\CS_{1+}(\lambda),\cdots,\CS_{N+}(\lambda)\big).$
For simplicity, let us set the notations
\begin{equation}\label{conv:f}
    \begin{aligned}
  f_1&=(f_{1\ell+}\mathds{1}_{\Omega_{+}}+f_{1\ell-}\mathds{1}_{\Omega_{-}}\,|\,\ell=1,\cdots,N),\\
  f_2&=(f_{2m\ell+}\mathds{1}_{\Omega_{+}}+f_{2m\ell- }\mathds{1}_{\Omega_{-}}\,|\,\ell,m=1,\cdots,N), \\
  f_3&=(f_{3\ell+}\mathds{1}_{\Omega_{+}}+f_{3\ell-}\mathds{1}_{\Omega_{-}}\,|\,\ell=1,\cdots,N),\\
  f_4&=(f_{4m\ell+}\mathds{1}_{\Omega_{+}}+f_{4m\ell-}\mathds{1}_{\Omega_{-}}\,|\,\ell,m=1,\cdots,N),\\
  f_5&=(f_{5 mn\ell+}\mathds{1}_{\Omega_{+}}+f_{5 mn\ell-}\mathds{1}_{\Omega_{-}}|\,\ell,m,n=1,\cdots ,N),
    \end{aligned}
\end{equation}
which correspond to the following variables
\begin{equation*}
 \begin{aligned}
   \lambda^{1/2} \bh&=\big(\lambda^{1/2}(h_{\ell+}\mathds{1}_{\Omega_{+}}+h_{\ell-}\mathds{1}_{\Omega_{-}})\,|\,\ell=1,\cdots,N\big),\\  
    \nabla \bh&=(\pd_mh_{\ell+}\mathds{1}_{\Omega_{+}}+\pd_mh_{\ell-}\mathds{1}_{\Omega_{-}}\,|\,m,\ell=1,\cdots,N), \\
    \lambda \bk &=\big(\lambda (k_{\ell+}\mathds{1}_{\Omega_{+}}+k_{\ell-}\mathds{1}_{\Omega_{-}})\,|\,\ell=1,\cdots,N\big),\\
    \lambda^{1/2}\nabla \bk & =\big(\lambda^{1/2} \pd_m(k_{\ell+}\mathds{1}_{\Omega_{+}}+k_{\ell-}\mathds{1}_{\Omega_{-}})\,|\,m,\ell=1,\cdots,N\big),\\ \nabla^2\bk&=\big(\pd_n\pd_m(k_{\ell+}\mathds{1}_{\Omega_{+}}+k_{\ell-}\mathds{1}_{\Omega_{-}})\,|\,n,m,\ell=1,\cdots,N\big).
 \end{aligned}   
\end{equation*}

Now, using \eqref{eq:2.27} and the convention \eqref{conv:f}, we will see that
\begin{equation}\label{eq:S+}
    \CS_{J+}(\lambda^{1/2}\bh, \bh,\lambda \bk,\lambda^{1/2}\bk,\bk)
   =\sum_{i=1}^2\CS_{ Ji+} (\lambda^{1/2}\bh, \bh,\lambda \bk,\lambda^{1/2}\bk,\bk)    =\sum_{i=1}^2\sum_{k=1}^5\CS_{Ji+}^k(\lambda)(f_k)
\end{equation}
for any $J=1,\cdots,N.$ Moreover, note that
\begin{equation}\label{eq:MeB+}
	\begin{aligned}
		\pd_NM_{+}(x_N+y_N)&=-\big(e^{- B_{+}(x_N+y_N)}+A_{+}M_{+}(x_N+y_N)\big),\\
		\pd_Ne^{- B_{+}(x_N+y_N)}&=- B_{+}e^{- B_{+}(x_N+y_N)}.
	\end{aligned}
\end{equation}
Then Lemma \ref{lem:Vol-1} and \eqref{eq:MeB+} give the form of the operators $\CS_{Ji+}^k,$ $i=1,2,$ $k=1,\dots,5,$ as follows:
\begin{equation*}
\begin{aligned}
    \CS_{J1+}^1(\lambda)(f_1)=&\sum_{\ell=1}^N\int_{0}^{\infty}\CF_{\xi'}^{-1}
    \Big[\big(e^{-B_+(x_N+y_N)}+A_+M_{+}(x_N+y_N)\big)
    \frac{R_{J\ell,0}^+\gamma_{1+}\lambda^{1/2}}{\mu_{+}B_{+}^2}
    \widehat{f_{1\ell+}}(\xi',y_N)\Big]dy_N\\
    &+\sum_{\ell=1}^N\int_{-b}^0\CF_{\xi'}^{-1}\Big[\big(\varphi'_{-0}(y_N)+\varphi_{-0}(y_N)A\big)e^{Ay_N}M_+(x_N)\frac{R_{J\ell,0}^+\gamma_{1+}\lambda^{1/2}}{\mu_{+}B_+^2}\widehat{f_{1\ell-}}(\xi',y_N)\Big]dy_N,
\end{aligned}
\end{equation*}

\begin{equation*}
    \begin{aligned}
    \CS_{J1+}^2(\lambda)(f_2)
    =&-\sum_{\ell=1}^N\sum_{m=1}^{N-1}\int_{0}^{\infty}\CF_{\xi'}^{-1}\Big[\big(e^{-B_+(x_N+y_N)}+A_+M_{+}(x_N+y_N)\big)
    \frac{R_{J\ell,0}^+i\xi_m}{B_{+}^2}\widehat{ f_{2m\ell+}}(\xi',y_N)\Big]dy_N\\
    &-\sum_{\ell=1}^N\int_{0}^{\infty}\CF_{\xi'}^{-1}\Big[
    M_{+}(x_N+y_N)R_{J\ell,0}^+\widehat{ f_{2N\ell+}}(\xi',y_N)\Big]dy_N\\
    &-\sum_{\ell=1}^N\sum_{m=1}^{N-1}\int_{-b}^0\CF_{\xi'}^{-1}\Big[\big(\varphi'_{-0}(y_N)+\varphi_{-0}(y_N)A\big)e^{Ay_N}M_+(x_N)\frac{R_{J\ell,0}^+i\xi_m}{B_+^2}\widehat{ f_{2m\ell-}}(\xi',y_N)\Big]dy_N\\
    &+\sum_{\ell=1}^N\int_{-b}^0\CF_{\xi'}^{-1}\Big[\varphi_{-0}(y_N)e^{Ay_N}M_+(x_N)
     R_{J\ell,0}^+\widehat{ f_{2N\ell-}}(\xi',y_N)\Big]dy_N,
\end{aligned}
\end{equation*}

\begin{equation*}
\begin{aligned}
    \CS_{J1+}^3(\lambda)(f_3)=&\sum_{\ell=1}^N\int_{0}^{\infty}\CF_{\xi'}^{-1}\Big[\big(e^{-B_+(x_N+y_N)}+A_+M_{+}(x_N+y_N)\big)\frac{R_{J\ell,1}^+\gamma_{1+}}{\mu_{+}B_{+}^2}\widehat{f_{3\ell+}}(\xi',y_N)\Big]dy_N\\
    &+\sum_{\ell=1}^N\int_{-b}^0\CF_{\xi'}^{-1}\Big[\big(\varphi'_{-0}(y_N)+\varphi_{-0}(y_N)A\big)e^{Ay_N}M_+(x_N)\frac{R_{J\ell,1}^+\gamma_{1+}}{\mu_{+}B_+^2}\widehat{f_{3\ell-}}(\xi',y_N)\Big]dy_N,
\end{aligned}
\end{equation*}

\begin{equation*}
\begin{aligned}
    \CS_{J1+}^4(\lambda)(f_4)=&-\sum_{\ell=1}^N\int_{0}^{\infty}\CF_{\xi'}^{-1}
    \Big[M_+(x_N+y_N)\frac{R_{J\ell,1}^+\gamma_{1+}\lambda^{1/2}}{\mu_{+}B_{+}^2}
    \widehat{ f_{4N\ell+}}(\xi',y_N)\Big]dy_N\\
    &+\sum_{\ell=1}^N\int_{-b}^0\CF_{\xi'}^{-1}\Big[\varphi_{-0}(y_N)e^{Ay_N}M_+(x_N)\frac{R_{J\ell,1}^+\gamma_{1+}\lambda^{1/2}}{\mu_{+}B_+^2}\widehat{ f_{4N\ell-}}(\xi',y_N)\Big]dy_N,
\end{aligned}
\end{equation*}

\begin{equation*}
    \begin{aligned}
    \CS_{J1+}^5(\lambda)(f_5)=&-\sum_{\ell=1}^N\sum_{m=1}^{N-1}\int_{0}^{\infty}\CF_{\xi'}^{-1}\Big[\big(e^{-B_+(x_N+y_N)}+A_+M_{+}(x_N+y_N)\big)\frac{R_{J\ell,1}^+}{B_{+}^2}\widehat{ f_{5mm\ell+}}(\xi',y_N)\Big]dy_N\\
    &+\sum_{\ell=1}^N\sum_{m=1}^{N-1}\int_{0}^{\infty}\CF_{\xi'}^{-1}\Big[ 
    M_{+}(x_N+y_N)\frac{R_{J\ell,1}^+i\xi_m}{B_+^2}\widehat{ f_{5mN\ell+}}(\xi',y_N)\Big]dy_N\\
    &-\sum_{\ell=1}^N\sum_{m=1}^{N-1}\int_{-b}^0\CF_{\xi'}^{-1}\Big[
    \big(\varphi'_{-0}(y_N)+\varphi_{-0}(y_N)A\big)e^{Ay_N}M_+(x_N)
    \frac{R_{J\ell,1}^+}{B_+^2}\widehat{ f_{5mm\ell-}}(\xi',y_N)\\
    &-\sum_{\ell=1}^N\sum_{m=1}^{N-1}\int_{-b}^0\CF_{\xi'}^{-1}\Big[
    \varphi_{-0}(y_N)e^{Ay_N}M_+(x_N) \frac{R_{J\ell,1}^+i\xi_m}{B_+^2}
    \widehat{ f_{5mN\ell-}}(\xi',y_N)\Big]dy_N,
\end{aligned}
\end{equation*}

\begin{equation*}
\begin{aligned}
  \CS_{J2+}^1(\lambda)(f_1)=&\sum_{\ell=1}^N\int_{0}^{\infty}\CF_{\xi'}^{-1}\Big[
  e^{-B_+(x_N+y_N)}\frac{S_{J\ell,-1}^+\gamma_{1+}\lambda^{1/2}}{\mu_{+}B_{+}}
  \widehat{f_{1\ell+}}(\xi',y_N)\Big]dy_N\\
    &+\sum_{\ell=1}^N\int_{-b}^0\CF_{\xi'}^{-1}\Big[\big(\varphi'_{-0}(y_N)+\varphi_{-0}(y_N)A\big)e^{Ay_N}e^{-B_+x_N}\frac{S_{J\ell,-1}^+\gamma_{1+}\lambda^{1/2}}{\mu_{+}B_+^2}\widehat{f_{1\ell-}}(\xi',y_N)\Big]dy_N,
\end{aligned}
\end{equation*}

\begin{equation*}
    \begin{aligned}
    \CS_{J2+}^2(\lambda)(f_2)=&-\sum_{\ell=1}^N\sum_{m=1}^{N-1}\int_{0}^{\infty}\CF_{\xi'}^{-1}\Big[
    e^{-B_+(x_N+y_N)}\frac{S_{J\ell,-1}^+i\xi_m}{B_{+}}
    \widehat{ f_{2m\ell+}}(\xi',y_N)\Big]dy_N\\
    &-\sum_{\ell=1}^N\int_{0}^{\infty}\CF_{\xi'}^{-1}\Big[
    e^{-B_+(x_N+y_N)}S_{J\ell,-1}^+\widehat{ f_{2N\ell+}}(\xi',y_N)\Big]dy_N\\
    &-\sum_{\ell=1}^N \sum_{m=1}^{N-1}\int_{-b}^0\CF_{\xi'}^{-1}\Big[
    \big(\varphi'_{-0}(y_N)+\varphi_{-0}(y_N)A\big)e^{Ay_N}e^{-B_+x_N}
   \frac{S_{J\ell,-1}^+i\xi_m}{B_+^2}
    \widehat{ f_{2m\ell-}}(\xi',y_N)\Big]dy_N\\
    &+\sum_{\ell=1}^N\int_{-b}^0\CF_{\xi'}^{-1}\Big[
    \varphi_{-0}(y_N)e^{Ay_N}e^{-B_+x_N}S_{J\ell,-1}^+
    \widehat{ f_{2N\ell-}}(\xi',y_N)\Big]dy_N,
\end{aligned}
\end{equation*}

\begin{equation*}
\begin{aligned}
    \CS_{J2+}^3(\lambda)(f_3)=&\sum_{\ell=1}^N\int_{0}^{\infty}\CF_{\xi'}^{-1}\Big[
    e^{-B_+(x_N+y_N)}\frac{S_{J\ell,0}^+\gamma_{1+}}{\mu_{+}B_{+}}
    \widehat{f_{3\ell+}}(\xi',y_N)\Big]dy_N\\
    &+\sum_{\ell=1}^N\int_{-b}^0\CF_{\xi'}^{-1}\Big[\big(\varphi'_{-0}(y_N)+\varphi_{-0}(y_N)A\big)e^{Ay_N}e^{-B_+x_N}\frac{S_{J\ell,0}^+\gamma_{1+}}{\mu_{+}B_+^2}\widehat{f_{3\ell-}}(\xi',y_N)\Big]dy_N,
\end{aligned}
\end{equation*}
\begin{equation*}
\begin{aligned}
    \CS_{J2+}^4(\lambda)(f_4)=&{-}\sum_{\ell=1}^N\int_{0}^{\infty}\CF_{\xi'}^{-1}\Big[
    e^{-B_+(x_N+y_N)}\frac{S_{J\ell,0}^+\gamma_{1+}\lambda^{1/2}}{\mu_{+}B_{+}^2}
    \widehat{ f_{4N\ell+}}(\xi',y_N)\Big]dy_N\\
    &+\sum_{\ell=1}^N\int_{-b}^0\CF_{\xi'}^{-1}\Big[
    \varphi_{-0}(y_N)e^{Ay_N} { e^{-B_+x_N}
    \frac{S_{J\ell,0}^+\gamma_{1+}\lambda^{1/2}}{\mu_{+}B_+^2} }
    \widehat{ f_{4N\ell-}}(\xi',y_N)\Big]dy_N,
\end{aligned}
\end{equation*}
\begin{equation*}
    \begin{aligned}
    \CS_{J2+}^5(\lambda)(f_5)=&-\sum_{\ell=1}^N\sum_{m=1}^{N-1}\int_{0}^{\infty}\CF_{\xi'}^{-1}\Big[
    e^{-B_+(x_N+y_N)}\frac{S_{J\ell,0}^+}{B_{+}}
    \widehat{ f_{5mm\ell+}}(\xi',y_N)\Big]dy_N\\
    &{+}\sum_{\ell=1}^N\sum_{m=1}^{N-1} \int_{0}^{\infty}\CF_{\xi'}^{-1}\Big[
    e^{-B_+(x_N+y_N)}\frac{S_{J\ell,0}^+i\xi_m}{B_+^2}
    \widehat{ f_{5mN\ell+}}(\xi',y_N)\Big]dy_N\\
    &-\sum_{\ell=1}^N\sum_{m=1}^{N-1}\int_{-b}^0\CF_{\xi'}^{-1}\Big[
    \big(\varphi'_{-0}(y_N)+\varphi_{-0}(y_N)A\big)e^{Ay_N}e^{-B_+x_N}
    \frac{S_{J\ell,0}^+}{B_+^2}\widehat{ f_{5mm\ell-}}(\xi',y_N)\Big]dy_N\\
    &-\sum_{\ell=1}^N\sum_{m=1}^{N-1}\int_{-b}^0\CF_{\xi'}^{-1}\Big[
    \varphi_{-0}(y_N)e^{Ay_N}e^{-B_+x_N}\frac{S_{J\ell,0}^+i\xi_m}{B_+^2}
    \widehat{ f_{5mN\ell-}}(\xi',y_N)\Big]dy_N.
\end{aligned}
\end{equation*}

At last, let us derive the bounds of the operator $\CS_+(\lambda).$
Corollary \ref{cor1} implies that all the multipliers
\begin{gather*}
	\frac{R_{J\ell,0}^+\gamma_{1+}\lambda^{1/2}}{\mu_{+}B_+^2A},\,\, \frac{R_{J\ell,0}^+i\xi_m}{B_+^2A}, \frac{R_{J\ell,0}^+}{AA_+},\,\,
	\frac{R_{J\ell,1}^+}{B_+^2A},\,\, \frac{R_{J\ell,1}^+\gamma_{1+}\lambda^{1/2}}{\mu_{+}B_+^2AA_+},\\ \frac{R_{J\ell,1}^+i\xi_m}{\mu_{+}B_+^2AA_+},\,\,
	\frac{R_{J\ell,0}^+}{A^2},\,\, \frac{R_{J\ell,1}^+\gamma_{1+}\lambda^{1/2}}{\mu_{+}B_+^2A^2},\,\, \frac{R_{J\ell,1}^+i\xi_m}{\mu_{+}B_+^2A^2},\\
	\frac{S_{J\ell,-1}^+\gamma_{1+}\lambda^{1/2}}{\mu_{+}B_+A},\,\, \frac{S_{J\ell,-1}^+i\xi_m}{B_+A},\,\, \frac{S_{J\ell,0}^+}{AB_+},\,\,
	\frac{S_{J\ell,-1}^+}{A},\,\, \frac{S_{J\ell,0}^+\gamma_{1+}\lambda^{1/2}}{\mu_{+}B_+^2A},\,\, \frac{S_{J\ell,0}^+i\xi_m}{B_+^2A},
\end{gather*}
which appear in formulations of $\CS_{Ji+}^k(\lambda)$ above, belong to the class
$\bM_{-2,2}(\widetilde{\Gamma}_{\ep, \lambda_0}).$ 
Then it is not hard to see from Lemmas \ref{lem:4.4} and \ref{lem:-+} that
\begin{align*}
	&\CR_{\CL\big(\CZ_q(\dot\Omega),H^{2-j}_q(\Omega_{+})^N\big)}\Big(\Big\{(\tau\pd_{\tau})^\ell\lambda^{j/2} \CS_{Ji+}(\lambda)\,\big|\,\lambda\in\Gamma_{\varepsilon,\lambda_0}\Big\}\Big)\leq r_b
\end{align*}
with $\ell=0,1, i=1,2, j=0,1,2$, which yields the $\CR$-bound of $\CS_{+}(\lambda)$ in \eqref{res:key_1}.

\subsection{Bounds of the operator families $\CS_-(\lambda)$ and $\CP_-(\lambda)$}
\label{subsec:S-P}
This subsection is dedicated to deriving the bounds of $\CS_-(\lambda)$ and $\CP_-(\lambda),$ namely, 
\begin{equation*}
	\begin{aligned}
		\bv_-&=\CS_-(\lambda)(\lambda^{1/2}\bh,\bh,\lambda\bk,\lambda^{1/2}\bk,\bk),\\
		\fp_-&=\CP_-(\lambda)\big(\lambda^{1/2}\bh,\bh,\lambda\bk,\lambda^{1/2}\bk,\bk\big),
	\end{aligned}
\end{equation*}
with $\CS_-(\lambda)=\big(\CS_{1-}(\lambda),\cdots,\CS_{N-}(\lambda)\big).$
However, according to the difference in the topology properties between the domains $\Omega_+$ and $\Omega_-,$
we use the following technical lemma to treat the traces of $(\bh_+,\bk_+).$
\begin{lem}[Volevich's trick]
\label{lem:Vol-2}
Let $\varphi_{+0}$ be defined in \eqref{def:cut-off-1} and $\Delta'=\sum_{\ell=1}^{N-1} \pd_{\ell}^2.$ Then we have 
\begin{equation}\label{Vol:1}
	\begin{aligned}
		\widehat{\bh_+}(\xi',0)
  =&-\int_{0}^b \big( \varphi'_{+0}(y_N) - \varphi_{+0}(y_N)A\big)
  e^{-Ay_N} \frac{\gamma_{1-}\lambda^{1/2}}{\mu_- B_-^2} \widehat{\lambda^{1/2}\bh_+}(\xi',y_N)dy_N\\
    &  + \sum_{\ell=1}^{N-1} \int_{0}^b \big(\varphi'_{+0}(y_N) -\varphi_{+0}(y_N)A\big)
  e^{-Ay_N} \frac{i\xi_{\ell}}{B_-^2} \widehat{\pd_{\ell}\bh_+}(\xi',y_N)dy_N\\
		&-\int_{0}^b\varphi_{+0}(y_N)e^{-Ay_N}\widehat{\pd_N \bh_+}(\xi',y_N)dy_N,
	\end{aligned}
\end{equation}
\begin{equation}\label{Vol:2}
	\begin{aligned}
		\widehat{\bk_+}(\xi',0)
  =&-\int_{0}^b \big( \varphi'_{+0}(y_N) -\varphi_{+0}(y_N)A\big)
  e^{-Ay_N} \frac{\gamma_{1-}}{\mu_- B_-^2} \widehat{\lambda\bk_+}(\xi',y_N)dy_N\\
 &-\int_{0}^b\varphi_{+0}(y_N)e^{-Ay_N}\frac{\gamma_{1-}\lambda^{1/2}}{\mu_- B_-^2}
 \widehat{\lambda^{1/2}\pd_N \bk_+}(\xi',y_N)dy_N\\
  &  +\int_{0}^b \big( \varphi'_{+0}(y_N)-\varphi_{+0}(y_N)A\big) e^{-Ay_N} 
    \frac{1}{B_-^2} \widehat{\Delta'\bk_+}(\xi',y_N)dy_N\\
		&+\sum_{\ell=1}^{N-1}\int_{0}^b\varphi_{+0}(y_N)e^{-Ay_N}
  \frac{i\xi_{\ell}}{B_-^2} \widehat{\pd_{\ell}\pd_N \bk_+}(\xi',y_N)dy_N.
	\end{aligned}
\end{equation}
\end{lem}
\begin{proof}
Analogous to \eqref{Vol:-1} and \eqref{Vol:-2}, it is not hard to prove the equalities \eqref{Vol:1} and \eqref{Vol:2}. Here, we omit the details of the proof.
\end{proof}

Now, recall \eqref{eq:2.27} and the convention \eqref{conv:f}. Then, analogous to \eqref{eq:S+}, we will also see that
\begin{equation}\label{eq:S-}
\begin{aligned}
    	\CS_{J-}(\lambda)
        (\lambda^{1/2}\bh, \bh,\lambda \bk,\lambda^{1/2}\bk,\bk)
        &=\sum_{i=1}^{4}\CS_{ Ji-}(\lambda)
       (\lambda^{1/2}\bh, \bh,\lambda \bk,\lambda^{1/2}\bk,\bk)\\
       &=\sum_{i=1}^{4}\sum_{k=1}^5\CS_{Ji-}^k(\lambda)(f_k),\\
\CP_{-}(\lambda)\big(\lambda^{1/2}\bh,\bh,\lambda\bk,\lambda^{1/2}\bk,\bk\big)
&=\sum_{k=1}^5\CP_{-}^k(\lambda)(f_k)
\end{aligned}
\end{equation}
for any $J=1,\cdots,N.$
In fact, Lemmas \ref{lem:Vol-1} and \ref{lem:Vol-2} furnish that 
\begin{equation*}
\begin{aligned}
    \CS_{J1-}^1(\lambda)(f_1)=&\sum_{\ell=1}^N\int_{0}^{b}\CF_{\xi'}^{-1}\Big[\big(-\varphi'_{+0}(y_N)+\varphi_{+0}(y_N)A\big)e^{-Ay_N}M_-(x_N)\frac{R_{J\ell,0}^-\gamma_{1-}\lambda^{1/2}}{\mu_{-}B_{-}^2}\widehat{f_{1\ell+}}(\xi',y_N)\Big]dy_N\\
    &+\sum_{\ell=1}^N\int_{-b}^0\CF_{\xi'}^{-1}\Big[\big(\varphi'_{-0}(y_N)+\varphi_{-0}(y_N)A\big)e^{Ay_N}M_-(x_N)\frac{R_{J\ell,0}^-\gamma_{1-}\lambda^{1/2}}{\mu_{-}B_-^2}\widehat{f_{1\ell-}}(\xi',y_N)\Big]dy_N,
\end{aligned}
\end{equation*}

\begin{equation*}
    \begin{aligned}
    \CS_{J1-}^2(\lambda)(f_2)=&-\sum_{\ell=1}^N\int_{0}^{b}\CF_{\xi'}^{-1}\Big[\big(-\varphi'_{+0}(y_N)+\varphi_{+0}(y_N)A\big)e^{-Ay_N}M_-(x_N)\sum_{m=1}^{N-1}\frac{R_{J\ell,0}^-i\xi_m}{B_{-}^2}\widehat{f_{2m\ell+}}(\xi',y_N)\\
    & \hspace{3cm} +\varphi_{+0}(y_N)e^{-Ay_N}M_{-}(x_N)R_{J\ell,0}^-\widehat{f_{2N\ell+}}(\xi',y_N)\Big]dy_N\\
    &-\sum_{\ell=1}^N\int_{-b}^0\CF_{\xi'}^{-1}\Big[\big(\varphi'_{-0}(y_N)+\varphi_{-0}(y_N)A\big)e^{Ay_N}M_-(x_N)\sum_{m=1}^{N-1}\frac{R_{J\ell,0}^-i\xi_m}{B_-^2}\widehat{f_{2m\ell-}}(\xi',y_N)\\
    &\hspace{3cm} -\varphi_{-0}(y_N)e^{Ay_N}M_-(x_N)R_{J\ell,0}^-\widehat{f_{2N\ell-}}(\xi',y_N)\Big]dy_N,
\end{aligned}
\end{equation*}

\begin{equation*}
\begin{aligned}
    \CS_{J1-}^3(\lambda)(f_3)=&\sum_{\ell=1}^N\int_{0}^{b}\CF_{\xi'}^{-1}\Big[\big(-\varphi'_{+0}(y_N)+\varphi_{+0}(y_N)A\big)e^{-Ay_N}M_-(x_N)\frac{R_{J\ell,1}^-\gamma_{1-}}{\mu_{-}B_{-}^2}\widehat{f_{3\ell+}}(\xi',y_N)\Big]dy_N\\
    &+\sum_{\ell=1}^N\int_{-b}^0\CF_{\xi'}^{-1}\Big[\big(\varphi'_{-0}(y_N)+\varphi_{-0}(y_N)A\big)e^{Ay_N}M_-(x_N)\frac{R_{J\ell,1}^-\gamma_{1-}}{\mu_{-}B_-^2}\widehat{f_{3\ell-}}(\xi',y_N)\Big]dy_N,
\end{aligned}
\end{equation*}

\begin{equation*}
\begin{aligned}
    \CS_{J1-}^4(\lambda)(f_4)=&-\sum_{\ell=1}^N\int_{0}^{b}\CF_{\xi'}^{-1}\Big[\varphi_{+0}(y_N)e^{-Ay_N}M_-(x_N)\frac{R_{J\ell,1}^-\gamma_{1-}\lambda^{1/2}}{\mu_{-}B_{-}^2}\widehat{f_{4N\ell+}}(\xi',y_N)\Big]dy_N\\
    &+\sum_{\ell=1}^N\int_{-b}^0\CF_{\xi'}^{-1}\Big[\varphi_{-0}(y_N)e^{Ay_N}M_-(x_N)\frac{R_{J\ell,1}^-\gamma_{1-}\lambda^{1/2}}{\mu_{-}B_-^2}\widehat{f_{4N\ell-}}(\xi',y_N)\Big]dy_N,
\end{aligned}
\end{equation*}

\begin{equation*}
    \begin{aligned}
    \CS_{J1-}^5(\lambda)(f_5)=&-\sum_{\ell=1}^N\sum_{m=1}^{N-1}\int_{0}^{b}\CF_{\xi'}^{-1}\Big[\big(-\varphi'_{+0}(y_N)+\varphi_{+0}(y_N)A\big)e^{-Ay_N}M_-(x_N)\frac{R_{J\ell,1}^-}{B_{-}^2}\widehat{f_{5mm\ell+}}(\xi',y_N)\\
    &\hspace{4cm} -\varphi_{+0}(y_N)e^{-Ay_N}M_{-}(x_N)\frac{R_{J\ell,1}^-i\xi_m}{B_-^2}\widehat{f_{5mN\ell+}}(\xi',y_N)\Big]dy_N\\
    &-\sum_{\ell=1}^N\sum_{m=1}^{N-1}\int_{-b}^0\CF_{\xi'}^{-1}\Big[\big(\varphi'_{-0}(y_N)+\varphi_{-0}(y_N)A\big)e^{Ay_N}M_-(x_N)\frac{R_{J\ell,1}^-}{B_-^2}\widehat{f_{5mm\ell-}}(\xi',y_N)\\
    &\hspace{4cm} +\varphi_{-0}(y_N)e^{Ay_N}M_-(x_N)\frac{R_{J\ell,1}^-i\xi_m}{B_-^2}\widehat{f_{5mN\ell-}}(\xi',y_N)\Big]dy_N,
\end{aligned}
\end{equation*}

\begin{equation*}
\begin{aligned}
    \CS_{J2-}^1(\lambda)(f_1)=&\sum_{\ell=1}^N\int_{0}^{b}\CF_{\xi'}^{-1}\Big[\big(-\varphi'_{+0}(y_N)+\varphi_{+0}(y_N)A\big)e^{-Ay_N}e^{B_-x_N}\frac{S_{J\ell,-1}^-\gamma_{1-}\lambda^{1/2}}{\mu_{-}B_{-}^2}\widehat{f_{1\ell+}}(\xi',y_N)\Big]dy_N\\
    &+\sum_{\ell=1}^N\int_{-b}^0\CF_{\xi'}^{-1}\Big[\big(\varphi'_{-0}(y_N)+\varphi_{-0}(y_N)A\big)e^{Ay_N}e^{B_-x_N}\frac{S_{J\ell,-1}^-\gamma_{1-}\lambda^{1/2}}{\mu_{-}B_-^2}\widehat{f_{1\ell-}}(\xi',y_N)\Big]dy_N,
\end{aligned}
\end{equation*}

\begin{equation*}
    \begin{aligned}
    \CS_{J2-}^2(\lambda)(f_2)=&-\sum_{\ell=1}^N\int_{0}^{b}\CF_{\xi'}^{-1}\Big[\big(-\varphi'_{+0}(y_N)+\varphi_{+0}(y_N)A\big)e^{-Ay_N}e^{B_-x_N}\sum_{m=1}^{N-1}\frac{S_{J\ell,-1}^-i\xi_m}{B_{-}^2}\widehat{f_{2m\ell+}}(\xi',y_N)\\
    &\hspace{3cm} +\varphi_{+0}(y_N)e^{-Ay_N}e^{B_-x_N}S_{J\ell,-1}^-\widehat{f_{2N\ell+}}(\xi',y_N)\Big]dy_N\\
    &-\sum_{\ell=1}^N\int_{-b}^0\CF_{\xi'}^{-1}\Big[\big(\varphi'_{-0}(y_N)+\varphi_{-0}(y_N)A\big)e^{Ay_N}e^{B_-x_N}\sum_{m=1}^{N-1}\frac{S_{J\ell,-1}^-i\xi_m}{B_-^2}\widehat{f_{2m\ell-}}(\xi',y_N)\\
    &\hspace{3cm} -\varphi_{-0}(y_N)e^{Ay_N}e^{B_-x_N}S_{J\ell,-1}^-\widehat{f_{2N\ell-}}(\xi',y_N)\Big]dy_N,
\end{aligned}
\end{equation*}

\begin{equation*}
\begin{aligned}
    \CS_{J2-}^3(\lambda)(f_3)=& \sum_{\ell=1}^N\int_{0}^{b}\CF_{\xi'}^{-1}\Big[\big(-\varphi'_{+0}(y_N)+\varphi_{+0}(y_N)A\big)e^{-Ay_N}e^{B_-x_N}\frac{S_{J\ell,0}^-\gamma_{1-}}{\mu_{-}B_{-}^2}\widehat{f_{3\ell+}}(\xi',y_N)\Big]dy_N\\
    &+\sum_{\ell=1}^N\int_{-b}^0\CF_{\xi'}^{-1}\Big[\big(\varphi'_{-0}(y_N)+\varphi_{-0}(y_N)A\big)e^{Ay_N}e^{B_-x_N}\frac{S_{J\ell,0}^-\gamma_{1-}}{\mu_{-}B_-^2}\widehat{f_{3\ell-}}(\xi',y_N)\Big]dy_N,
\end{aligned}
\end{equation*}

\begin{equation*}
\begin{aligned}
    {\CS_{J2-}^4(\lambda)} (f_4)=&-\sum_{\ell=1}^N\int_{0}^{b}\CF_{\xi'}^{-1}\Big[\varphi_{+0}(y_N)e^{-Ay_N}e^{B_-x_N}\frac{S_{J\ell,0}^-\gamma_{1-}\lambda^{1/2}}{\mu_{-}B_{-}^2}\widehat{f_{4N\ell+}}(\xi',y_N)\Big]dy_N\\
    &+\sum_{\ell=1}^N\int_{-b}^{0}\CF_{\xi'}^{-1}\Big[\varphi_{-0}(y_N)e^{Ay_N}e^{B_-x_N}\frac{S_{J\ell,0}^-\gamma_{1-}\lambda^{1/2}}{\mu_{-}B_-^2}\widehat{f_{4N\ell-}}(\xi',y_N)\Big]dy_N,
\end{aligned}
\end{equation*}

\begin{equation*}
    \begin{aligned}
    \CS_{J2-}^5(\lambda)(f_5)=&-\sum_{\ell=1}^N\sum_{m=1}^{N-1}\int_{0}^{b}\CF_{\xi'}^{-1}
    \Big[\big(-\varphi'_{+0}(y_N)+\varphi_{+0}(y_N)A\big)e^{-Ay_N}e^{B_-x_N}\frac{S_{J\ell,0}^-}{B_{-}^2}\widehat{f_{5mm\ell+}}(\xi',y_N)\\
    &\hspace{4cm} -\varphi_{+0}(y_N)e^{-Ay_N}e^{B_-x_N}\frac{S_{J\ell,0}^-i\xi_m}{B_-^2}\widehat{f_{5mN\ell+}}(\xi',y_N)\Big]dy_N\\
    &-\sum_{\ell=1}^N\sum_{m=1}^{N-1}\int_{-b}^0\CF_{\xi'}^{-1}\Big[\big(\varphi'_{-0}(y_N)+\varphi_{-0}(y_N)A\big)e^{Ay_N}e^{B_-x_N}\frac{S_{J\ell,0}^-}{B_-^2}\widehat{f_{5mm\ell-}}(\xi',y_N)\\
    &\hspace{4cm} +\varphi_{-0}(y_N)e^{Ay_N}e^{B_-x_N}\frac{S_{J\ell,0}^-i\xi_m}{B_-^2}\widehat{f_{5mN\ell-}}(\xi',y_N)\Big]dy_N,
\end{aligned}
\end{equation*}

\begin{equation*}
\begin{aligned}
    \CS_{J3-}^1(\lambda)(f_1)=&\sum_{\ell=1}^N\int_{0}^{b}\CF_{\xi'}^{-1}\Big[\big(-\varphi'_{+0}(y_N)+\varphi_{+0}(y_N)A\big)e^{-Ay_N}M_-(-x_N-b)\frac{T_{J\ell,0}^-\gamma_{1-}\lambda^{1/2}}{\mu_{-}B_{-}^2}\widehat{f_{1\ell+}}(\xi',y_N)\Big]dy_N\\
    &+\sum_{\ell=1}^N\int_{-b}^0\CF_{\xi'}^{-1}\Big[\big(\varphi'_{-0}(y_N)+\varphi_{-0}(y_N)A\big)e^{Ay_N}M_-(-x_N-b)\frac{T_{J\ell,0}^-\gamma_{1-}\lambda^{1/2}}{\mu_{-}B_-^2}\widehat{f_{1\ell-}}(\xi',y_N)\Big]dy_N,
\end{aligned}
\end{equation*}

\begin{equation*}
    \begin{aligned}
    \CS_{J3-}^2(\lambda)(f_2)=&-\sum_{\ell=1}^N\int_{0}^{b}\CF_{\xi'}^{-1}\Big[\big(-\varphi'_{+0}(y_N)+\varphi_{+0}(y_N)A\big)e^{-Ay_N}M_-(-x_N-b)\sum_{m=1}^{N-1}\frac{T_{J\ell,0}^-i\xi_m}{B_{-}^2}\widehat{f_{2m\ell+}}(\xi',y_N)\\
    &\hspace{3cm} +\varphi_{+0}(y_N)e^{-Ay_N}M_{-}(-x_N-b)T_{J\ell,0}^-\widehat{f_{2N\ell+}}(\xi',y_N)\Big]dy_N\\
    &-\sum_{\ell=1}^N\int_{-b}^0\CF_{\xi'}^{-1}\Big[\big(\varphi'_{-0}(y_N)+\varphi_{-0}(y_N)A\big)e^{Ay_N}M_-(-x_N-b)\sum_{m=1}^{N-1}\frac{T_{J\ell,0}^-i\xi_m}{B_-^2}\widehat{f_{2m\ell-}}(\xi',y_N)\\
    &\hspace{3cm} -\varphi_{-0}(y_N)e^{Ay_N}M_-(-x_N-b)T_{J\ell,0}^-\widehat{f_{2N\ell-}}(\xi',y_N)\Big]dy_N,
\end{aligned}
\end{equation*}

\begin{equation*}
\begin{aligned}
    \CS_{J3-}^3(\lambda)(f_3)=& \sum_{\ell=1}^N\int_{0}^{b}\CF_{\xi'}^{-1}\Big[\big(-\varphi'_{+0}(y_N)+\varphi_{+0}(y_N)A\big)e^{-Ay_N}M_-(-x_N-b)\frac{T_{J\ell,1}^-\gamma_{1-}}{\mu_{-}B_{-}^2}\widehat{f_{3\ell+}}(\xi',y_N)\Big]dy_N\\
    &+\sum_{\ell=1}^N\int_{-b}^0\CF_{\xi'}^{-1}\Big[\big(\varphi'_{-0}(y_N)+\varphi_{-0}(y_N)A\big)e^{Ay_N}M_-(x_N)\frac{T_{J\ell,1}^-\gamma_{1-}}{\mu_{-}B_-^2}\widehat{f_{3\ell-}}(\xi',y_N)\Big]dy_N,
\end{aligned}
\end{equation*}

\begin{equation*}
\begin{aligned}
    \CS_{J3-}^4(\lambda)(f_4)=&-\sum_{\ell=1}^N\int_{0}^{b}\CF_{\xi'}^{-1}\Big[\varphi_{+0}(y_N)e^{-Ay_N}M_-(-x_N-b)\frac{T_{J\ell,1}^-\gamma_{1-}\lambda^{1/2}}{\mu_{-}B_{-}^2}\widehat{f_{4N\ell+}}(\xi',y_N)\Big]dy_N\\
    &+\sum_{\ell=1}^N\int_{-b}^0\CF_{\xi'}^{-1}\Big[\varphi_{-0}(y_N)e^{Ay_N}M_-(-x_N-b)\frac{T_{J\ell,1}^-\gamma_{1-}\lambda^{1/2}}{\mu_{-}B_-^2}\widehat{f_{4N\ell-}}(\xi',y_N)\Big]dy_N,
\end{aligned}
\end{equation*}

\begin{equation*}
    \begin{aligned}
    \CS_{J3-}^5(\lambda)(f_5)=&-\sum_{\ell=1}^N\sum_{m=1}^{N-1}\int_{0}^{b}\CF_{\xi'}^{-1}\Big[\big(-\varphi'_{+0}(y_N)+\varphi_{+0}(y_N)A\big)e^{-Ay_N}M_-(-x_N-b)\frac{T_{J\ell,1}^-}{B_{-}^2}\widehat{f_{5mm\ell+}}(\xi',y_N)\\
    &\hspace{4cm}-\varphi_{+0}(y_N)e^{-Ay_N}M_{-}(-x_N-b)\frac{T_{J\ell,1}^-i\xi_m}{B_-^2}\widehat{f_{5mN\ell+}}(\xi',y_N)\Big]dy_N\\
    &-\sum_{\ell=1}^N\sum_{m=1}^{N-1}\int_{-b}^0\CF_{\xi'}^{-1}\Big[\big(\varphi'_{-0}(y_N)+\varphi_{-0}(y_N)A\big)e^{Ay_N}M_-(-x_N-b)\frac{T_{J\ell,1}^-}{B_-^2}\widehat{f_{5mm\ell-}}(\xi',y_N)\\
    &\hspace{4cm} +\varphi_{-0}(y_N)e^{Ay_N}M_-(-x_N-b)\frac{T_{J\ell,1}^-i\xi_m}{B_-^2}\widehat{f_{5mN\ell-}}(\xi',y_N)\Big]dy_N,
\end{aligned}
\end{equation*}

\begin{equation*}
\begin{aligned}
    \CS_{J4-}^1(\lambda)(f_1)=&\sum_{\ell=1}^N\int_{0}^{b}\CF_{\xi'}^{-1}\Big[\big(-\varphi'_{+0}(y_N)+\varphi_{+0}(y_N)A\big)e^{-Ay_N}e^{-B_-(x_N+b)}\frac{S_{J\ell,-1}^b\gamma_{1-}\lambda^{1/2}}{\mu_{-}B_{-}^2}\widehat{f_{1\ell+}}(\xi',y_N)\Big]dy_N\\
    &+\sum_{\ell=1}^N\int_{-b}^0\CF_{\xi'}^{-1}\Big[\big(\varphi'_{-0}(y_N)+\varphi_{-0}(y_N)A\big)e^{Ay_N}e^{-B_-(x_N+b)}\frac{S_{J\ell,-1}^b\gamma_{1-}\lambda^{1/2}}{\mu_{-}B_-^2}\widehat{f_{1\ell-}}(\xi',y_N)\Big]dy_N,
\end{aligned}
\end{equation*}

\begin{equation*}
    \begin{aligned}
    \CS_{J4-}^2(\lambda)(f_2)=&-\sum_{\ell=1}^N\int_{0}^{b}\CF_{\xi'}^{-1}\Big[\big(-\varphi'_{+0}(y_N)+\varphi_{+0}(y_N)A\big)e^{-Ay_N}e^{-B_-(x_N+b)}\sum_{m=1}^{N-1}\frac{S_{J\ell,-1}^bi\xi_m}{B_{-}^2}\widehat{f_{2m\ell+}}(\xi',y_N)\\
    &\hspace{3cm} +\varphi_{+0}(y_N)e^{-Ay_N}e^{-B_-(x_N+b)}S_{J\ell,-1}^b\widehat{f_{2N\ell+}}(\xi',y_N)\Big]dy_N\\
    &-\sum_{\ell=1}^N\int_{-b}^0\CF_{\xi'}^{-1}\Big[\big(\varphi'_{-0}(y_N)+\varphi_{-0}(y_N)A\big)e^{Ay_N}e^{-B_-(x_N+b)}\sum_{m=1}^{N-1}\frac{S_{J\ell,-1}^bi\xi_m}{B_-^2}\widehat{f_{2m\ell-}}(\xi',y_N)\\
    &\hspace{3cm} -\varphi_{-0}(y_N)e^{Ay_N}e^{-B_-(x_N+b)}S_{J\ell,-1}^b\widehat{f_{2N\ell-}}(\xi',y_N)\Big]dy_N,
\end{aligned}
\end{equation*}

\begin{equation*}
\begin{aligned}
    \CS_{J4-}^3(\lambda)(f_3)=&\sum_{\ell=1}^N\int_{0}^{b}\CF_{\xi'}^{-1}\Big[\big(-\varphi'_{+0}(y_N)+\varphi_{+0}(y_N)A\big)e^{-Ay_N}e^{-B_-(x_N+b)}\frac{S_{J\ell,0}^b\gamma_{1-}}{\mu_{-}B_{-}^2}\widehat{f_{3\ell+}}(\xi',y_N)\Big]dy_N\\
    &+\sum_{\ell=1}^N\int_{-b}^0\CF_{\xi'}^{-1}\Big[\big(\varphi'_{-0}(y_N)+\varphi_{-0}(y_N)A\big)e^{Ay_N}e^{-B_-(x_N+b)}\frac{S_{J\ell,0}^b\gamma_{1-}}{\mu_{-}B_-^2}\widehat{f_{3\ell-}}(\xi',y_N)\Big]dy_N,
\end{aligned}
\end{equation*}

\begin{equation*}
\begin{aligned}
    \CS_{J4-}^4(\lambda)  (f_4)=&-\sum_{\ell=1}^N\int_{0}^{b}\CF_{\xi'}^{-1}\Big[\varphi_{+0}(y_N)e^{-Ay_N}e^{-B_-(x_N+b)}\frac{S_{J\ell,0}^b\gamma_{1-}\lambda^{1/2}}{\mu_{-}B_{-}^2}\widehat{f_{4N\ell+}}(\xi',y_N)\Big]dy_N\\
    &+\sum_{\ell=1}^N\int_{-b}^0\CF_{\xi'}^{-1}\Big[\varphi_{-0}(y_N)e^{Ay_N} e^{-B_-(x_N+b)}\frac{S_{J\ell,0}^b\gamma_{1-}\lambda^{1/2}}{\mu_{-}B_-^2}\widehat{f_{4N\ell-}}(\xi',y_N)\Big]dy_N,
\end{aligned}
\end{equation*}

\begin{equation*}
    \begin{aligned}
    \CS_{J4-}^5(\lambda)(f_5)=&-\sum_{\ell=1}^N\sum_{m=1}^{N-1}\int_{0}^{b}\CF_{\xi'}^{-1}\Big[\big(-\varphi'_{+0}(y_N)+\varphi_{+0}(y_N)A\big)e^{-Ay_N}e^{-B_-(x_N+b)}\frac{S_{J\ell,0}^b}{B_{-}^2}\widehat{f_{5mm\ell+}}(\xi',y_N)\\
    &\hspace{4cm} -\varphi_{+0}(y_N)e^{-Ay_N}e^{-B_-(x_N+b)}\frac{S_{J\ell,0}^bi\xi_m}{B_-^2}\widehat{f_{5mN\ell+}}(\xi',y_N)\Big]dy_N\\
    &-\sum_{\ell=1}^N\sum_{m=1}^{N-1}\int_{-b}^0\CF_{\xi'}^{-1}\Big[\big(\varphi'_{-0}(y_N)+\varphi_{-0}(y_N)A\big)e^{Ay_N}e^{-B_-(x_N+b)}\frac{S_{J\ell,0}^b}{B_-^2}\widehat{f_{5mm\ell-}}(\xi',y_N)\\
    &\hspace{4cm} +\varphi_{-0}(y_N)e^{Ay_N}e^{-B_-(x_N+b)}\frac{S_{J\ell,0}^bi\xi_m}{B_-^2}\widehat{f_{5mN\ell-}}(\xi',y_N)\Big]dy_N.
\end{aligned}
\end{equation*}

To bound the operators $\CS_{ Ji-}(\lambda),$ $i=1,\dots,4,$ we note from Corollary \ref{cor1}  that the multipliers 
\begin{gather*}
	\frac{R_{J\ell,0}^-\gamma_{1-}\lambda^{1/2}}{\mu_{-}B_-^3},\,\, 
    \frac{R_{J\ell,0}^-i\xi_m}{B_-^3},\,\, \frac{R_{J\ell,1}^-}{B_-^3},\,\,
	\frac{R_{J\ell,0}^-}{AB_-},\,\, 
    \frac{R_{J\ell,1}^-\gamma_{1-}}{\mu_{-}AB_-^3},\,\, 
    \frac{R_{J\ell,1}^-i\xi_m}{AB_-^3},\\
	\frac{R_{J\ell,0}^-\gamma_{1-}\lambda^{1/2}}{\mu_{-}AB_-^2},\,\, 
    \frac{R_{J\ell,0}^-i\xi_m}{AB_-^2},\,\, 
    \frac{R_{J\ell,1}^-}{AB_-^2},\,\,
    \frac{R_{J\ell,0}^-}{A^2},\,\, 
    \frac{R_{J\ell,1}^- \gamma_{1-}}{\mu_{-}A^2B_-^2},\,\, 
    \frac{R_{J\ell,1}^-i\xi_m}{A^2B_-^2},\\
	\frac{S_{J\ell,-1}^- \gamma_{1-}\lambda^{1/2}}{\mu_{-}B_-^2},\,\, 
    \frac{S_{J\ell,-1}^- i\xi_m}{B_-^2},\,\,
    \frac{S_{J\ell,0}^-}{B_-^2},\,\,
	\frac{S_{J\ell,-1}^-}{A},\,\, 
    \frac{S_{J\ell,0}^-\gamma_{1-}}{\mu_{-}AB_-^2},\\ 
    \frac{S_{J\ell,0}^- i\xi_m}{AB_-^2},\,\,
	\frac{S_{J\ell,-1}^- \gamma_{1-}\lambda^{1/2}}{\mu_{-}AB_-},\,\,
    \frac{S_{J\ell,-1}^- i\xi_m}{AB_-},\,\, 
    \frac{S_{J\ell,0}^-}{AB_-},\\
	\frac{T_{J\ell,0}^-\gamma_{1-}\lambda^{1/2}}{\mu_{-}B_-^3},\,\, 
    \frac{T_{J\ell,0}^-i\xi_m}{B_-^3},\,\, 
    \frac{T_{J\ell,1}^-}{B_-^3},\,\,
	\frac{T_{J\ell,0}^-}{AB_-},\,\, 
    \frac{T_{J\ell,1}^-\gamma_{1-}}{\mu_{-}AB_-^3},\,\,
    \frac{T_{J\ell,1}^-i\xi_m}{AB_-^3},\\
    \frac{T_{J\ell,0}^-\gamma_{1-}\lambda^{1/2}}{\mu_{-}AB_-^2},\,\, 
    \frac{T_{J\ell,0}^-i\xi_m}{AB_-^2},\,\, 
    \frac{T_{J\ell,1}^-}{AB_-^2},\,\,
     \frac{T_{J\ell,0}^-}{A^2},\,\, 
     \frac{T_{J\ell,1}^-\gamma_{1-}}{\mu_{-}A^2B_-^2},\,\, 
     \frac{T_{J\ell,1}^-i\xi_m}{A^2B_-^2},\\
	\frac{S_{J\ell,-1}^b\gamma_{1-}\lambda^{1/2}}{\mu_{-}B_-^2},\,\, 
    \frac{S_{J\ell,-1}^bi\xi_m}{B_-^2},\,\,
    \frac{S_{J\ell,0}^b}{B_-^2},\,\,
	\frac{S_{J\ell,-1}^b}{A},\,\, 
    \frac{S_{J\ell,0}^b \gamma_{1-}\lambda^{1/2}}{\mu_{-}AB_-^2},\\ 
    \frac{S_{J\ell,0}^bi\xi_m}{AB_-^2},\,\,
	\frac{S_{J\ell,-1}^b\gamma_{1-}}{\mu_{-}AB_-},\,\,
    \frac{S_{J\ell,-1}^bi\xi_m}{AB_-},\,\, 
    \frac{S_{J\ell,0}^b}{AB_-},
\end{gather*}
which appear in formulations of $\CS_{Ji-}^k(\lambda)$ above, belong to the class
$\bM_{-2,2}(\widetilde{\Gamma}_{\ep, \lambda_0}).$ Then Lemmas \ref{lem:4.5} and \ref{lem:+-} yield
\begin{align*}
	&\CR_{\CL\big(\CZ_q(\dot\Omega),H^{2-j}_q(\Omega_{-})^N\big)}\Big(\Big\{(\tau\pd_{\tau})^\ell\lambda^{j/2} \CS_{Ji-}(\lambda)\,\big|\,\lambda\in\Gamma_{\varepsilon,\lambda_0}\Big\}\Big)\leq r_b
\end{align*}
with $\ell=0,1, i=1,\dots,4, j=0,1,2$, where we have also used the embedding 
$L_q(\Omega_{+})\hookrightarrow L_q(\Omega_{+b}).$ This gives the $\CR$-bound of $\CS_{-}(\lambda)$ in \eqref{res:key_1}.
\medskip

At last, we give the formulation of $\CP_{-}^k(\lambda)$ and prove the $\CR$-bound of $\CP_{-}(\lambda)$ in \eqref{res:key_1}. According to \eqref{eq:pressure}, \eqref{conv:f} and Lemma \ref{lem:Vol-2}, we have 
\begin{equation*}
\begin{aligned}
    \CP_{-}^1(\lambda)(f_1)=&\sum_{\ell=1}^N\int_{0}^{b}\CF_{\xi'}^{-1}\Big[\big(-\varphi'_{+0}(y_N)+\varphi_{+0}(y_N)A\big)e^{-Ay_N}e^{Ax_N}\frac{p_{\ell,0}^0\gamma_{1-}\lambda^{1/2}}{\mu_{-}B_{-}^2}\widehat{f_{1\ell+}}(\xi',y_N)\Big]dy_N\\
    &+\sum_{\ell=1}^N\int_{-b}^0\CF_{\xi'}^{-1}\Big[\big(\varphi'_{-0}(y_N)+\varphi_{-0}(y_N)A\big)e^{Ay_N}e^{Ax_N}\frac{p_{\ell,0}^0\gamma_{1-}\lambda^{1/2}}{\mu_{-}B_-^2}\widehat{f_{1\ell-}}(\xi',y_N)\Big]dy_N\\
    &+\sum_{\ell=1}^N\int_{0}^{b}\CF_{\xi'}^{-1}\Big[\big(-\varphi'_{+0}(y_N)+\varphi_{+0}(y_N)A\big)e^{-Ay_N}e^{-A(x_N+b)}\frac{p_{\ell,0}^b\gamma_{1-}\lambda^{1/2}}{\mu_{-}B_{-}^2}\widehat{f_{1\ell+}}(\xi',y_N)\Big]dy_N\\
    &+\sum_{\ell=1}^N\int_{-b}^0\CF_{\xi'}^{-1}\Big[\big(\varphi'_{-0}(y_N)+\varphi_{-0}(y_N)A\big)e^{Ay_N}e^{-A(x_N+b)}\frac{p_{\ell,0}^b\gamma_{1-}\lambda^{1/2}}{\mu_{-}B_-^2}\widehat{f_{1\ell-}}(\xi',y_N)\Big]dy_N,
\end{aligned}
\end{equation*}

\begin{equation*}
    \begin{aligned}
    \CP_{-}^2(\lambda)(f_2)=&-\sum_{\ell=1}^N\int_{0}^{b}\CF_{\xi'}^{-1}\Big[\big(-\varphi'_{+0}(y_N)+\varphi_{+0}(y_N)A\big)e^{-Ay_N}e^{Ax_N}\sum_{m=1}^{N-1}\frac{p_{\ell,0}^0i\xi_m}{B_{-}^2}\widehat{f_{2m\ell+}}(\xi',y_N)\\
    &\qquad\qquad+\varphi_{+0}(y_N)e^{-Ay_N}e^{Ax_N}p_{\ell,0}^0\widehat{f_{2N\ell+}}(\xi',y_N)\Big]dy_N\\
    &-\sum_{\ell=1}^N\int_{-b}^0\CF_{\xi'}^{-1}\Big[\big(\varphi'_{-0}(y_N)+\varphi_{-0}(y_N)A\big)e^{Ay_N}e^{Ax_N}\sum_{m=1}^{N-1}\frac{p_{\ell,0}^0i\xi_m}{B_-^2}\widehat{f_{2m\ell-}}(\xi',y_N)\\
    &\qquad\qquad-\varphi_0(y_N)e^{Ay_N}e^{Ax_N}p_{J\ell,0}^0\widehat{f_{2N\ell-}}(\xi',y_N)\Big]dy_N\\
    &-\sum_{\ell=1}^N\int_{0}^{b}\CF_{\xi'}^{-1}\Big[\big(-\varphi'_{+0}(y_N)+\varphi_{+0}(y_N)A\big)e^{-Ay_N}e^{-A(x_N+b)}\sum_{m=1}^{N-1}\frac{p_{\ell,0}^bi\xi_m}{B_{-}^2}\widehat{f_{2m\ell+}}(\xi',y_N)\\
    &\qquad\qquad+\varphi_{+0}(y_N)e^{-Ay_N}e^{-A(x_N+b)}p_{\ell,0}^b\widehat{f_{2N\ell+}}(\xi',y_N)\Big]dy_N\\
    &-\sum_{\ell=1}^N\int_{-b}^0\CF_{\xi'}^{-1}\Big[\big(\varphi'_{-0}(y_N)+\varphi_{-0}(y_N)A\big)e^{Ay_N}e^{-A(x_N+b)}\sum_{m=1}^{N-1}\frac{p_{\ell,0}^bi\xi_m}{B_-^2}\widehat{f_{2m\ell-}}(\xi',y_N)\\
    &\qquad\qquad-\varphi_{-0}(y_N)e^{Ay_N}e^{-A(x_N+b)}p_{J\ell,0}^b\widehat{f_{2N\ell-}}(\xi',y_N)\Big]dy_N
\end{aligned}
\end{equation*}
\begin{equation*}
\begin{aligned}
    \CP_{-}^3(\lambda)(f_3)=&\sum_{\ell=1}^N\int_{0}^{b}\CF_{\xi'}^{-1}\Big[\big(-\varphi'_{+0}(y_N)+\varphi_{+0}(y_N)A\big)e^{-Ay_N}e^{Ax_N}\frac{p_{\ell,1}^0\gamma_{1-}}{\mu_{-}B_{-}^2}\widehat{f_{3\ell+}}(\xi',y_N)\Big]dy_N\\
    &+\sum_{\ell=1}^N\int_{-b}^0\CF_{\xi'}^{-1}\Big[\big(\varphi'_{-0}(y_N)+\varphi_{-0}(y_N)A\big)e^{Ay_N}e^{Ax_N}\frac{p_{\ell,1}^0\gamma_{1-}}{\mu_{-}B_-^2}\widehat{f_{3\ell-}}(\xi',y_N)\Big]dy_N\\
    &+\sum_{\ell=1}^N\int_{0}^{b}\CF_{\xi'}^{-1}\Big[\big(-\varphi'_{+0}(y_N)+\varphi_{+0}(y_N)A\big)e^{-Ay_N}e^{-A(x_N+b)}\frac{p_{\ell,1}^b\gamma_{1-}}{\mu_{-}B_{-}^2}\widehat{f_{3\ell+}}(\xi',y_N)\Big]dy_N\\
    &+\sum_{\ell=1}^N\int_{-b}^0\CF_{\xi'}^{-1}\Big[\big(\varphi'_{-0}(y_N)+\varphi_{-0}(y_N)A\big)e^{Ay_N}e^{-A(x_N+b)}\frac{p_{\ell,1}^b\gamma_{1-}}{\mu_{-}B_-^2}\widehat{f_{3\ell-}}(\xi',y_N)\Big]dy_N,
\end{aligned}
\end{equation*}
\begin{equation*}
\begin{aligned}
    \CP_{-}^4(\lambda)(f_4)=&-\sum_{\ell=1}^N\int_{0}^{b}\CF_{\xi'}^{-1}\Big[\varphi_{+0}(y_N)e^{-Ay_N}e^{Ax_N}\frac{p_{\ell,1}^0\gamma_{1-}\lambda^{1/2}}{\mu_{-}B_{-}^2}\widehat{f_{4N\ell+}}(\xi',y_N)\Big]dy_N\\
    &+\sum_{\ell=1}^N\int_{-b}^0\CF_{\xi'}^{-1}\Big[\varphi_{-0}(y_N)e^{Ay_N}e^{Ax_N}\frac{p_{\ell,1}^0\gamma_{1-}\lambda^{1/2}}{\mu_{-}B_-^2}\widehat{f_{4N\ell-}}(\xi',y_N)\Big]dy_N\\
    &-\sum_{\ell=1}^N\int_{0}^{b}\CF_{\xi'}^{-1}\Big[\varphi_{+0}(y_N)e^{-Ay_N}e^{-A(x_N+b)}\frac{p_{\ell,1}^b\gamma_{1-}\lambda^{1/2}}{\mu_{-}B_{-}^2}\widehat{f_{4N\ell+}}(\xi',y_N)\Big]dy_N\\
    &+\sum_{\ell=1}^N\int_{-b}^0\CF_{\xi'}^{-1}\Big[\varphi_{-0}(y_N)e^{Ay_N}e^{-A(x_N+b)}\frac{p_{\ell,1}^b\gamma_{1-}\lambda^{1/2}}{\mu_{-}B_-^2}\widehat{f_{4N\ell-}}(\xi',y_N)\Big]dy_N,
\end{aligned}
\end{equation*}
\begin{equation*}
    \begin{aligned}
    \CP_{-}^5(\lambda)(f_5)=&-\sum_{\ell=1}^N\int_{0}^{b}\CF_{\xi'}^{-1}\Big[\big(-\varphi'_{+0}(y_N)+\varphi_{+0}(y_N)A\big)e^{-Ay_N}e^{Ax_N}\sum_{m=1}^{N-1}\frac{p_{\ell,1}^0}{B_{-}^2}\widehat{f_{5mm\ell+}}(\xi',y_N)\\
    &\qquad\qquad-\varphi_{+0}(y_N)e^{-Ay_N}e^{Ax_N}\sum_{m=1}^{N-1}\frac{p_{\ell,1}^0i\xi_m}{B_-^2}\widehat{f_{5mN\ell+}}(\xi',y_N)\Big]dy_N\\
    &-\sum_{\ell=1}^N\int_{-b}^0\CF_{\xi'}^{-1}\Big[\big(\varphi'_{-0}(y_N)+\varphi_{-0}(y_N)A\big)e^{Ay_N}e^{Ax_N}\sum_{m=1}^{N-1}\frac{p_{\ell,1}^0}{B_-^2}\widehat{f_{5mm\ell-}}(\xi',y_N)\\
    &\qquad\qquad+\varphi_{-0}(y_N)e^{Ay_N}e^{Ax_N}\sum_{m=1}^{N-1}\frac{p_{\ell,1}^0i\xi_m}{B_-^2}\widehat{f_{5mN\ell-}}(\xi',y_N)\Big]dy_N\\
    &-\sum_{\ell=1}^N\int_{0}^{b}\CF_{\xi'}^{-1}\Big[\big(-\varphi'_{+0}(y_N)+\varphi_{+0}(y_N)A\big)e^{-Ay_N}e^{-A(x_N+b)}\sum_{m=1}^{N-1}\frac{p_{\ell,1}^b}{B_{-}^2}\widehat{f_{5mm\ell+}}(\xi',y_N)\\
    &\qquad\qquad-\varphi_{+0}(y_N)e^{-Ay_N}e^{-A(x_N+b)}\sum_{m=1}^{N-1}\frac{p_{\ell,1}^bi\xi_m}{B_-^2}\widehat{f_{5mN\ell+}}(\xi',y_N)\Big]dy_N\\
    &-\sum_{\ell=1}^N\int_{-b}^0\CF_{\xi'}^{-1}\Big[\big(\varphi'_{-0}(y_N)+\varphi_{-0}(y_N)A\big)e^{Ay_N}e^{-A(x_N+b)}\sum_{m=1}^{N-1}\frac{p_{\ell,1}^b}{B_-^2}\widehat{f_{5mm\ell-}}(\xi',y_N)\\
    &\qquad\qquad+\varphi_{-0}(y_N)e^{Ay_N}e^{-A(x_N+b)}\sum_{m=1}^{N-1}\frac{p_{\ell,1}^bi\xi_m}{B_-^2}\widehat{f_{5mN\ell-}}(\xi',y_N)\Big]dy_N.
\end{aligned}
\end{equation*}

Finally, for the multipliers appearing in $\CP^-(\lambda),$ we have 
\begin{gather*}
	\frac{p_{\ell,0}^zA\gamma_{1-}\lambda^{1/2}}{\mu_{-}B_-^2},\,\, \frac{p_{\ell,0}^zAi\xi_m}{B_-^2},\,\, \frac{p_{\ell,1}^zA}{B_-^2},\,\,
	{p_{\ell,0}^z},\,\, \frac{p_{\ell,1}^z\gamma_{1-}\lambda^{1/2}}{\mu_{-}B_-^2},\,\, \frac{p_{\ell,1}^zi\xi_m}{B_-^2}\in\bM_{0,2}(\widetilde{\Gamma}_{\ep, \lambda_0})
\end{gather*}
with $z\in\{0,b\}.$
Then Lemmas \ref{lem:4.5} and \ref{lem:+-2} furnish that
\begin{align*}
	&\CR_{\CL\big(\CZ_q(\dot\Omega),L_q(\Omega_{-})^{N}\big)}\Big(\Big\{(\tau\pd_{\tau})^\ell \nabla \CP_{-}(\lambda)\,\big|\,\lambda\in\Gamma_{\varepsilon,\lambda_0}\Big\}\Big)\leq r_b\quad (\ell=0,1).
\end{align*}
This completes the proof of Theorem \ref{thm:5.1}.


\section{Resolvent operators of the two-phase problem}
\label{sec:resolvent}
In this section, we consider the resolvent problem \eqref{eq:res_full} and the main result on the problem \eqref{eq:res_full} reads as follows.
\begin{thm}\label{thm:main2}
	Assume that $0<\varepsilon<\pi/2$, $1<q<\infty$ and the constants $\mu_{\pm},$ $\nu_+,$ $\gamma_{1\pm},$ $\gamma_{2+}>0.$
	Let us introduce the functional spaces for the given terms in \eqref{eq:res_full} as follows 
	\begin{equation*}
		\begin{aligned}
			X_q (\dot \Omega)&=\big\{ (f_1,\Bf_2,\Bf_3,f_4,\Bf_5,\Bf_6,\Bf_7) :
			f_1\in H^1_q(\Omega_+),
			\Bf_2\in L_q(\Omega_+)^N, \Bf_3,\Bf_5\in L_q(\Omega_-)^N,\\
			& \hspace{5.5cm} f_4\in H^1_q(\Omega_-), 
			\Bf_6\in H^1_q(\dot\Omega)^N, \Bf_7\in H^2_q(\dot\Omega)^N, 
			f_4=\dv \Bf_5 \big\}, \\
			\CX_q(\dot\Omega)&=\big\{ (f_1,\Bf_2,\Bf_3,g_4, f_4,\Bf_5,\bg_6,\Bf_6,\bh_7,\bg_7,\Bf_7) :
			(f_1,\Bf_2,\Bf_3,f_4,\Bf_5,\Bf_6,\Bf_7) \in X_q(\dot \Omega),\\
			& \hspace{6cm} g_4\in L_q(\Omega_-), \,\Bg_6\in L_q(\dot\Omega)^N,\,\bg_7\in H^1_q(\dot\Omega)^N,\,\bh_7\in L_q(\dot\Omega)^N\big\}.
		\end{aligned}
	\end{equation*} 
	In addition, we also write for short that 
	\begin{equation*}
		\begin{aligned}
			\bF&=(f_+, \bg_+,\bg_-,g_d, \fg_d,\bh,\bk),\\
			\bF_{\lambda}&=(f_+,\bg_+,\bg_-,\lambda^{1/2} g_d, g_d, \lambda \fg_d, 
			\lambda^{1/2} \bh,\bh,\lambda\bk,\lambda^{1/2}\bk,\bk).
		\end{aligned}
	\end{equation*}
	For any $\bF\in X_q(\dot\Omega),$ there exist constants $\lambda_0, r_b$ and the operator families
	\begin{equation*}
		\begin{aligned}
			\Theta_+(\lambda) &\in \Hol\Big( \Sigma_{\ep, \lambda_0}  \cap K_{\ep} ; 
			\CL\big( \CX_q(\dot \Omega) ; H^1_q(\Omega_+) \big) \Big),\\
			\CA_\pm(\lambda) &\in \Hol\Big( \Sigma_{\ep, \lambda_0}  \cap K_{\ep} ; 
			\CL\big( \CX_q(\dot \Omega) ;H^2_q(\Omega_\pm)^N \big) \Big),\\
			\CB_-(\lambda) &\in \Hol\Big( \Sigma_{\ep, \lambda_0}  \cap K_{\ep}; 
			\CL\big( \CX_q(\dot \Omega) ;\wh H^1_q(\Omega_-) \big) \Big)
		\end{aligned}
	\end{equation*}
	such that $(\rho_+, \bv_+,\bv_-,\fp_-)
	= \big(\Theta_+(\lambda) , \CA_+(\lambda), \CA_-(\lambda),\CB_-(\lambda)\big) \,\bF_{\lambda}$ 
	solves \eqref{eq:res_full} uniquely.  Here the regions $\Sigma_{\ep, \lambda_0}$ and $K_{\ep}$ are defined by \eqref{eq:KSigma}. Moreover, we have 
	\begin{gather*}
		\CR_{\CL\big(\CX_q(\dot \Omega); H^{1}_q(\Omega_+) \big)}
		\Big( \Big\{ (\tau \pa_{\tau})^{\ell}\big( \lambda^{k}
		\Theta_+(\lambda)\big) : 
		\lambda \in \Sigma_{\ep, \lambda_0} \cap K_{\ep} \Big\}\Big) \leq r_b,\\
		\CR_{\CL\big(\CX_q(\dot\Omega); H^{2-j}_q(\Omega_\pm)^N \big)}
		\Big( \Big\{ (\tau \pa_{\tau})^{\ell}\big( \lambda^{j\slash 2}
		\CA_\pm(\lambda)\big) : 
		\lambda \in \Sigma_{\ep, \lambda_0} \cap K_{\ep} \Big\}\Big) \leq r_b,\\
		\CR_{\CL\big(\CX_q(\dot\Omega); L_q(\Omega_-)^N \big)}
		\Big( \Big\{ (\tau \pa_{\tau})^{\ell}   \nabla \CB_-(\lambda) : 
		\lambda \in \Sigma_{\ep, \lambda_0} \cap K_{\ep} \Big\}\Big) \leq r_b
	\end{gather*}
	for $k,\ell =0,1,$ $j=0,1,2,$ and $\tau = \Im \lambda.$ 
	Above the choice of  $r_b$ depends solely on 
	$\varepsilon,$ $q,$ $N,$ $\gamma_{1\pm},$ $\gamma_{2+},$ $\mu_{\pm},$ 
	$\nu_+,$  $\lambda_0$, $\delta_0$ and $b.$
\end{thm}

To prove Theorem \ref{thm:main2}, it suffices to study the following reduced problem:
\begin{equation}\label{eq:rr}
	\left\{ \begin{aligned}
		&\gamma_{1+}\lambda\bv_+-\Di\bS_{\delta+}(\bv_+)=\bg_+
		&&\quad&\text{in}& \quad  \Omega_+, \\
		&\gamma_{1-}\lambda\bv_--\Di \bT_-(\bv_-,\fp_-)=\bg_-
		&&\quad&\text{in}& \quad  \Omega_-, \\ 
		&\di\bv_-= g_d=\di \fg_d
		&&\quad&\text{in}& \quad  \Omega_-, \\ 
		&\bS_{\delta+}(\bv_+)\bn
		-\bT_-(\bv_-,\fp_-)\bn =\bh
		&&\quad&\text{on}& \quad \Gamma,\\
		& \bv_+-\bv_-=\bk
		&&\quad&\text{on}& \quad \Gamma ,	\\
		&\bv_-= 0
		&&\quad&\text{on}& \quad S
	\end{aligned}
	\right.
\end{equation}
with $\bS_{\delta+}(\bv_+)$ defined in \eqref{eq:S_delta}.
Concerning the problem \eqref{eq:rr}, it is not hard to prove the following result by using Theorem \ref{thm:5.1}.
\begin{thm}\label{thm:rr}
Assume that $0<\varepsilon<\pi/2$, $1<q<\infty$ and the constants $\mu_{\pm},$ $\nu_+,$ $\gamma_{1\pm},$ $\gamma_{2+}>0.$
With the spaces $Z_q(\dot\Omega)$ and $\CZ_q(\dot\Omega)$ introduced by Theorem \ref{thm:5.1}, we denote that 
	\begin{equation*}
		\begin{aligned}
			Y_q (\dot \Omega)&=\big\{ (\Bf_1,\Bf_2,f_3,\Bf_4,\Bf_5,\Bf_6) :
			\Bf_1\in L_q(\Omega_+)^N, \Bf_2,\Bf_4\in L_q(\Omega_-)^N,
			f_3\in H^1_q(\Omega_-),\\
			& \hspace{7.5cm} (\Bf_5,\Bf_6)\in Z_q(\dot\Omega), f_3=\dv \Bf_4 \big\}, \\
			\CY_q(\dot\Omega)&=\big\{ (\Bf_1,\Bf_2,g_3, f_3,\Bf_4,\bg_5,\Bf_5,\bh_6,\bg_6,\Bf_6) :
			(\Bf_1,\Bf_2,f_3,\Bf_4,\Bf_5,\Bf_6) \in Y_q(\dot \Omega),\\
			& \hspace{5.5cm}  g_3\in L_q(\Omega_-), \, 
			(\bg_5,\Bf_5,\bh_6,\bg_6,\Bf_6)\in \CZ_q(\dot\Omega)\big\}.
		\end{aligned}
	\end{equation*} 
	In addition, we also write for short that 
	\begin{equation*}
		\begin{aligned}
			\bG&=(\bg_+,\bg_-,g_d, \fg_d,\bh,\bk),\\
			\bG_{\lambda}&=(\bg_+,\bg_-,\lambda^{1/2} g_d, g_d, \lambda \fg_d, 
			\lambda^{1/2} \bh,\bh,\lambda\bk,\lambda^{1/2}\bk,\bk).
		\end{aligned}
	\end{equation*}
	For any $\bG\in Y_q(\dot \Omega),$ there exist constants $\lambda_0, r_b$ and the operator families
	\begin{equation*}
		\begin{aligned}
			\CA_\pm^0(\lambda) &\in \Hol\Big( \Gamma_{\ep,\lambda_0} ; 
			\CL\big( \CY_q(\dot \Omega) ;H^2_q(\Omega_\pm)^N \big) \Big),\\
			\CB_-^0(\lambda) &\in \Hol\Big( \Gamma_{\ep,\lambda_0} ; 
			\CL\big( \CY_q(\dot \Omega) ;\wh H^1_q(\Omega_-) \big) \Big)
		\end{aligned}
	\end{equation*}
	such that $(\bv_+,\bv_-,\fp_-)
	= \big(\CA_+^0(\lambda), \CA_-^0(\lambda),\CB_-^0(\lambda)\big) \,\bG_{\lambda}$
	solves \eqref{eq:rr} uniquely.  Here the set $\Gamma_{\ep,\lambda_0}$ is defined by \eqref{def:Gamma} for some $\delta_0>0.$ Moreover, we have 
	\begin{gather*}
		\CR_{\CL\big(\CY_q(\dot \Omega); H^{2-j}_q(\Omega_\pm)^N \big)}
		\Big( \Big\{ (\tau \pa_{\tau})^{\ell}\big( \lambda^{j\slash 2}
		\CA_\pm^0(\lambda)\big) : 
		\lambda \in \Gamma_{\ep,\lambda_0} \Big\}\Big) \leq r_b,\\
		\CR_{\CL\big(\CY_q(\dot \Omega); L_q(\Omega_-)^N \big)}
		\Big( \Big\{ (\tau \pa_{\tau})^{\ell}   \nabla \CB_-^0(\lambda) : 
		\lambda \in \Gamma_{\ep,\lambda_0} \Big\}\Big) \leq r_b
	\end{gather*}
	for $\ell =0,1,$ $j=0,1,2,$ and $\tau = \Im \lambda.$ 
	Above the choice of  $r_b$ depends solely on 
	$\varepsilon,$ $q,$ $N,$ $\gamma_{1\pm},$ $\gamma_{2+},$ $\mu_{\pm},$ 
	$\nu_+,$  $\lambda_0$, $\delta_0$ and $b.$
\end{thm}

\begin{rmk}
It is not hard to prove Theorem \ref{thm:main2} by imposing $\delta=\gamma_{1+}\gamma_{2+}\lambda^{-1}$ in \eqref{eq:rr} and taking advantage of Theorem \ref{thm:rr}. In addition, the uniqueness part of \eqref{eq:rr} follows from the existence part as one may consider the dual problem of \eqref{eq:rr} in the functional spaces with respect to the conjugate index $q'=q/(q-1)$.
\end{rmk}

To prove Theorem \ref{thm:rr}, we recall some results on the one-phase flows in \cite{GS2014,Saito2015}.
\begin{prop}[{\cite[Theorem 2.3]{GS2014}}]
\label{prop:v+1}
	Assume that $0<\varepsilon<\pi/2$, $1<q<\infty$ and the constants $\mu_{+},$ $\nu_+,$ $\gamma_{1+},$ $\gamma_{2+}>0.$
	For any $\bg_+ \in  L_q(\Omega_+)^N,$ there exist constants $\lambda_0, r_b$ and the operator family
	\begin{equation*}
		\CA_+^1(\lambda) \in \Hol\Big( \Gamma_{\ep,\lambda_0} ; 
		\CL\big( L_q(\Omega_+)^N;H^2_q(\Omega_+)^N \big) \Big)
	\end{equation*}
	such that $\Bv_+^1= \CA_+^1(\lambda) \, \bg_+$ solves the following system
	\begin{equation}\label{eq:rr_1}
		\left\{ \begin{aligned}
			&\gamma_{1+}\lambda\bv_+^1-\Di\bS_{\delta+}(\bv_+^1)=\bg_+
			&&\quad&\text{in}& \quad  \Omega_+, \\
			&\bS_{\delta+}(\bv_+^1)\bn_{\Gamma}=0
			&&\quad&\text{on}& \quad \Gamma.
		\end{aligned}
		\right.
	\end{equation} 
	Here the set $\Gamma_{\ep,\lambda_0}$ is defined by \eqref{def:Gamma} for some $\delta_0>0$ and $\bn_{\Gamma}=(0,\cdots,0,-1)^{\top}.$
	Moreover, we have 
	\begin{equation*}
		\CR_{\CL\big(L_q(\Omega_+)^N; H^{2-j}_q(\Omega_+)^N \big)}
		\Big( \Big\{ (\tau \pa_{\tau})^{\ell}\big( \lambda^{j\slash 2}\CA_+^1(\lambda)\big) : 
		\lambda \in \Gamma_{\ep,\lambda_0} \Big\}\Big) \leq r_b
	\end{equation*}
	for $\ell =0,1,$ $j=0,1,2,$ and $\tau = \Im \lambda.$ 
	Above the choice of  $r_b$ depends solely on 
	$\varepsilon,$ $q,$ $N,$ $\gamma_{1+},$ $\gamma_{2+},$ $\mu_{+},$ 
	$\nu_+,$  $\lambda_0$ and $\delta_0.$
\end{prop}

\begin{prop}[{\cite[Theorem 1.5]{Saito2015}}]
\label{prop:6.3}
Assume that $0<\varepsilon<\pi/2$, $1<q<\infty$ and the constants $\mu_{-},$  $\gamma_{1-}>0.$ Let us set
	\begin{equation*}
		\begin{aligned}
			W_q (\Omega_-)&=\big\{ (\Bf_1,f_2,\Bf_3) \in  L_q(\Omega_-)^N \times 
			H^1_q(\Omega_-)\times L_q(\Omega_-)^N : f_2=\dv \Bf_3 \big\},\\
			\CW_q(\Omega_-)&=\big\{ (\Bf_1,g_2,f_2,\Bf_3) : 
			(\Bf_1,f_2,\Bf_3)  \in  W_q (\Omega_-),\, g_2 \in L_q(\Omega_-)\big\}.
		\end{aligned}
	\end{equation*}
	For any $(\bg_-,g_d,\fg_d) \in  W_q(\Omega_-),$ 
	there exist constants $\lambda_0, r_b$ and the operator families
	\begin{equation*}
		\begin{aligned}
			\CA_-^1(\lambda) &\in \Hol\Big( \Sigma_{\ep,\lambda_0} ; 
			\CL\big( \CW_q(\Omega_-) ;H^2_q(\Omega_-)^N \big) \Big),\\
			\CP_-^1(\lambda) &\in \Hol\Big( \Sigma_{\ep,\lambda_0} ; 
			\CL\big( \CW_q(\Omega_-) ;\wh H^1_q(\Omega_-) \big) \Big)
		\end{aligned}
	\end{equation*}
	such that $(\Bv_-^1,\fp_-^1)= \big(\CA_-^1(\lambda),\CP_-^1(\lambda)\big)  \,
	(\bg_-,\lambda^{1/2} g_d, g_d, \lambda \fg_d)$ solves the following system
	\begin{equation}\label{eq:rr_2}
		\left\{ \begin{aligned}
			&\gamma_{1-}\lambda\bv_-^1-\Di \bT_-(\bv_-^1,\fp_-^1)=\bg_-
			&&\quad&\text{in}& \quad  \Omega_-, \\ 
			&\di\bv_-^1= g_d =\dv \fg_d
			&&\quad&\text{in}& \quad  \Omega_-, \\ 
			&\bT_-(\bv_-^1,\fp_-^1)\bn =0
			&&\quad&\text{on}& \quad \Gamma,\\
			&\bv_-^1=0
			&&\quad&\text{on}& \quad S.
		\end{aligned}
		\right.
	\end{equation}
	Moreover, we have 
	\begin{gather*}
		\CR_{\CL\big(\CX_q(\Omega_-); H^{2-j}_q(\Omega_-)^N \big)}
		\Big( \Big\{ (\tau \pa_{\tau})^{\ell}\big( \lambda^{j\slash 2}\CA_-^1(\lambda)\big) : 
		\lambda \in \Sigma_{\ep,\lambda_0} \Big\}\Big) \leq r_b,\\
		\CR_{\CL\big(\CX_q(\Omega_-); L_q(\Omega_-)^N \big)}
		\Big( \Big\{ (\tau \pa_{\tau})^{\ell}   \nabla \CP_-^1(\lambda) : 
		\lambda \in \Sigma_{\ep,\lambda_0} \Big\}\Big) \leq r_b
	\end{gather*}
	for $\ell =0,1,$ $j=0,1,2,$ and $\tau = \Im \lambda.$ 
	Above the choice of  $r_b$ depends solely on 
	$\varepsilon,$ $q,$ $N,$ $\gamma_{1-},$ $\mu_{-},$ $\lambda_0$ and $b.$
\end{prop}

Inspired by Propositions \ref{prop:v+1} and \ref{prop:6.3},  we decompose the solution  of \eqref{eq:rr} as 
\begin{equation*}
	\bv_+= \bv_+^1+\bv_+^2,\quad 
	\bv_-=\bv_-^1+\bv_-^2, \quad \fp_-=\fp_-^1 +\fp_-^2
\end{equation*}
with 
$$\Bv_+^1= \CA_+^1(\lambda) \, \bg_+ 
\quad \text{and} \quad 
(\Bv_-^1,\fp_-^1)= \big(\CA_-^1(\lambda),\CP_-^1(\lambda)\big)  
\,(\bg_-,\lambda^{1/2} g_d, g_d, \lambda \fg_d).$$ 
Thus the compensated part $(\bv_+^2,\bv_-^2, \fp_-^2)$ satisfies 
\begin{equation}\label{eq:rr_3}
	\left\{ \begin{aligned}
		&\gamma_{1+}\lambda\bv_+^2-\Di\bS_{\delta+}(\bv_+^2)=0
		&&\quad&\text{in}& \quad  \Omega_+, \\
		&\gamma_{1-}\lambda\bv_-^2-\Di \bT_-(\bv_-^2,\fp_-^2)=0
		&&\quad&\text{in}& \quad  \Omega_-, \\ 
           &\di\bv_-^2=0
		&&\quad&\text{in}& \quad  \Omega_-, \\ 
		&\bS_{\delta+}(\bv_+^2)\bn
		-\bT_-(\bv_-^2,\fp_-^2)\bn=\bh
		&&\quad&\text{on}& \quad \Gamma,\\
		& \bv_+^2-\bv_-^2=\bk-(\bv_+^1-\bv_-^1)
		&&\quad&\text{on}& \quad \Gamma,	\\
		&\bv_-^2= 0
		&&\quad&\text{on}& \quad S.   
	\end{aligned}
	\right.
\end{equation}

Next, for $\bG_{\lambda}=(\bg_+,\bg_-,\lambda^{1/2} g_d, g_d, \lambda \fg_d, 
\lambda^{1/2} \bh,\bh,\lambda\bk,\lambda^{1/2}\bk,\bk)\in \CY_q(\dot \Omega),$ we denote 
\begin{equation*}
	\begin{aligned}
		\bH_{\lambda} &=(\bg_+,\bg_-,\lambda^{1/2} g_d, g_d, \lambda \fg_d, \lambda\bk,\lambda^{1/2}\bk,\bk),\\
		\CK(\lambda) \bH_{\lambda} &= \bk-\CA_+^1(\lambda) \bg_+
		+\CA_-^1(\lambda)\,(\bg_-,\lambda^{1/2} g_d, g_d, \lambda \fg_d),\\
		\CT(\lambda) \bG_{\lambda}&= 
		(\lambda^{1/2} \bh,\bh,\lambda\CK(\lambda) \bH_{\lambda},\lambda^{1/2}\CK(\lambda) \bH_{\lambda},\CK(\lambda) \bH_{\lambda}).
	\end{aligned}
\end{equation*}
Then the operator family $\CT(\lambda) $ can be regarded as the $\CR$-bounded mappings from $\CY_q(\dot \Omega)$ into $\CZ_q(\dot \Omega)$  thanks to Propositions \ref{prop:v+1} and \ref{prop:6.3}.
Thus 
\begin{equation*}
	(\bv_{\pm}^2,\fp_-^2)=\big(\CS_{\pm}(\lambda),\CP_-(\lambda)\big) \circ \CT(\lambda) \bG_{\lambda}
\end{equation*}
solves \eqref{eq:rr_3} clearly. Moreover, Remark \ref{rmk:R-bdd} and Theorem \ref{thm:5.1} yield 
\begin{equation}\label{es:vp-2}
    \begin{aligned}
       	\CR_{\CL\big(\CY_q(\dot \Omega); H^{2-j}_q(\Omega_\pm)^N \big)}
	\Big( \Big\{ (\tau \pa_{\tau})^{\ell}\big( \lambda^{j/2}
	\CS_{\pm}(\lambda) \circ \CT(\lambda) \big) : 
	\lambda \in \Gamma_{\ep,\lambda_0} \Big\}\Big) \lesssim 1,\\
	\CR_{\CL\big(\CY_q(\dot \Omega); L_q(\Omega_-)^N \big)}
	\Big( \Big\{ (\tau \pa_{\tau})^{\ell}   \nabla \CP_-(\lambda)\circ \CT(\lambda) : 
	\lambda \in \Gamma_{\ep,\lambda_0} \Big\}\Big) \lesssim 1. 
    \end{aligned}
\end{equation}

Finally, we define 
\begin{align*}
	\CA_{+}^0 (\lambda) \bG_{\lambda}&=\CA_+^1(\lambda)\bg_+
	+\CS_{+}(\lambda) \circ \CT(\lambda) \bG_{\lambda},\\
	\CA_{-}^0(\lambda)\bG_{\lambda}
	&=\CA_-^1(\lambda) (\bg_-,\lambda^{1/2} g_d, g_d, \lambda \fg_d)
	+\CS_{-}(\lambda) \circ \CT(\lambda) \bG_{\lambda},\\
	\CB_-^0(\lambda)\bG_{\lambda}&=
	\CP_-^1(\lambda) (\bg_-,\lambda^{1/2} g_d, g_d, \lambda \fg_d)
	+\CP_-(\lambda)\circ \CT(\lambda)\bG_{\lambda}.
\end{align*}
According to Remark \ref{rmk:R-bdd}, \eqref{es:vp-2}, Propositions \ref{prop:v+1} and \ref{prop:6.3}, we have the desired estimates for the operator families $\CA_\pm^0(\lambda)$ and  $\CB_-^0(\lambda).$

 \section{Linear theory on the maximal $L_p$-$L_q$ regularity property}
\label{sec:mr}

In this section, we shall prove Theorem \ref{thm:mr}.
As we will see later, the main difficulty is the case where source terms $f_+=g_{\pm}=g_d=\fg_d=\bh=0$ in \eqref{eql:1}. This essentially leads to the semigroup theory.
However, the model \eqref{eql:1} contains the incompressible part in the region $\Omega_-$. 
Thus we need to introduce some reduced two-phase gas-liquid operator to overcome the difficulty from the pressure term $\fp_-$ as the one-phase incompressible viscous flow (see for instance \cite{Shibata2020} and the references therein). 
Subsections \ref{subsec:reduced_1} and \ref{subsec:reduced_2} are dedicated to the reduced problem of \eqref{eql:1}. In Subsection \ref{subsec:reduced_3}, we will derive the maximal $L_p$-$L_q$ regularity property in view of Theorem \ref{thm:main2} and the Weis' theory \cite{Weis2001}.

\subsection{Resolvent problem for the reduced two-phase operator} 
\label{subsec:reduced_1}
Let us review some result from the one-phase incompressible problem in the layer $\Omega_-.$ 
According to \cite[p.1493]{Saito2018}, we have the following technical lemma.
\begin{lem}\label{lem:weak}
Let $1<q<\infty$ and $q'=q/(q-1).$
For any $\Bf \in L_q(\Omega_-)^N$ and $g\in W^{1-1/q}_q(\Gamma)$
there exists a solution $\fq \in H^1_{q}(\Omega_-)+H^1_{q,\Gamma}(\Omega_-)$ satisfying 
\begin{equation*}
\left\{ \begin{aligned}
    & (\nabla \fq, \nabla \varphi)_{\Omega_-}=(\Bf,\nabla \varphi)_{\Omega_-}, \,\,\,\forall \,\,\varphi \in H^1_{q',\Gamma}(\Omega_-),\\
    & \fq|_{\Gamma}=g.\\
    \end{aligned} \right.
\end{equation*}
Moreover, there exists a contant $C$ such that 
\begin{equation*}
    \|\fq\|_{H^1_q(\Omega_-)} \leq C \big( \|\Bf\|_{L_q(\Omega_-)} + \|g\|_{W^{1-1/q}_q(\Gamma)} \big).
\end{equation*}
\end{lem}

Let $\rho_+ \in H^1_q(\Omega_+)$ and $\bv_{\pm}\in H^2_q(\Omega_{\pm})^N.$ According to Lemma \ref{lem:weak},  
there exists 
$$\theta=\CK_-(\rho_+,\bv_+,\bv_-) \in H^1_{q}(\Omega_-)+H^1_{q,\Gamma}(\Omega_-)$$ 
solving the following variational problem
\begin{equation}\label{def:K}
\left\{ \begin{aligned}
  &  (\nabla \theta, \nabla \varphi)_{\Omega_-} 
    = \big(\DV  \bS_- (\bv_-) -\nabla \dv \bv_-, \nabla \varphi \big)_{\Omega_-},\,\,\,\forall \,\,\varphi \in H^1_{q',\Gamma}(\Omega_-),\\
  &      \theta|_{\Gamma} =  \bS_-(\bv_-)\bn \cdot \bn -\dv \bv_- 
    -\bS_{+}(\bv_+) \bn \cdot \bn +\gamma_{2+} \rho_+.
        \end{aligned} \right.
\end{equation}
Moreover, we have 
\begin{equation*}
    \|\CK_-(\rho_+,\bv_+,\bv_-)\|_{H^1_{q}(\Omega_-)} \lesssim \|\bv_+\|_{H^2_q(\Omega_+)}+\|\bv_-\|_{H^2_q(\Omega_-)}
    +\|\rho_+\|_{H^{1}_q(\Omega_+)}.
\end{equation*}

With the operator $\CK_-$ defined by \eqref{def:K} at hand, we consider the following generalized reduced resolvent problem:
\begin{equation}\label{gres:2}
	\left\{ \begin{aligned}
		&\lambda\rho_++\gamma_{1+}\di \bv_+=f_+
		&&\quad&\text{in}& \quad \Omega_+, \\
		&\lambda\bv_+-\gamma_{1+}^{-1}\DV\big(\bS_+(\bv_+)-\gamma_{2+}\rho_+\bI\big)=\bg_+
		&&\quad&\text{in}& \quad  \Omega_+, \\
		&\lambda\bv_--\gamma_{1-}^{-1} \Di\big(\bS_-(\bv_{-})-\CK_-(\rho_+,\bv_+,\bv_-)\bI\big)=\bg_-
		&&\quad&\text{in}& \quad  \Omega_-, \\ 
		&\big(\bS_+(\bv_+)-\gamma_{2+}\rho_+\bI\big)\bn
		=\big(\bS_-(\bv_{-})-\CK_-(\rho_+, \bv_+,\bv_-)\bI\big)\bn
		&&\quad&\text{on}& \quad \Gamma,\\
		& \bv_+-\bv_-=\bk
		&&\quad&\text{on}& \quad \Gamma,	\\
		&\bv_-= 0
		&&\quad&\text{on}& \quad S.   
	\end{aligned}
	\right.
\end{equation}
For \eqref{gres:2}, we have 
\begin{thm}\label{thm:gres-2}
	Assume that $0<\varepsilon<\pi/2$, $1<q<\infty$ and the constants $\mu_{\pm},$ $\nu_+,$ $\gamma_{1\pm},$ $\gamma_{2+}>0.$
	Let us introduce the functional spaces for source terms in \eqref{gres:2} as follows 
	\begin{equation*}
		\begin{aligned}
			X_q^1 (\dot \Omega)&=\big\{ (f_1,\Bf_2,\Bf_3,\Bf_4) :
			f_1\in H^1_q(\Omega_+),
			\Bf_2\in L_q(\Omega_+)^N, \Bf_3\in J_q(\Omega_-)^N,
		   \Bf_4\in H^{2}_q(\dot \Omega)^N\big\}, \\
			\CX_q^1(\dot\Omega)&=\big\{ (f_1,\Bf_2,\Bf_3,\bh_4,\bg_4,\Bf_4) :
			(f_1,\Bf_2,\Bf_3,\Bf_4) \in X^1_q(\dot \Omega),\,
			\bh_4\in L_q(\dot\Omega)^N,\bg_4\in H^1_q(\dot\Omega)^N\big\}.
		\end{aligned}
	\end{equation*} 
	In addition, we write for short that 
	\begin{equation*}
			\bF^1=(f_+, \bg_+,\bg_-,\bk),\quad
			\bF^1_{\lambda}=(f_+,\bg_+,\bg_-,\lambda\bk, \lambda^{1/2} \bk,\bk).
	\end{equation*}
	For any $\bF^1\in X_q^1(\dot\Omega),$ there exist constants $\lambda_0, r_b$ and the operator families
	\begin{equation*}
		\begin{aligned}
			\Theta^1_+(\lambda) &\in \Hol\Big( \Sigma_{\ep, \lambda_0}  \cap K_{\ep} ; 
			\CL\big( \CX^1_q(\dot \Omega) ; H^1_q(\Omega_+) \big) \Big),\\
			\CA^1_\pm(\lambda) &\in \Hol\Big( \Sigma_{\ep, \lambda_0}  \cap K_{\ep} ; 
			\CL\big( \CX_q^1(\dot \Omega) ;H^2_q(\Omega_\pm)^N \big) \Big)
		\end{aligned}
	\end{equation*}
	such that $(\rho_+, \bv_+,\bv_-)
	= \big(\Theta_+^1(\lambda) , \CA^1_+(\lambda), \CA^1_-(\lambda)\big) \,\bF^1_{\lambda}$
	solves \eqref{gres:2}.  Here the regions $\Sigma_{\ep, \lambda_0}$ and $K_{\ep}$ are defined by \eqref{eq:KSigma}. Moreover, we have 
	\begin{gather*}
		\CR_{\CL\big(\CX^1_q(\dot \Omega); H^{1}_q(\Omega_+) \big)}
		\Big( \Big\{ (\tau \pa_{\tau})^{\ell}\big( \lambda^{k}
		\Theta^1_+(\lambda)\big) : 
		\lambda \in \Sigma_{\ep, \lambda_0}  \cap K_{\ep} \Big\}\Big) \leq r_b,\\
		\CR_{\CL\big(\CX^1_q(\dot\Omega); H^{2-j}_q(\Omega_\pm)^N \big)}
		\Big( \Big\{ (\tau \pa_{\tau})^{\ell}\big( \lambda^{j\slash 2}
		\CA^1_\pm(\lambda)\big) : 
		\lambda \in \Sigma_{\ep, \lambda_0}  \cap K_{\ep} \Big\}\Big) \leq r_b
	\end{gather*}
	for $k,\ell =0,1,$ $j=0,1,2,$ and $\tau = \Im \lambda.$ 
	Above the choice of  $r_b$ depends solely on 
	$\varepsilon,$ $q,$ $N,$ $\gamma_{1\pm},$ $\gamma_{2+},$ $\mu_{\pm},$ 
	$\nu_+,$  $\lambda_0$ and $b.$
\end{thm}

In fact, it is not hard to prove Theorem \ref{thm:gres-2} by Theorem \ref{thm:main2} and the following proposition.
\begin{prop}
Let $1<q<\infty$ and $\lambda \in \Sigma_{\ep,\lambda_0} \cap K_{\ep}.$
Let us introduce
	\begin{equation}\label{def:Jq}
		J_q(\Omega_-)=\big\{\bu_-\in L_q(\Omega_-)^N\,|\,(\bu_-,\nabla\varphi)_{\Omega_-}=0,\quad\forall\,\,\varphi\in H^1_{q',\Gamma}(\Omega_-)\big\}.
	\end{equation}
Suppose that $f_+\in H^1_q(\Omega_+),$ $\bg_{+} \in L_q(\Omega_{+})^N,$ $\bg_{-} \in J_q(\Omega_{-})^N,$ and $\bk\in H^{2}_q(\dot \Omega)^N.$ For simplicity, let us denote 
\begin{equation*}
\begin{aligned}
      E^{1}_q (\dot \Omega)&=H^1_{q}(\Omega_+) \times H^2_{q}(\Omega_+)^N
    \times \big( H^2_{q}(\Omega_-)^N \cap J_q(\Omega_-) \big)\times H^1_q(\Omega_-), \\
          E^{2}_q (\dot \Omega)&=H^1_{q}(\Omega_+) \times H^2_{q}(\Omega_+)^N \times H^2_{q}(\Omega_-)^N.
\end{aligned}
\end{equation*}
Then the following assertions hold true.
\begin{enumerate}[label=(\arabic*)]	
    \item If $(\rho_+,\bv_+,\bv_-,\fp_-) \in E^1_q(\dot \Omega) $ solves the problem
    \begin{equation}\label{gres:1}
	\left\{ \begin{aligned}
		&\lambda\rho_++\gamma_{1+}\di \bv_+=f_+
		&&\quad&\text{in}& \quad \Omega_+, \\
		&\lambda\bv_+-\gamma_{1+}^{-1}\DV\big(\bS_+(\bv_+)-\gamma_{2+}\rho_+\bI\big)=\bg_+
		&&\quad&\text{in}& \quad  \Omega_+, \\
		&\lambda\bv_--\gamma_{1-}^{-1} \Di\big(\bS_-(\bv_{-})-\fp_-\bI\big)=\bg_-
		&&\quad&\text{in}& \quad  \Omega_-, \\ 
  		& \di\bv_-=0
		&&\quad&\text{in}& \quad  \Omega_-, \\ 
		&\big(\bS_+(\bv_+)-\gamma_{2+}\rho_+\bI\big)\bn
		=\big(\bS_-(\bv_{-})-\fp_-\bI\big)\bn
		&&\quad&\text{on}& \quad \Gamma,\\
		& \bv_+-\bv_-=\bk
		&&\quad&\text{on}& \quad \Gamma,	\\
		&\bv_-= 0
		&&\quad&\text{on}& \quad S,  
	\end{aligned}
	\right.
\end{equation} then $(\rho_+,\bv_+,\bv_-)$ satisfies \eqref{gres:2}.
    \item If $(\rho_+,\bv_+,\bv_-) \in E^2_q(\dot \Omega)$ solves \eqref{gres:2}, then $\big(\rho_+,\bv_+,\bv_-, \CK_-(\rho_+, \bv_+,\bv_-)\big)$ is a solution of \eqref{gres:1}.
\end{enumerate}

\end{prop}
\begin{proof}
Suppose that $(\rho_+,\bv_+,\bv_-,\fp_-) \in E^1_q(\dot \Omega)$  satisfies \eqref{gres:1}.
Let us denote $$\theta=\CK_-(\rho_+,\bv_+,\bv_-) \in H^1_q(\Omega_-).$$ 
For any $\varphi \in H^1_{q',\Gamma}(\Omega_-),$ taking the inner product of $\eqref{gres:1}_3$ with $\nabla \varphi$ yields 
\begin{equation*}
\begin{aligned}
 0=(\bg_-,\nabla \varphi)_{\Omega_-} &= (\lambda\bv_-,\nabla \varphi)_{\Omega_-}
   -\big(\gamma_{1-}^{-1}\Di \bS_-(\bv_-),\nabla \varphi\big)_{\Omega_-}  
   +(\gamma_{1-}^{-1}\nabla \fp_-,\nabla \varphi)_{\Omega_-} \\
   &= \gamma_{1-}^{-1}\big(\nabla (\fp_--\theta),\nabla \varphi\big)_{\Omega_-},
\end{aligned}
\end{equation*}
 where we have used the facts 
\begin{equation*}
  \bg_- \in J_q(\Omega_-), \quad  \dv \bv_-=0\quad 
   \hbox{and}\quad  \bv_-|_{S}=0.
\end{equation*}
Moreover, we have 
\begin{equation*}
   (\fp_--\theta)|_{\Gamma}=0.
\end{equation*}
Then the uniqueness of the variational problem (see \cite[Proposition 3.4]{Saito2018} for instance) yields $\fp_-=\theta.$
\medskip

On the other hand, suppose that $(\rho_+,\bv_+,\bv_-) \in E^2_q(\dot \Omega)$ solves \eqref{gres:2}.
It suffices to verify that  $\dv \bv_-$ vanishes in $\Omega_-.$
The definition of the operator $\theta=\CK_-(\rho_+, \bv_+,\bv_-)$ implies that $\dv \bv_-|_{\Gamma}=0.$
Then taking the inner product of $\eqref{gres:2}_3$ with $\nabla \varphi$ and integration by parts yield
\begin{equation*}
\begin{aligned}
 0=-(\bg_-,\nabla \varphi)_{\Omega_-} &= -(\lambda\bv_-,\nabla \varphi)_{\Omega_-}
   +\big(\gamma_{1-}^{-1}\Di \bS_-(\bv_-),\nabla \varphi\big)_{\Omega_-}  
   -(\gamma_{1-}^{-1}\nabla \theta,\nabla \varphi)_{\Omega_-} \\
   &= \lambda (\dv \bv_-,\varphi) + \gamma_{1-}^{-1}\big(\nabla \dv \bv_-,\nabla \varphi\big)_{\Omega_-}.
\end{aligned}
\end{equation*}
Then \cite[Proposition 4.2]{Saito2018} furnishes that $\bv_-$ is divergence-free in $\Omega_-.$
\end{proof}

\subsection{Semigroup generated by the reduced two-phase operator}
\label{subsec:reduced_2}
In this subsection, we consider the following evolution problem:
\begin{equation}\label{ibeq:0}
	\left\{ \begin{aligned}
		&\partial_t \rho_++\gamma_{1+}\di \bv_+=0
		&&\quad&\text{in}& \quad \Omega_+\times \BR_+, \\
		&\pd_t \bv_+-\gamma_{1+}^{-1}\DV\big(\bS_+(\bv_+)-\gamma_{2+}\rho_+\bI\big)=0
		&&\quad&\text{in}& \quad  \Omega_+\times \BR_+, \\
		&\pd_t\bv_--\gamma_{1-}^{-1} \Di\big(\bS_-(\bv_{-})-\CK_-(\rho_+,\bv_+,\bv_-)\bI\big)=0
		&&\quad&\text{in}& \quad  \Omega_-\times \BR_+, \\ 
		&\big(\bS_+(\bv_+)-\gamma_{2+}\rho_+\bI\big)\bn
		=\big(\bS_-(\bv_{-})-\CK_-(\rho_+, \bv_+,\bv_-)\bI\big)\bn
		&&\quad&\text{on}& \quad \Gamma\times \BR_+,\\
		& \bv_+=\bv_-
		&&\quad&\text{on}& \quad \Gamma\times \BR_+,	\\
		&\bv_-= 0
		&&\quad&\text{on}& \quad S\times \BR_+,\\   
        		&(\rho_+,\bv_+)\big|_{t=0}=(\rho_{+,0},\bv_{+,0})
		&&\quad&\text{in}& \quad \Omega_+,\\
		&\bv_-\big|_{t=0}=\bv_{-,0}
		&&\quad&\text{in}& \quad \Omega_-
	\end{aligned}
	\right.
\end{equation}
for $\BR_+=(0,\infty)$ and $\bn=(0,\dots,0,1)^{\top}.$ In fact, the system \eqref{ibeq:0} can be written as some abstract Cauchy problem associated to some reduced two-phase operator.
Let us introduce some notations.
\begin{dfn}\label{def:tp-A}
Let $1<q<\infty$ and $\dot\Omega=\Omega_+ \cup \Omega_-.$ Let us set the underlying space 
\footnote{See \eqref{def:Jq} for the definition of $J_q(\Omega_-)$.}
\begin{equation}\label{def:under}
    \CX_q(\dot\Omega)=H^1_q(\Omega_+) \times L_q(\Omega_+)^N \times J_q(\Omega_-).
\end{equation}
Then, according to \eqref{ibeq:0}, we define the reduced two-phase gas-liquid operator 
\begin{equation*}
    \CA \begin{pmatrix}
        \theta_+\\\bu_+\\\bu_-
    \end{pmatrix}=\begin{pmatrix}
        -\gamma_{1+} \dv \bu_+ \\
        \gamma_{1+}^{-1}\DV\big(\bS_+(\bu_+)-\gamma_{2+}\theta_+\bI\big)\\
        \gamma_{1-}^{-1} \Di\big(\bS_-(\bu_{-})-\CK_-(\theta_+,\bu_+,\bu_-)\bI\big)
    \end{pmatrix}
\end{equation*}
whenever $(\theta_+,\bu_+,\bu_-)$ lies in the space
    \begin{equation}\label{cdt:intial0}
    \begin{aligned}
\CD_q(\CA)=\Big\{(\theta_+,\bu_+,\bu_-)\in \CX_q(\dot \Omega):\, &
        \bu_\pm \in H^2_q(\Omega_\pm)^N, \, \bu_+|_{\Gamma}=\bu_-|_{\Gamma},\, \bu_-|_{S}=0,\\
       & \quad \text{and}\,\,\, \CT_{\bn} \big(\mu_+\bD(\bu_{+})\bn\big)|_{\Gamma}
         =\CT_{\bn} \big(\mu_-\bD(\bu_{-})\bn \big) |_{\Gamma} \Big\}.
    \end{aligned}
    \end{equation}
Here the projection $\CT_{\bn}$ is defined in \eqref{def:proj}.
\end{dfn}

Now, writing $V=(\rho_{+},\bv_{+},\bv_{-})$ and $V_0=(\rho_{+,0},\bv_{+,0},\bv_{-,0}),$ we may regard \eqref{ibeq:0} as the following abstract Cauchy problem 
\begin{equation*}
    \pd_t \bV -\CA \bV =0, \quad \bV|_{t=0}=\bV_0.
\end{equation*}
According to Theorem \ref{thm:gres-2}, there exists a $\lambda_0>0$ such that 
$\Sigma_{\ep,\lambda_0}\cap K_\varepsilon$ belongs to the resolvent set $\rho(\CA)$ of $\CA.$
Thus the standard semigroup theory allows us to establish the following result:
\begin{thm}\label{thm:sg}
Let $1<q<\infty.$ Then the reduced two-phase gas-liquid operator $\CA$ in Definition \ref{def:tp-A}  generates a $C^0$-semigroup $\{T(t)\}_{t\geq 0}$ on $\CX_q(\dot \Omega)$ in \eqref{def:under},  which is also analytic.
\end{thm}
Moreover, Theorem \ref{thm:sg} and the interpolation argument (see \cite[Theorem 3.9]{ShiShi2008} for instance) imply that the solution $(\rho_+,\bv_{+}, \bv_{-}, \fp_-)$ of \eqref{ibeq:0} satisfies
\begin{equation}\label{es:sg}
\CN_{p,q,\gamma_0}(\rho_+,\bv_{+}, \bv_{-}, \fp_-) \leq 
C \|(\rho_{+,0}, \bv_{+,0},\bv_{-,0})\|_{\CD_{p,q}(\dOm)} 
\end{equation}
with $\fp_-=\CK_-(\rho_+,\bv_+,\bv_-),$ where the notations $\CD_{p,q}(\dOm)$ and $\CN_{p,q,\gamma_0}$ have been defined in Definition \ref{def:FD} and Theorem \ref{thm:mr} respectively.

\subsection{Proof of Theorem \ref{thm:mr}}
\label{subsec:reduced_3}
In this subsection, we begin with the following problem 
\begin{equation}\label{eq:MR_1}
	\left\{ \begin{aligned}
		&\pd_t\rho_+^1 + \gamma_{1+}\, \dv \bv_+^1=f_+
		&&\quad&\text{in} &\quad \Omega_+ \times \BR, \\
		&\gamma_{1+}\pd_t\bv_+^1 - \DV\bT_+(\bv_+^1,\gamma_{2+}\rho_+^1)  = \bg_+
		&&\quad&\text{in}& \quad \Omega_+ \times \BR, \\
		&\gamma_{1-}\pd_t\bv_-^1  - \Di\bT_-(\bv_{-}^1,\fp_-^1) = \bg_-
		&&\quad&\text{in}& \quad \Omega_- \times \BR, \\
		& \di \bv_-^1 = g_d=\di \fg_d
		&&\quad&\text{in}& \quad \Omega_- \times \BR, \\
		&\big(\bT_+(\bv_+^1,\gamma_{2+}\rho_+^1)- \bT_-(\bv_-^1,\fp_-^1)\big)\bn
		= \bh
		&&\quad&\text{on}& \quad \Gamma \times \BR,\\
		&\bv_{+}^1=\bv_{-}^1
		&&\quad&\text{on}& \quad \Gamma \times \BR,  \\
		& \bv_-^1=0
		&&\quad&\text{on}& \quad   S \times \BR.
	\end{aligned}
	\right.
\end{equation}
Then applying the Laplace transform \eqref{eq:Laplace}  to \eqref{eq:MR_1} implies that 
\begin{equation*}
	\left\{ \begin{aligned}
		&\lambda \wh{\rho_+^1} + \gamma_{1+}\, \dv \wh{\bv_+^1}=\wh{f_+}
		&&\quad&\text{in} &\quad \Omega_+, \\
		&\gamma_{1+}\lambda \wh{\bv_+^1} - \DV\bT_+(\wh{\bv_+^1},\gamma_{2+}\wh{\rho_+^1})  = \wh{\bg_+}
		&&\quad&\text{in}& \quad \Omega_+, \\
		&\gamma_{1-}\lambda \wh{\bv_-^1}  - \Di\bT_-(\wh{\bv_{-}^1},\wh{\fp_-^1}) = \wh{\bg_-}
		&&\quad&\text{in}& \quad \Omega_- , \\
		& \di \wh{\bv_-^1} = \wh{g_d}=\di \wh{\fg_d}
		&&\quad&\text{in}& \quad \Omega_- , \\
		&\big(\bT_+(\wh{\bv_+^1},\gamma_{2+}\wh{\rho_+^1})- \bT_-(\wh{\bv_-^1},\wh{\fp_-^1})\big)\bn
		= \wh{\bh}
		&&\quad&\text{on}& \quad \Gamma,\\
		&\wh{\bv_{+}^1}=\wh{\bv_{-}^1}
		&&\quad&\text{on}& \quad \Gamma ,  \\
		&\wh{\bv_-^1}=0
		&&\quad&\text{on}& \quad   S.
	\end{aligned}
	\right.
\end{equation*}
Here $\widehat{f}=\CF_{L} [f]$ denotes the Laplace transform of $f$. In view of \eqref{eq:half} and Theorem \ref{thm:main2}, the system \eqref{eq:MR_1} admits a solution given by
\begin{equation*}
\begin{aligned}
(\rho_+^1, \bv_+^1,\bv_-^1,\fp_-^1)
	&=e^{\gamma t} \CF\Big[ \big(\Theta_+(\lambda) , \CA_+(\lambda), \CA_-(\lambda),\CB_-(\lambda)\big) \,\bG_{\lambda}\Big]  \\
    &=e^{\gamma t} \CF\Big[ \big(\Theta_+(\lambda) , \CA_+(\lambda), \CA_-(\lambda),\CB_-(\lambda)\big)
    \,\CF[e^{-\gamma t}\bF_{\gamma}]\Big] 
\end{aligned}
\end{equation*}
for 
\begin{equation*}
\begin{aligned}
\bG_{\lambda}&=(\wh{f_+},\wh{\bg_+},\wh{\bg_-},\lambda^{1/2} \wh{g_d}, \wh{g_d}, \lambda \wh{\fg_d}, 
			\lambda^{1/2} \wh{\bh},\wh{\bh}),\\
 \bF_{\gamma}&=(f_+,\bg_+,\bg_-,\Lambda^{1/2}_{\gamma} g_d, g_d,  \partial_t \fg_d, 
			\Lambda^{1/2}_{\gamma} \bh,\bh).
\end{aligned}
\end{equation*}
Thus keeping $\CF_{p,q,\gamma_0}(\dot \Omega)$ in Definition \ref{def:FD} in mind,
Weis' theory \cite{Weis2001} yields 
 \begin{equation}\label{es:v1}
\CN_{p,q,\gamma_0}(\rho_+^1,\bv_{+}^1, \bv_{-}^1, \fp_-^1) 
\lesssim \|(f_+, \bg_+, \bg_-, g_d,\fg_d,\bh) \|_{\CF_{p,q,\gamma_0}(\dot \Omega)}
\end{equation}
for some $\gamma_0>0$.
Moreover, by the trace method and \eqref{es:v1}, we have 
\begin{equation}\label{es:v10}
    \begin{aligned}
        \|\rho_+^1|_{t=0}\|_{H^1_q(\Omega_+)} &\lesssim \|e^{-\gamma t}\rho_+^1\|_{H^1_p(\BR_+;H^1_q(\Omega_+))} 
        \lesssim \|(f_+, \bg_+, \bg_-, g_d,\fg_d,\bh) \|_{\CF_{p,q,\gamma_0}(\dot \Omega)},\\
        \|\bv_{\pm}^1|_{t=0}\|_{B^{2-2/p}_{q,p}(\Omega_{\pm})}& \lesssim 
\|e^{-\gamma t}\partial_t \bv_{\pm}^1\|_{L_p(\BR_+;L_q(\Omega_{\pm}))} 
+\|e^{-\gamma t}\bv_{\pm}^1\|_{L_p(\BR_+;H^2_q(\Omega_{\pm}))} \\
& \lesssim \|(f_+, \bg_+, \bg_-, g_d,\fg_d,\bh) \|_{\CF_{p,q,\gamma_0}(\dot \Omega)}.
    \end{aligned}
\end{equation}

Now, suppose that the solution of \eqref{eql:1} can be decomposed as follows
\begin{equation*}
  (\rho_+, \bv_+,\bv_-,\fp_-)
  = (\rho_+^1, \bv_+^1,\bv_-^1,\fp_-^1)|_{t\in \BR_+}+(\rho_+^2, \bv_+^2,\bv_-^2,\fp_-^2).
\end{equation*}
with 
\begin{equation*}
(\rho_{+,0}^2, \bv_{+,0}^2,\bv_{-,0}^2)=
(\rho_{+,0}, \bv_{+,0},\bv_{-,0})-(\rho_{+}^1, \bv_{+}^1,\bv_{-}^1)|_{t=0}.    
\end{equation*}
According to \eqref{eq:1.4}, \eqref{eq:MR_1} and \eqref{es:v10},  it is not hard to find that  $(\rho_{+,0}^2, \bv_{+,0}^2,\bv_{-,0}^2)\in \CD_q(\CA)$ (see \eqref{cdt:intial0}).
Thus $(\rho_+^2, \bv_+^2,\bv_-^2,\fp_-^2)= T(t) (\rho_{+,0}^2, \bv_{+,0}^2,\bv_{-,0}^2)$ satisfies 
\begin{equation}\label{es:rhovp_2}
\CN_{p,q,\gamma_0}(\rho_+^2,\bv_{+}^2, \bv_{-}^2, \fp_-^2) \lesssim  
\|(\rho_{+,0}, \bv_{+,0},\bv_{-,0})\|_{\CD_{p,q}(\dOm)} 
+\|(f_+, \bg_+, \bg_-, g_d,\fg_d,\bh) \|_{\CF_{p,q,\gamma_0}(\dot \Omega)}
\end{equation}
for $\fp_-^2=\CK_-(\rho_+^2,\bv_+^2,\bv_-^2)$ thanks to \eqref{es:sg} and \eqref{es:v10}. 
Then combining the estimates \eqref{es:v1} and \eqref{es:rhovp_2} implies the desired bound of $  (\rho_+, \bv_+,\bv_-,\fp_-)$.
\smallbreak 

At last, the uniqueness part of Theorem \ref{thm:mr} also holds true by applying \eqref{es:sg} to some dual problem (see \cite[Subsection 4.5]{EvBS2014} for instance). This completes the proof of Theorem \ref{thm:mr}.

\newpage
\appendix

\section{Nonlinear gas-liquid two-phase problem} 
\label{app:non}
As we mentioned in the introduction, the system \eqref{eql:1} comes from some liquid-gas two-phase flow problem in some unbounded domain
$\Omega_{t,-} \cup \Omega_{t,-} \subset \mathbb{R}^N$ for $N\geq 2$ (see Figure \ref{Fig:Omega2}). More precisely, the motion of the viscous gases in $\Omega_{t,+}$ and the viscous liquid in $\Omega_{t,-}$ is separated by a moving sharp interface $\Gamma_t$ with the fixed bottom $S$. 
Furthermore, the viscous fluids obey the Compressible-Incompressible Navier-Stokes equations in view of the \emph{Eulerian} coordinates as follows:
\begin{equation}\label{eq:NS}
	\left\{ \begin{aligned}
		& \pd_t\rho_++\di(\rho_+\bu_+)=0
		&&\quad&\text{in}& \quad  \bigcup_{0<t<T}\Omega_{t,+}\times \{t\}, \\
		&\rho_+(\pd_t\bu_++\bu_+\cdot\nabla\bu_+)
		-\Di \bT_+\big(\bu_+,P(\rho_+)\big)=0
		&&\quad&\text{in}& \quad \bigcup_{0<t<T} \Omega_{t,+}\times \{t\}, \\
		&\rho_-(\pd_t\bu_- +\bu_- \cdot\nabla\bu_-)-\Di \bT_-(\bu_-,\pi_-)=0,
		\quad \di\bu_-=0
		&&\quad&\text{in}& \quad \bigcup_{0<t<T} \Omega_{t,-}\times \{t\}, \\
		&\bT_+\big(\bu_+,P(\rho_+)\big)\bn_t- \bT_-(\bu_-,\pi_-)\bn_t=0
		&&\quad&\text{on}& \quad \bigcup_{0<t<T} \Gamma_t\times \{t\}, \\
		& \bu_+-\bu_-=0,\quad V_{\Gamma_t}=\bu_-\cdot\bn_t
		&&\quad&\text{on} & \quad \bigcup_{0<t<T} \Gamma_t\times \{t\} , \\
		& \bu_-=0
		&&\quad&\text{on}& \quad   S \times (0,T), \\
		&(\rho_+,\bu_+)\big|_{t=0}=(\rho_{+,0},\bu_{+,0})
		&&\quad&\text{in}& \quad \Omega_+,\\
		&\bu_-\big|_{t=0}=\bu_{-,0}
		&&\quad&\text{in}& \quad \Omega_-.
	\end{aligned}
	\right.
\end{equation} 
In \eqref{eq:NS}, the given positive constant $\rho_-$ denotes the mass density of the liquid,
and the unknowns $\bu_{\pm}$, $\rho_+$ and $\pi_-$ are velocity fields of the fluids, the mass density of the gas and the scalar pressure. Besides, $P(\cdot)$ is a smooth function on $[0,\infty)$, $\bn_t$ is the unit normal from $\Omega_{t,-}$ to $\Omega_{t,+}$, $V_{\Gamma_t}$ denotes the normal velocity of the moving surface $\Gamma_t,$
and $\rho_{+,0},\bu_{\pm,0}$ are given initial data.
For simplicity, we assume that the initial state of the unbounded domain $\Omega_{t,-} \cup \Omega_{t,-}$ coincides with the reference domain $\Omega_+ \cup \Omega_-$ in \eqref{dfn:Omega}.

\subsection{Lagrangian transformation}
\label{ss:Lag}

For the application of the mathematical analysis, we usually transfer the free boundary value problem \eqref{eq:NS} to some system defined in the fixed domain. Without taking surface tension into account, it is standard to write down \eqref{eq:NS} in the Lagrangian coordinates, which allows us to pull back \eqref{eq:NS} to the initial domain along the trajectory of fluid particles (See Figure \ref{fig:LT}). 

\begin{center}
	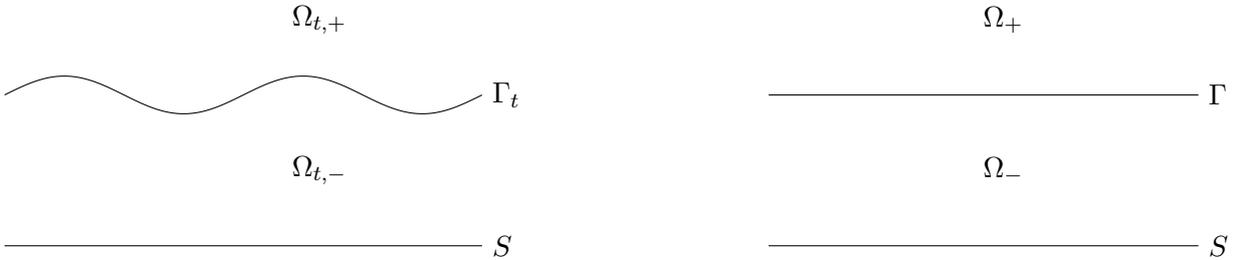
\begin{figure}[h]
		\centering
		\begin{tikzpicture}[>=stealth]
			
			\draw[ domain=-pi:pi, samples=200] plot (\x, {sin(2*\x r)/4}) node[right] {$\Gamma_t$};
			
			\draw (-pi,-2) -- (pi,-2) node[right] {$S$};
			
			\node at (1,1) {$\Omega_{t,+}$};
			\node at (1,-1) {$\Omega_{t,-}$};
			
			\draw (2.2*pi,0) -- (4*pi,0) node[right] {$\Gamma$};
			
			\draw (2.2*pi,-2) -- (4*pi,-2) node[right] {$S$};
			
			\node at (10,1) {$\Omega_{+}$};
			\node at (10,-1) {$\Omega_{-}$};
			
			
		\end{tikzpicture}
		\caption{Role of the Lagrangian transform}
  \label{fig:LT}
	\end{figure}
\end{center}

Now, let us introduce the so-called \emph{Lagrangian transformations} for the two-phase problem \eqref{eq:NS}:
\begin{equation}\label{def:LT}
	X_{\bv_\pm}(y,t)= y + \int^t_0 \bv_\pm(y, s)\,ds,
 \quad \forall \,\,y\in \Omega= \Omega_{+} \cup \Gamma \cup \Omega_- 
 \,\,\,\text{and}\,\,\,0<t<T,
\end{equation}
where \emph{Lagrangian velocities} $\bv_\pm=\bv_\pm(y,t).$ 
Moreover, we suppose that 
\begin{equation*}
\begin{aligned}
   \Omega_{t,\pm} &= \{x =X_{\bv_\pm}(y,t): y \in \Omega_{\pm}\}, \\
	\Gamma_t &= \{x =X_{\bv_+}(y,t){ = X_{\bv_-}(y,t)} : y \in \Gamma\}.
\end{aligned}
\end{equation*}
To use the inverse mappings $X_{\bv_{\pm}}^{-1}(x, t)$ of the transformations in \eqref{def:LT}, we assume that 
\begin{equation}\label{asmp:small_1}
	\int^T_0\|\bv_\pm(\cdot, s)\|_{H^1_\infty(\Omega_\pm)}\,ds \leq \delta <1.
\end{equation}

In what follows, we will derive the nonlinear system in the Lagrangian coordinates.
Let us denote the positive constant $\gamma_{1+}$ for the equilibrium state of the density, namely, 
\begin{equation*}
  \rho_{+,0}=\gamma_{1+}+\eta_{+,0}  
\end{equation*}
for some function $\eta_{+,0.}$
Then we set the new unknowns:
\begin{gather*}
\rho_+(x, t) =  \gamma_{1+} +\eta_+ \big(X^{-1}_{\bv_+}(x, t), t\big), \quad
\bu_{\pm}(x, t) = \bv_{\pm}\big(X^{-1}_{\bv_\pm}(x, t), t\big),\quad 
\fp_-(x,t)=\pi_-\big(X^{-1}_{\bv_-}(x,t), t\big).
\end{gather*}
Then the kinematic condition $V_{\Gamma_t}=\bu_-\cdot\bn_t$ is automatically satisfied under
\eqref{def:LT}. 
\smallbreak 
Note that
$$\nabla_y X_{\bv_\pm}= \bI + \int^t_0\nabla_y\bv_\pm(y, s)\,ds.$$
Then, by the assumption \eqref{asmp:small_1}, the inverses of $\nabla_y X_{\bv_\pm}$ are given by 
\begin{equation}\label{eq:Vk}
\big(\nabla_y X_{\bv_\pm}\big)^{-1}
= \bI + \bV_0(\bk_\pm) ,    
\end{equation}
where $\bV_0(\bk_\pm)$ are matrix-valued functions defined by 
$$\bV_0(\bk_\pm)= \sum_{j=1}^{\infty} (-\bk_\pm)^j
\quad \text{for}\quad  
\bk_\pm =\int^t_0\nabla_y\bv_\pm(y, s)\,ds.$$
In particular, $\bV_0(0) = 0.$
Furthermore, let us set the Jacobians $J_{\bv_{\pm}}=\det (\nabla_y X_{\bv_{\pm}}).$ 
In fact, $J_{\bv_{-}}=1$ due to the divergence free condition of $\bu_-$.

In view of the chain rule, we introduce the gradient, divergence and stress tensor operators with respect to \eqref{def:LT},
\begin{equation}\label{eq:1.5}
	\begin{aligned}
		\nabla_{\bv_{\pm}}\bv_{\pm} &= \big(\bI + \bV_0(\bk_{\pm}) \big)\nabla_y\bv_{\pm},\\
		\dv_{\bv_{\pm}} \bv_{\pm} &= \big(\bI+\bV_0(\bk_{\pm})\big):\nabla_y \bv_{\pm}, \\
  		\DV_{\bv_+} \bA &= J_{\bv_+}^{-1} \DV_{y} \Big(J_{\bv_+}\bA \big(\bI + \bV_0(\bk_+) \big) \Big),\\
 \DV_{\bv_-} \bA& =  \DV_{y} \Big(\bA \big(\bI + \bV_0(\bk_-) \big) \Big),\\
   \bT_{\bv_{\pm}}(\bv_{\pm},\pi_{\pm})&=\bS_{\bv_{\pm}}(\bv_{\pm})-\pi_{\pm}\bI\\
  \bS_{\bv_+}(\bv_+)&=\mu_+\bD_{\bv_+} (\bv_+)
		+(\nu_+-\mu_+) (\dv_{\bv_+} \bv_+) \bI,\\ 
  \bS_{\bv_-}(\bv_-)&=\mu_{-}\bD_{\bv_-}(\bv_-),\\
  \bD_{\bv_{\pm}} (\bv_{\pm}) &=\big(\bI + \bV_0(\bk_{\pm}) \big)\nabla \bv_{\pm}
		+(\nabla \bv_{\pm})^{\top}\big(\bI + \bV_0(\bk_{\pm}) \big)^{\top}.
	\end{aligned}
\end{equation}
In addition, $\DV_{\bv_\pm} \bA$ can be also written via
\begin{equation*}
	\DV_{\bv_\pm} \bA= \DV_y \bA+\big(\bV_0(\bk_\pm) \nabla \,|\, \bA\big)
\end{equation*}
with the $i$th component 
$(\bB \nabla | \bA)_{i}= \sum_{j,k=1}^3 B_{jk}\pd_k A_{ij}$ for $\bB=[B_{ij}]_{3\times 3},$ $i=1,2,3.$
Therefore, like the standard operator $\DV,$ $\DV_{\bv} \bA=0$ if $\bA$ is a constant matrix.
\smallbreak 

According to \eqref{eq:NS} and \eqref{eq:1.5}, the knew unknowns $\eta_+,\bv_\pm$ and $\fp_-$ satisfy the following system:
\begin{equation}\label{eq:LL_CNS_1}
	\left\{ \begin{aligned}
		&\pd_t\eta_+ + (\gamma_{1+}+\eta_+)\, \dv_{\bv_+} \bv_+=0
		&&\quad&\text{in} &\quad \Omega_+ \times (0,T), \\
		&(\gamma_{1+}+\eta_+)\pd_t\bv_+-\DV_{\bv_+}\bT_{\bv_+}\big(\bv_+,P(\gamma_{1+}+\eta_+)\big)  = 0 
		&&\quad&\text{in}& \quad \Omega_+ \times (0,T), \\
		&\rho_-\pd_t\bv_-  - \Di_{\bv_-}\bT_{\bv_-}(\bv_{-},\fp_-) = 0,\quad \di_{\bv_-} \bv_- = 0,
		&&\quad&\text{in}& \quad \Omega_- \times (0,T), \\
		&\Big(\bT_{\bv_+}\big(\bv_+,P(\gamma_{1+}+\eta_+)\big) - \bT_{\bv_-}(\bv_-,\fp_-)\Big)\bn_{\bv}
		= 0
		&&\quad&\text{on}& \quad \Gamma \times (0,T),\\
		&\bv_{+}-\bv_{-}
		= 0
		&&\quad&\text{on}& \quad \Gamma \times (0,T),  \\
		& \bv_-=0
		&&\quad&\text{on}& \quad   S \times (0,T), \\
		&(\eta_+, \bv_+)|_{t=0} = (\eta_{+,0}, \bu_{+,0})
		&&\quad&\text{in}& \quad \Omega_+, \\
		&\bv_-\big|_{t=0}=\bu_{-,0}
		&&\quad&\text{in}& \quad \Omega_-,
	\end{aligned}
	\right.
\end{equation}
where the $\bn_\bv$ is defined by
\begin{equation}\label{eq:nv}
  \bn_{\bv} = { \frac{\big(\bI + \bV_0(\bk_+)\big)\bn }{\big|\big(\bI + \bV_0(\bk_+)\big)\bn \big|} }
= { \frac{\big(\bI + \bV_0(\bk_-)\big)\bn }{\big|\big(\bI + \bV_0(\bk_-)\big)\bn \big|} }
\end{equation}
for $\bn=(0,\dots,0,1)^{\top}.$  In fact, the second equality in \eqref{eq:nv} is due to the following lemma.
\begin{lem}
Let us denote
    \begin{equation*}
        \bn_{\bv_+} =  \frac{\big(\bI + \bV_0(\bk_+)\big)\bn }{\big|\big(\bI + \bV_0(\bk_+)\big)\bn \big|} 
        \quad \text{and}\quad 
        \bn_{\bv_-}=\frac{\big(\bI + \bV_0(\bk_-)\big)\bn }{\big|\big(\bI + \bV_0(\bk_-)\big)\bn \big|} \cdot
    \end{equation*}
    Then we have $\bn_{\bv_+}=\bn_{\bv_-}$ if  $\delta$ in \eqref{asmp:small_1} is small enough.
\end{lem}
\begin{proof}
Let $\mathcal{A}_{\bv_\pm}$ stand for the co-factor matrices of $\nabla_y X_{\bv_\pm}.$ 
Then by \eqref{eq:Vk} we get
\begin{equation*}
       \CA_{\bv_+}= J_{\bv+} \big( \bI+\bV_0(\bk_+) \big),\quad 
       \CA_{\bv_-}= \bI+\bV_0(\bk_-)
\end{equation*}
due to $J_{\bv_-}=0.$ Assuming $\delta$ in \eqref{asmp:small_1} small enough, we have $J_{\bv+}>0.$ 
Thus it suffices to prove 
    \begin{equation}\label{bn2}
        \CA_{\bv_+}\bn=\CA_{\bv_-}\bn \,\,\, \text{on}\,\,\, \Gamma.
    \end{equation}
    
    Note that the tangential derivatives of $X_{\bv_\pm}$ are continuous across the boundary $\Gamma$ since $\bv_{+}=\bv_{-}$ on $\Gamma,$ namely, 
    \begin{equation*}
        \nabla_{y'} X_{\bv_+} = \nabla_{y'} X_{\bv_-} \,\,\,  \text{on}\,\,\, \Gamma
    \end{equation*}
for $y'=(y_1,\dots,y_{N-1}).$ Then we see 
    \begin{equation*}
        \begin{aligned}
           &\text{the i-th component of}\,\,\, \CA_{\bv_+}\bn\\
           =&\, (-1)^{N+i} \det ( \nabla_{y'} X_{\bv_+,1} \wedge  \nabla_{y'} X_{\bv_+,2} 
           \wedge \cdots\wedge\nabla_{y'} X_{\bv_+,i-1} \wedge \nabla_{y'} X_{\bv_+,i+1} 
           \wedge\cdots\wedge \nabla_{y'} X_{\bv_+,N} )\\
          =& \,(-1)^{N+i} \det ( \nabla_{y'} X_{\bv_-,1} \wedge  \nabla_{y'} X_{\bv_-,2} 
           \wedge \cdots\wedge\nabla_{y'} X_{\bv_-,i-1} \wedge \nabla_{y'} X_{\bv_-,i+1} 
           \wedge\cdots\wedge \nabla_{y'} X_{\bv_-,N} )\\
           =&\, \text{the i-th component of}\,\,\, \CA_{\bv_-}\bn.
        \end{aligned}
    \end{equation*}
    Here $(\ba_1' \wedge \cdots \wedge \ba_{N-1}')$ denotes $(N-1)\times (N-1)$ matrix spanned by the column vectors $\ba_{j}'\in \mathbb{R}^{N-1}$ for $j=1,\dots,N-1$. This proves the claim \eqref{bn2}. 
\end{proof}

\subsection{Linearization of \eqref{eq:LL_CNS_1}}
\label{section7.1}
In this subsection, we shall rewrite \eqref{eq:LL_CNS_1} into the form of \eqref{eql:1}. 
To this end, suppose that 
\begin{equation*}
    \rho_{+,0} = \gamma_{1+} + \theta_{+,0}, \quad 
    \gamma_{1-}=\rho_- \quad \text{and}\quad 
    \gamma_{2+}=P'(\gamma_{1+}).
\end{equation*}
Then, by Taylor's theorem, we arrive
\begin{equation*}
	P(\gamma_{1+}+\eta_+)-P(\gamma_{1+})=\gamma_{2+} \eta_++Q(\eta_+)
\end{equation*}
with 
$$Q(\eta_+)=\eta_+^2\int_{0}^{1}P''(\gamma_{1+}+\theta\eta_+)(1-\theta)d\theta.$$
For simplicity, let us introduce 
\begin{equation*}
    \begin{aligned}
    \fq_{-}&=\fp_{-}-P(\gamma_{1+}),\\
        \CD_{\bD_\pm}(\bv_\pm)&=\bD_{\bv_\pm}(\bv_\pm)-\bD(\bv_\pm)=\bV_0(\bk_\pm)\nabla{\bv_\pm}+(\bV_0(\bk_\pm)\nabla{\bv_\pm})^{\top},\\
          \CS_{\bD_+} (\bv_+) & = \mu_+\CD_{\bD_+} (\bv_+) +(\nu_+-\mu_+) \big(\bV_{0}(\bk_+):\nabla\bv_{+}\big) \bI.
    \end{aligned}
\end{equation*}
Then we see from \eqref{eq:LL_CNS_1} that
\begin{equation}\label{eq:NS3}
	\left\{ \begin{aligned}
		&\pd_t\eta_+ + \gamma_{1+}\, \dv \bv_+=f_+(\eta_+,\bv_{+})
		&&\quad&\text{in} &\quad \Omega_+ \times (0,T), \\
		&\gamma_{1+}\pd_t\bv_+-\DV\bT_+(\bv_+,\gamma_{2+}\eta_+)  = \bg_+(\eta_+,\bv_{+})
		&&\quad&\text{in}& \quad \Omega_+ \times (0,T), \\
		&\gamma_{1-}\pd_t\bv_-  - \Di\bT_-(\bv_{-},\fq_-) = \bg_-(\bv_{-},\fq_-)
		&&\quad&\text{in}& \quad \Omega_- \times (0,T), \\
  		& \di \bv_- = g_d (\bv_{-})=\dv \fg_d (\bv_{-})
		&&\quad&\text{in}& \quad \Omega_- \times (0,T), \\
		&\big(\bT_+(\bv_+,\gamma_{2+}\eta_+) - \bT_-(\bv_-,\fq_-)\big)\bn
		= \bh(\eta_+,\bv_{\pm},\fq_-)
		&&\quad&\text{on}& \quad \Gamma \times (0,T),\\
		&\bv_{+}-\bv_{-}
		= 0
		&&\quad&\text{on}& \quad \Gamma \times (0,T),  \\
		& \bv_-=0
		&&\quad&\text{on}& \quad   S \times (0,T), \\
		&(\eta_+, \bv_+)|_{t=0} = (\theta_{+,0}, \bu_{+,0})
		&&\quad&\text{in}& \quad \Omega_+, \\
		&\bv_-\big|_{t=0}=\bu_{-,0}
		&&\quad&\text{in}& \quad \Omega_-,
	\end{aligned}
	\right.
\end{equation}
where we have defined that 
\begin{equation}\label{eq:22}
	\begin{split}
		f_+(\eta_+,\bv_+)&=-\eta_{+}\di\bv_{+}-(\gamma_{1+}+\eta_+) \bV_0(\bk_+):\nabla\bv_{+},\\
		\bg_+(\eta_+,\bv_+)&=-\eta_+\pd_t\bv_{+}+\Di\big(\CS_{\bD_+} (\bv_+)-Q(\eta_+)\,\bI\big)
		 \\&\quad +\Big(\bV_{0}(\bk_+)\nabla \,\big|\, \bS_{\bv_+}(\bv_{+})-\gamma_{2+}\eta_+\bI-Q(\eta_+)\bI\Big),\\
		\bg_-(\bv_-,\fq_-)&=\Di\big(\mu_-\CD_{\bD_-}(\bv_-)\big)+\big(\bV_0(\bk_-)\nabla\,\big|\,\mu_-\bD_{\bv_-}(\bv_-)-\fq_-\bI\big),\\
		g_d (\bv_-)&=-\bV_0(\bk_-):\nabla\bv_-,\qquad  
  \fg_d (\bv_{-}) = -\bV_0(\bk_-)^{\top} \bv_-,
		\\\bh(\eta_+,\bv_\pm, \fq_-)&=\bT_+(\bv_+,\gamma_{2+}\eta_+)\bn-\bT_{\bv_+}(\bv_+,\gamma_{2+}\eta_++Q(\eta_+))\Big(\frac{\big(\bI + \bV_0(\bk_+)\big)\bn }{\big|\big(\bI + \bV_0(\bk_+)\big)\bn \big|}-\bI \Big)\\
  &\quad -\bT_{\bv_+}(\bv_+,\gamma_{2+}\eta_++Q(\eta_+))\bn-\bT_-(\bv_-,\fq_-)\bn\\
  &\quad 
        +\bT_{\bv_-}(\bv_-,\fq_-)\Big(\frac{\big(\bI + \bV_0(\bk_-)\big)\bn }{\big|\big(\bI + \bV_0(\bk_-)\big)\bn \big|}-\bI \Big)+\bT_{\bv_-}(\bv_-,\fq_-)\bn.
	\end{split}
\end{equation}

\section{Some technical results on the multipliers}
\label{proof}
Lemmas \ref{lem:-+}, \ref{lem:+-} and \ref{lem:+-2} can be proved by the following three lemmas respectively. For simplicity, we will provide the detailed proof of Lemma \ref{lem:-+2} and Lemmas \ref{lem:+-1} and \ref{lem:B3} follow analogously.
\begin{lem}
\label{lem:-+2}
	Assume that $0<\varepsilon<\pi/2$, $1<q<\infty$ and the constants 
    $\mu_{+},$ $\nu_+,$ $\gamma_{1+},$ $\gamma_{2+}>0.$ Let the set $\Gamma_{\ep,\lambda_0}$ be defined by \eqref{def:Gamma} for some $\delta_0>0,$ and let the multiplier $m(\lambda,\xi')$ belong to $\bM_{0,2}(\widetilde{\Gamma}_{\varepsilon,\lambda_0}).$ For $0<x_N<+\infty$, we define $K_{-,+}^{j,n}\,(j,n=1,2)$ by
	\begin{equation*}
		\begin{aligned}
			\big[K_{-,+}^{1,n}(\lambda)f\big](x',x_N)=&\int_{-b}^0\CF_{\xi'}^{-1}\Big[
			E_n^-(y_N)m(\lambda,\xi')Ae^{-B_+x_N}e^{Ay_N}\widehat{f}(\xi',y_N)\Big](x')dy_N,\\
			\big[K_{-,+}^{2,n}(\lambda)f\big](x',x_N)=&\int_{-b}^0\CF_{\xi'}^{-1}\Big[
			E_n^-(y_N)m(\lambda,\xi')A^2M_+(x_N)e^{Ay_N}\widehat{f}(\xi',y_N)\Big](x')dy_N.
		\end{aligned}
	\end{equation*}
	Then there exists a constant $r_b$ depending solely on 
	$\varepsilon,$ $q,$ $N,$ $\gamma_{1+},$ $\gamma_{2+},$ $\mu_{+},$ 
	$\nu_+,$  $\lambda_0$, $\delta_0$ and $b$ such that
    \begin{equation}\label{Rbdd:K_-+}
        \mathcal{R}_{\mathcal{L}\big(L_q(\Omega_{-}),L_q(\Omega_{+})\big)}
	\Big(\Big\{(\tau\pd_{\tau})^\ell K_{-,+}^{j,n}(\lambda):
	\lambda\in \Gamma_{\varepsilon,\lambda_0} \Big\}\Big)\leq r_b
    \end{equation}
	for any $\ell=0,1$ and $j,n=1,2.$
\end{lem}

\begin{lem}
\label{lem:+-1}
	Assume that $0<\varepsilon<\pi/2$, $1<q<\infty$ and the constants 
     $\nu_+,$ $\mu_-,$ $\gamma_{1\pm},$ $\gamma_{2+},$ $\lambda_0,$ $\delta_0$ and $b>0.$ 
     Let the set $\Gamma_{\ep,\lambda_0}$ be defined by \eqref{def:Gamma} and let the multiplier $m(\lambda,\xi')$ belong to $\bM_{0,2}(\widetilde{\Gamma}_{\varepsilon,\lambda_0}).$ For $-b<x_N<0$, we define $K_{+,-}^{j,n}\,(j=1,2,3,4,\,n=1,2)$ by
	\begin{equation*}
		\begin{aligned}
			\big[K^{1,n}_{+,-}(\lambda)f\big](x',x_N)=&\int_{0}^b\mathcal{F}_{\xi'}^{-1}\Big[
			E_n^+(y_N)m(\lambda,\xi')A^2e^{-Ay_N}M_-(x_N)\mathcal{F}_{y'}[f](\xi',y_N)\Big](x')dy_N,\\
			\big[K^{2,n}_{+,-}(\lambda)f\big](x',x_N)=&\int_{0}^b\mathcal{F}_{\xi'}^{-1}\Big[
			E_n^+(y_N)m(\lambda,\xi')Ae^{-Ay_N}e^{B_-x_N}\mathcal{F}_{y'}[f](\xi',y_N)\Big](x')dy_N,\\
			\big[K^{3,n}_{+,-}(\lambda)f\big](x',x_N)=&\int_{0}^b\mathcal{F}_{\xi'}^{-1}\Big[
			E_n^+(y_N)m(\lambda,\xi')A^2e^{-Ay_N}M_-(-x_N-b)\mathcal{F}_{y'}[f](\xi',y_N)\Big](x')dy_N,\\
			\big[K^{4,n}_{+,-}(\lambda)f\big](x',x_N)=&\int_{0}^b\mathcal{F}_{\xi'}^{-1}\Big[
			E_n^+(y_N)m(\lambda,\xi')Ae^{-Ay_N}e^{B_-(-x_N-b)}\mathcal{F}_{y'}[f](\xi',y_N)\Big](x')dy_N.
		\end{aligned}
	\end{equation*}
	Then there exists a constant $r_b$ depending solely on 
	$\varepsilon,$ $q,$ $N,$  $\nu_+,$ $\mu_-,$ $\gamma_{1\pm},$ $\gamma_{2+},$ $\lambda_0,$ $\delta_0$ and $b$ 
	  such that
     \begin{equation}\label{eq:R2}
         \mathcal{R}_{\mathcal{L}\big(L_q(\Omega_{+b}); L_q(\Omega_-)\big)}
	\Big(\Big\{(\tau\pd_{\tau})^\ell  K^{j,n}_{+,-}(\lambda):
	\lambda\in \Gamma_{\varepsilon,\lambda_0} \Big\}\Big)\leq r_b
     \end{equation}
	for $\ell=0,1$ and $j=1,2,3,4$, $n=1,2$.
\end{lem}
\begin{lem}\label{lem:B3}
	Assume that $0<\varepsilon<\pi/2$, $1<q<\infty$ and the constants $\mu_{-},$ $\gamma_{1-}>0.$ Let the set $\Gamma_{\ep,\lambda_0}$ be defined by \eqref{def:Gamma} for some $\delta_0>0,$ and let the multiplier $m(\lambda,\xi')$ belong to $\bM_{1,2}(\widetilde{\Gamma}_{\varepsilon,\lambda_0}).$ For $-b<x_N<0$, we define $K_{+,-}^{j,n}\,(j=5,6,\,n=1,2)$ by
	\begin{equation*}
		\begin{aligned}
			\big[K^{5,n}_{+,-}(\lambda)f\big](x',x_N)=&\int_{0}^b\mathcal{F}_{\xi'}^{-1}\Big[
			E_n^+(y_N)m(\lambda,\xi')e^{-Ay_N}e^{Ax_N}\mathcal{F}_{y'}[f](\xi',y_N)\Big](x')dy_N,\\
			\big[K^{6,n}_{+,-}(\lambda)f\big](x',x_N)=&\int_{0}^b\mathcal{F}_{\xi'}^{-1}\Big[
			E_n^+(y_N)m(\lambda,\xi')e^{-Ay_N}e^{A(-x_N-b)}\mathcal{F}_{y'}[f](\xi',y_N)\Big](x')dy_N.
		\end{aligned}
	\end{equation*}
	Then there exists a constant $r_b$ depending solely on 
	$\varepsilon,$ $q,$ $N,$ $\gamma_{1-},$ $\mu_{-},$ 
  $\lambda_0$, $\delta_0$ and $b$ such that
	$$\mathcal{R}_{\mathcal{L}\big(L_q(\Omega_{+b}),L_q(\Omega_-)\big)}
	\Big(\Big\{(\tau\pd_{\tau})^\ell  K^{j,n}_{+,-}(\lambda):
	\lambda\in \Gamma_{\varepsilon,\lambda_0} \Big\}\Big)\leq r_b$$
	for $\ell=0,1$.
\end{lem}

\begin{proof}[Proof of Lemma \ref{lem:-+2}]
For any $\lambda\in\Gamma_{\ep, \lambda_0},$ 
$x=(x',x_N) \in \Omega_+,$ $y=(y',y_N) \in \Omega_-$ and $j,n=1,2,$ let us write 
\begin{equation*}
\big[K_{-,+}^{j,n}(\lambda)f\big](x',x_N)
=\int_{\Omega_-}E_n^{-}(y_N) k_{-,+}^{j}(\lambda,x',y';x_N,y_N)f(y)dy
\end{equation*}
with 
\begin{equation*}
\begin{aligned}
k_{-,+}^{1}(\lambda, x',y';x_N,y_N)&=\CF_{\xi'}^{-1}\big[m(\lambda,\xi')Ae^{-B_+x_N}e^{Ay_N}\big](x'-y'),\\
k_{-,+}^{2}(\lambda,x',y';x_N,y_N)&=\CF_{\xi'}^{-1}\big[m(\lambda,\xi')A^2M_+(x_N)e^{Ay_N}\big](x'-y').
\end{aligned}
\end{equation*}

{\bf Step 1.} 
We first prove that there exists some positive constant $C_1$ depending on $N,\varepsilon,\lambda_0,\mu_{+},\nu_+$ such that
\begin{equation}\label{eq:est:k1}
\big|(\tau\pd_{\tau})^{\ell}k_{-,+}^{j,n}(\lambda,x',y';x_N,y_N)\big|
\leq C_1 \big(|x'-y'|+x_N-y_N\big)^{-N}
\end{equation}
for any $\lambda\in\Gamma_{\ep, \lambda_0}$, $\ell=0,1,$ $j,n=1,2,$ and $0<x_N<\infty$.
\smallbreak 

In fact, by the assumption on $m(\lambda,\xi'),$ Leibniz's rule and Lemma \ref{lem:4.3},  there exist positive constants 
$C=C_{\alpha',\varepsilon,\lambda_0,\mu_{+},\gamma_{1+}}$ and $d=d_{\varepsilon,\mu_{+},\nu_+,\gamma_{1+}}$ such that 
\begin{equation}\label{eq:42}
\begin{aligned}
\big|\pd^{\alpha'}_{\xi'}(\tau\pd_{\tau})^{\ell} \big(m(\lambda,\xi')Ae^{-B_+x_N}e^{Ay_N}\big)\big|
&\leq  CA^{1-|\alpha'|}
e^{-d(|\lambda|+A)x_N}
e^{A y_N/2}\\
& \leq C A^{1-|\alpha'|} e^{-c A(x_N-y_N)},\\
\big|\pd^{\alpha'}_{\xi'}(\tau\pd_{\tau})^{\ell} \big(m(\lambda,\xi')A^2M_+(x_N)e^{Ay_N}\big)\big|
&\leq  C A^{2-|\alpha'|}
x_N e^{-d A x_N}
e^{A y_N/2}\\
& \leq  C A^{1-|\alpha'|}
e^{-c A(x_N-y_N)}
\end{aligned}
\end{equation}
for $c=\min \{d,1\}/2$ and for any multi-index $\alpha'\in\BN_0^{N-1}.$
Then we apply \cite[Theorem 2.3]{shishi2001} 
\footnote{One may take the parameters $n=N-1$ and $\alpha=\sigma=1$ in \cite[Theorem 2.3]{shishi2001}.}
to obtain
\begin{equation}\label{eq:43}
	\big|(\tau\pd_{\tau})^{\ell}k_{-,+}^{j,n}(\lambda,x',y';x_N,y_N)\big|
    \lesssim |x'-y'|^{-N}.
\end{equation}

On the other hand, using \eqref{eq:42} with $\alpha'=0$ and the change of variables $\eta'=c(x_N-y_N)\xi'$ we have
\begin{equation*}
	\begin{aligned}
		\big|(\tau\pd_{\tau})^{\ell}k_{-,+}^{j,n}(\lambda,x',y';x_N,y_N)\big|
        &\lesssim \int_{\mathbb{R}^{N-1}}Ae^{-cA(x_N-y_N)}d\xi'\\
		&\lesssim  \big(c (x_N-y_N)\big)^{-N}\int_{\mathbb{R}^{N-1}}|\eta'|e^{-|\eta'|}d\eta'\\
		& \lesssim (x_N-y_N)^{-N},
	\end{aligned}
\end{equation*}
which combined with \eqref{eq:43} implies \eqref{eq:est:k1}. 
\medskip

{\bf Step 2.} 
In views of the property \eqref{eq:est:k1}, we introduce the operator
\begin{equation*}
	[K_0^nf](x',x_N)=\int_{\Omega_-} \big| E_n^{-}(y_N)\big| \,k_0(x'-y',x_N-y_N)\,|f(y)|\,dy.
\end{equation*}
with $k_0(x)=C_{1}|x|^{-N}.$ In what follows, we prove $K_0^n\in \CL\big(L_q(\Omega_{-});L_q(\Omega_{+})\big)$ for $n=1,2$. 
\smallbreak

By Young's inequality, we have
\begin{equation}\label{es:K_0}
	\begin{aligned}
		\|[K_0^nf](\cdot,x_N)\|_{L_q(\BR^{N-1})}
        &\leq \int_{-b}^0\big|E_n^{-}(y_N)\big|\,\|k_0(\cdot,x_N-y_N)\|_{L_1(\BR^{N-1})}
        \|f(\cdot,y_N)\|_{L_q(\BR^{N-1})}dy_N\\
		&\lesssim  \int_{-b}^0 \big|E_n^{-}(y_N)\big|\,\frac{\|f(\cdot,y_N)\|_{L_q(\BR^{N-1})}}{x_N-y_N}dy_N.
	\end{aligned}
\end{equation}

\begin{itemize}
    \item Case $n=1$. Using the change of variable $y_N=-x_Nt$, we get from \eqref{es:K_0} that
\begin{equation}\label{eq:45}
	\begin{aligned}
		\|[K_0^1f](\cdot,x_N)\|_{L_q(\BR^{N-1})}
		&\lesssim \int_{0}^{b/x_N}\varphi_{-0}(-x_Nt)\frac{\|f(\cdot,-x_Nt)\|_{L_q(\BR^{N-1})}}{1+t}dt,
	\end{aligned}
\end{equation}
which combined with Minkowski's inequality implies that 
\begin{equation*}
	\begin{aligned}
		\|K_0^1f\|_{L_q(\Omega_+)}
        &\lesssim\int_{0}^\infty\frac{1}{1+t}
        \Big[\int_{0}^{b/t}\big(\varphi_{-0}(-x_Nt)
        \|f(\cdot,-x_Nt)\|_{L_q(\BR^{N-1})}\big)^qdx_N\Big]^{1/q}dt\\
		&\lesssim\int_{0}^\infty\frac{1}{1+t} \Big[\int_{-b}^0
        \big(\varphi_{-0}(y_N)\|f(\cdot,y_N)\|_{L_q(\BR^{N-1})}\big)^q\frac{dy_N}{t}\Big]^{1/q}dt\\
		&\lesssim \|f\|_{L_q(\Omega_{-})}\int_{0}^\infty\frac{1}{(1+t)t^{1/q}}dt\\
        &\lesssim \|f\|_{L_q(\Omega_{-})}.
	\end{aligned}
\end{equation*}

\item Case $n=2.$ Analogous to the \eqref{eq:45}, we obtain from the change of variable $y_N=-x_Nt$ and \eqref{es:K_0} that
\begin{equation}\label{eq:46}
	\begin{aligned}
		\|[K_0^2f](\cdot,x_N)\|_{L_q(\BR^{N-1})}
        &\lesssim \int_{-2b/3}^{-b/3}\big| \varphi'_{-0}(y_N)\big| \frac{\|f(\cdot,y_N)\|_{L_q(\BR^{N-1})}}{x_N-y_N}dy_N\\
		&\lesssim \int_{b/(3x_N)}^{2b/(3x_N)}
        \big| \varphi'_{-0}(-x_Nt)\big|\frac{\|f(\cdot,-x_Nt)\|_{L_q(\BR^{N-1})}}{1+t}dt,
	\end{aligned}
\end{equation}
which yields 
\begin{equation*}
	\begin{aligned}
		\|K_0^2f\|_{L_q(\Omega_+)}
        &\lesssim \int_{0}^\infty\frac{1}{1+t}\Big[\int_{b/(3t)}^{2b/(3t)}
        \big( \big| \varphi'_{-0}(-x_Nt) \big| \,\|f(\cdot,-x_Nt)\|_{L_q(\BR^{N-1})}\big)^q \,dx_N\Big]^{1/q}dt\\
		&\lesssim \int_{0}^\infty\frac{1}{1+t}\Big[\int_{-2b/3}^{-b/3} \big(\big| \varphi'_{-0}(y_N)\big|\, 
        \|f(\cdot,y_N)\|_{L_q(\BR^{N-1})}\big)^q\frac{dy_N}{t}\Big]^{1/q}dt\\
		&\lesssim \|\varphi'_{-0}\|_{L_\infty(\BR)}\|f\|_{L_q(\Omega_{-})}\int_{0}^\infty\frac{1}{(1+t)t^{1/q}}dt\\
         &\lesssim \|f\|_{L_q(\Omega_{-})}.
	\end{aligned}
\end{equation*}
\end{itemize}
\medskip

{\bf Step 3.} Now, we prove \eqref{Rbdd:K_-+} by verifying Definition \ref{def:R}. 
Assume that  $M\in\BN, \lambda_k\in\Gamma_{\ep,\lambda_0}, f_k\in L_q(\Omega_{-}),$ $r_j(u)$ be the Rademacher functions  and $k=1,\dots,M.$
By Khintchine inequality (see \cite[(3.4)]{Denk2003} for instance), we obtain
\begin{equation*}
	\begin{aligned}
	&\Big(\int_0^1 \Big\|\sum_{k=1}^M r_k(u)(\tau\pd_{\tau})^{\ell} K_{-,+}^{j,n}(\lambda_k)f_k\Big\|_{L_q(\Omega_{+})} \Big)^{1/q}\\
    \lesssim  & \Big\|\Big(\sum_{k=1}^M |(\tau\pd_{\tau})^{\ell} 
    K_{-,+}^{j,n}(\lambda_k)f_k|^2\Big)^{1/2} \Big\|_{L_q(\Omega_{+})}\\
    \lesssim & \Big\|\big(\sum_{k=1}^M (K_0^{n} f_k )^2\big)^{1/2} \Big\|_{L_q(\Omega_{+})}\\
    \lesssim & \Big\| K_0^n \Big(\sum_{k=1}^M |f_k|^2\Big)^{1/2}\Big\|_{L_q(\Omega_{+})}\\
\lesssim  &\|K_0^n\|_{\CL\big(L_q(\Omega_{-}); L_q(\Omega_{+})\big)}
\Big\|\Big(\sum_{k=1}^M |f_k|^2\Big)^{1/2}\Big\|_{L_q(\Omega_{-})}\\
\lesssim & \|K_0^n\|_{\CL\big(L_q(\Omega_{-}),L_q(\Omega_{+})\big)}
\Big(\int_0^1 \Big\|\sum_{k=1}^M r_k(u)f_k\Big\|_{L_q(\Omega_{-})}^q \,du \Big)^{1/q}.
	\end{aligned}
\end{equation*}
This gives the bound \eqref{Rbdd:K_-+}. 
\end{proof}

\section*{Acknowledgement}
X.Z. was partially supported by NSF of China under Grant 12101457 and the Fundamental Research Funds for the Central Universities.


\end{document}